\pdfoutput=1
\documentclass[11pt]{article_ET_m}
\usepackage{lmodern}
\usepackage[utf8]{inputenc}
\usepackage[T1]{fontenc}
\usepackage{amsmath, amsfonts, amssymb}
\usepackage{subcaption, mwe, bbm, stmaryrd}
\usepackage{graphicx}
\usepackage{longtable}
\usepackage{array}
\graphicspath{{plots/}} 
\usepackage{geometry}
\usepackage{amsthm}
\usepackage{etoolbox}
\usepackage{changepage}
\usepackage{enumitem}
\usepackage{dsfont}
\usepackage{tikz}
\usepackage{pgfplots}
\pgfplotsset{compat=1.18}
\usetikzlibrary{calc}
\usetikzlibrary{automata, positioning, arrows.meta}

\usepackage{caption}

\usepackage{mathtools}

\definecolor{db}{RGB}{0, 0, 130}
\usepackage[colorlinks=true,citecolor=red,linkcolor=db,urlcolor=blue,pdfstartview=FitH]{hyperref}
\usepackage{cleveref}

\makeatletter
\newtheoremstyle{nonitalic}
  {\topsep}   
  {\topsep}   
  {\normalfont}  
  {}          
  {\bfseries} 
  {.}         
  {.5em}      
  {}          
\makeatother
\theoremstyle{plain}

\newtheorem{theoreme}{Theorem}[section]
\newtheorem{proposition}[theoreme]{Proposition}
\newtheorem{lemma}[theoreme]{Lemma}
\newtheorem{corollaire}[theoreme]{Corollary}
\theoremstyle{definition}

\newtheorem{remarque}[theoreme]{Remark}

\newcommand{\CR}{\mathrm{CR}}
\newcommand{\dd}{\mathrm{d}}

\newcommand{\LC}{\mathcal{L}_{\mathfrak{C}}^{\theta, \tau}}

\newcommand{\LCtc}{\mathcal{L}_{\mathfrak{C}}^{\theta, \tau, \circ}}
\newcommand{\Lp}{\overline{\Lambda}{}^{\theta, \tau}}
\newcommand{\Lm}{\underline{\Lambda}{}^{\theta, \tau}}

\newcommand{\bpi}{\boldsymbol{\pi}}
\newcommand{\LCl}{\mathcal{L}_\mathfrak{C}^{\theta, \le \varepsilon}}
\newcommand{\LCg}{\mathcal{L}_\mathfrak{C}^{\theta, > \varepsilon}}

\makeatletter
\newcommand{\@footnotenomark}[1]{%
  \insert\footins{%
    \reset@font\footnotesize
    \interlinepenalty\interfootnotelinepenalty
    \splittopskip\footnotesep
    \splitmaxdepth\dp\strutbox \floatingpenalty\@MM
    \hsize\columnwidth \@parboxrestore
    \noindent\rule\z@\footnotesep\ignorespaces#1\@finalstrut\strutbox}%
}
\newcommand{\subjclass}[1]{\@footnotenomark{\textit{MSC2020 subject classifications:} #1}}
\newcommand{\keywords}[1]{\@footnotenomark{\textit{Keywords:} #1}}
\makeatother

\author{Julien Enguehard\footnote{Universit\'e Paris Cité, France. \texttt{enguehard@lpsm.paris}.} 
}

\title{Sharp crossing exponents, truncated percolation, and boundary excursions in the Brownian loop soup}

\begin{document}
\date{}

\maketitle
\begin{abstract} 
Motivated by the recent work of Jego, Lupu, and Qian, we first sharpen the crossing probability for clusters of the Brownian loop soup, establishing the order $|\log r|^{-1+\theta}$. We then study the truncated Brownian loop soup on the semi-infinite cylinder with loops of duration below a fixed threshold $\tau$ removed. We confirm the existence of a phase transition at the intensity $\theta = 1$. For a longitudinal distance $x$, we compute the asymptotic crossing probabilities, confirming a polynomial decay $x^{-1+\theta}$ in the subcritical regime $\theta \in (0, 1)$. At criticality ($\theta = 1$), we establish a logarithmic decay $(\log x)^{-1}$, proving the absence of percolation and thus the continuity of the phase transition. In the super-critical regime $\theta > 1$, we prove that there are almost surely only finitely many bounded clusters. 

In an independent and self-contained part, we study the excursions from the outer boundary of a loop-soup cluster. We derive an explicit formula for their first-order intensity measure, and reveal a surprising non-monotonicity with respect to the soup's intensity $\theta$.

Finally, our sharper crossing estimate allows us to refine several consequences: we sharpen the diagonal singularity of the two-point covariance function for the conformally invariant fields $h_\theta$ introduced in \cite{jego2023conformallyinvariantfieldsbrownian}, we sharpen the Minkowski gauge of loop soup clusters, and give a necessary condition for their polar sets in terms of logarithmic capacity.
\end{abstract}


\section{Introduction}

Given a planar domain $D \subset \mathbb{C}$, the Brownian loop soup $\mathcal{L}_D^\theta$ with intensity $\theta > 0$ is a Poissonian ensemble of Brownian loops in $D$ with intensity measure $\theta \mu_D^{\mathrm{loop}}$. The underlying reference measure, called the Brownian loop measure $\mu_D^{\mathrm{loop}}$, is defined as,
\begin{equation} \label{eq:def_mu_D_loop}
\mu_D^{\mathrm{loop}}(\mathrm{d}\wp) = \int_D \mathrm{d}z \int_0^\infty \frac{\mathrm{d}t}{t} p_D(t, z, z) \mathbb{P}_{D}^{t, z, z}(\mathrm{d}\wp),
\end{equation}
where $p_D$ is the heat kernel in $D$ and $\mathbb{P}_{D}^{t, z, z}$ is the law of a Brownian bridge in $D$ of duration $t$ starting and ending at $z$. Since its introduction by Lawler and Werner \cite{lawler2003brownianloopsoup}, the Brownian loop soup has emerged as a central object in two-dimensional conformal geometry. It exhibits deep connections with the Gaussian free field (GFF) \cite{jan2010markovpathsloopsfields, Lupu_2016, Qian_2019,sheffield2006gaussianfreefieldsmathematicians} and provides a fundamental geometric construction for Conformal Loop Ensembles (CLE) and Schramm-Loewner Evolutions (SLE) \cite{miller2016imaginarygeometryiinteracting, schramm2010contourlinecontinuumgaussian, sheffield2011conformalloopensemblesmarkovian, sheffield2006explorationtreesconformalloop}. This framework is also intimately linked to the First Passage Sets (FPS) of the GFF \cite{Aru_2019, Aru_2020}. Furthermore, the scope of the loop soup continues to expand; recently, it has been used to construct a new class of conformally invariant fields denoted $h_\theta$ \cite{jego2023conformallyinvariantfieldsbrownian} that generalize the GFF for $\theta \in (0,1/2)$.

A key feature of the Brownian loop soup is its percolative behavior, which is known to undergo a phase transition at the critical intensity $\theta = 1/2$ \cite{sheffield2011conformalloopensemblesmarkovian}. Note that much of the literature (e.g. \cite{lawler2003brownianloopsoup}) uses the convention $c = 2\theta$, where $c$ corresponds to the central charge in conformal field theory. In the subcritical regime $\theta \in (0, 1/2]$, the loops almost surely do not percolate. They merge into infinitely many disjoint connected components, called \emph{clusters}, and the outer boundaries of the outermost clusters form a $\mathrm{CLE}_\kappa$ \cite{sheffield2011conformalloopensemblesmarkovian}, where $\kappa \in (8/3, 4]$ is related to the intensity via
\begin{equation}
2\theta \,=\, c \,=\, (3\kappa - 8)(6 - \kappa)/(2\kappa).
\end{equation}

Recently, Jego, Lupu, and Qian \cite{Jego_2025} investigated the crossing exponent for clusters of the Brownian loop soup in the subcritical regime $\theta \in (0, 1/2]$. They showed that the probability for a macroscopic cluster to reach a microscopic disc exhibits a logarithmic decay: for any fixed macroscopic radius $\rho \in (0, 1)$,
\begin{equation} \label{eq:crossing_asymptotics_jlq}
\mathbb{P}\Big( \rho \partial \mathbb{D} \stackrel{\mathcal{L}_{\mathbb{D}}^\theta}{\longleftrightarrow} r \partial \mathbb{D} \Big) = |\log r|^{-1+\theta+o(1)} \quad \text{as } r \to 0.
\end{equation}
However, their method leaves an undetermined $o(1)$ error in the exponent. They conjectured in \cite{jego2023conformallyinvariantfieldsbrownian} that this asymptotic behavior is in fact sharp, meaning that the crossing probability in Equation~\eqref{eq:crossing_asymptotics_jlq} behaves as $c |\log r|^{-1+\theta}$ as $r \to 0$ for some positive constant $c = c(\theta, \rho) $.

Surprisingly, their method to establish \eqref{eq:crossing_asymptotics_jlq} revealed that large crossings of clusters penetrating deep towards the origin are driven exclusively by large loops, essentially ignoring microscopic loops with very short durations. To formalise this, they considered $\mathcal{L}_\mathfrak{C}^\theta := f(\mathcal{L}_{\mathbb{D}}^\theta)$ the Brownian loop soup on the semi-infinite cylinder $\mathfrak{C} = (0,\infty) \times \mathbb{S}^1$, which is conformally equivalent to the punctured unit disc $\mathbb{D} \setminus \{0\}$ via the map $f : re^{i\theta} \mapsto (|\log r|, e^{i\theta})$. To avoid the infinite density of small loops and to focus on the connectivity targeting a point (or, equivalently, going towards infinity on the cylinder), they introduced the \emph{truncated loop soup}: for a fixed time threshold $\tau > 0$, they denoted $\LC$ as the subset of $\mathcal{L}_{\mathfrak{C}}^{\theta}$ consisting only of loops with duration $T(\wp) \ge \tau$.

Apart from breaking conformal invariance, removing these small loops drastically changes the percolation behavior. Because the density of loops is no longer infinite, the critical intensity $\theta = 1/2$ is no longer the threshold for percolation. Instead, the authors in \cite{Jego_2025} revealed the existence of a second critical intensity at $\theta = 1$. They established that in the subcritical regime $\theta \in (0, 1)$, the probability of a cluster crossing a longitudinal distance $x$ along the cylinder decays polynomially as $x^{-1+\theta+o(1)}$. This proves the absence of an unbounded cluster for $\theta < 1$.

Their method has naturally opened up new directions and led to several conjectures.
First, for $\theta \in (1/2, 1)$, their proof required the threshold $\tau$ to be larger than a specific value $\tau_\theta$. They conjectured that the phase transition at $\theta = 1$ should hold for any threshold $\tau > 0$, no matter how small. 
Second, the exact behavior of the crossing probability at the critical point $\theta = 1$ was explicitly left as an open question, as the connection with the Bessel process degenerates and does not allow one to easily rule out percolation.
Finally, in the supercritical regime $\theta > 1$, the authors conjectured that, almost surely, the loops form an unbounded connected component regardless of the threshold $\tau$, and that this giant unbounded cluster should coexist with \emph{infinitely many} bounded clusters. 

The present paper aims to provide answers to these questions and conjectures. First, we manage to remove the $o(1)$ error in the non-truncated loop soup crossing probability. Secondly, when the loops of small duration are removed, we sharpen the crossing probabilities, we establish the exact critical behavior at $\theta = 1$, we extend their results to any arbitrary threshold $\tau > 0$, and we provide a surprising negative answer to the conjecture in the supercritical phase. 

Finally, in an independent and self-contained part, we study the excursions from the boundary of a loop soup cluster \cite{Qian_2019,qian2018lawpointprocessbrownian, qian2018conditioningbrownianloopsoupcluster}, compute via a mass-conservation argument their exact first-order intensity measure and reveal a surprising non-monotonic behavior.
\newline

\subsection{Main results.}

The main body of this paper (Theorems \ref{th:crossing_exponent}–\ref{th:n_delta}) focuses on percolation and crossing probabilities for the Brownian loop soup and its truncated version on the cylinder. Moreover, we also include an independent result (Theorem~\ref{th:intensity}) on the excursions from the boundary of loop soup clusters. It is independent of the rest of the paper, but reveals a surprising non-monotonicity property and has a short, self-contained proof.

\paragraph{Percolation and crossing probabilities.}
Our first result sharpens the crossing exponent in the non-truncated Brownian loop soup established in \cite{Jego_2025}. Let us first formalize the connectivity event: for any loop configuration $\mathcal{L}$ and subsets $A, B \subset \mathfrak{C}$, we write $\{A \stackrel{\mathcal{L}}{\longleftrightarrow} B\}$ to denote the event that a connected cluster of loops in $\mathcal{L}$ intersects both $A$ and $B$. Throughout this paper, $f \asymp g$ will mean that there exist $c_1, c_2 > 0$ such that $c_1 g(x) \leq f(x) \leq c_2g(x)$. Depending on the context, these constants may depend on the intensity $\theta$, the cutoffs $\tau, \varepsilon$, on the macroscopic radius $\rho$, and on the macroscopic distance $y$, but they will always be independent of the microscopic radius $r$ or the longitudinal distance $x$.

\begin{theoreme} \label{th:crossing_exponent} For all $\theta \in (0,1/2]$ and all fixed $\rho \in (0,1)$, we have, as $r \to 0$,
$$
\mathbb{P} \Big( \rho\partial\mathbb{D} \stackrel{\mathcal{L}_\mathbb{D}^\theta}{\longleftrightarrow} r\partial \mathbb{D} \Big) \,\asymp\, |\log r|^{-1+\theta}.
$$
\end{theoreme}
Although we are only able to prove that the crossing probability is of order $|\log r|^{-1+\theta}$ rather than deriving an exact asymptotic, we expect, as conjectured in \cite{jego2023conformallyinvariantfieldsbrownian} that such a precise asymptotic formula holds. 
\medskip

Our second result describes the macroscopic connectivity of the truncated Brownian loop soup on the cylinder $\mathfrak{C} = (0,\infty) \times \mathbb{S}^1$.
We consider the truncated Brownian loop soup $\LC$, obtained by keeping only the loops of $\mathcal{L}^\theta_{\mathfrak{C}}$ with a duration $T(\wp) \ge \tau$, for a fixed threshold $\tau > 0$. While \cite{Jego_2025} identified $\theta = 1$ as the critical intensity for this model, their bounds on the subcritical crossing probabilities contained $o(1)$ error terms in the exponent and, in the regime $\theta \in (1/2, 1)$, required the threshold $\tau$ to be sufficiently large. Furthermore, the percolation behavior exactly at the critical point $\theta = 1$ was left as an open problem. Our second theorem settles these questions. We obtain sharper asymptotic crossing probabilities for all $\theta \le 1$, confirming in passing that the phase transition holds for any arbitrary threshold $\tau > 0$. Most importantly, we completely resolve the critical case $\theta = 1$.

\begin{theoreme} \label{th:phase_transition_loopsoup}  
Let $\mathcal{L}_{\mathfrak{C}}^{\theta}$ be a Brownian loop soup on the cylinder $\mathfrak{C} = (0, \infty) \times \mathbb{S}^1$. Let $\tau > 0$ be any threshold. Then, for any fixed $y > 0$, as $x \to \infty$, we have,
\begin{itemize} [itemsep=0pt, topsep=3pt, parsep=0pt]
\item if $\theta \in (0, 1)$, each cluster formed by the remaining loops is a.s. bounded and, 
$$
\mathbb{P} \Big( \{y\} \times \mathbb{S}^1  \stackrel{\LC}{\longleftrightarrow} \{x\} \times \mathbb{S}^1 \Big) \,\asymp\, x^{-1 + \theta}.
$$
\item if $\theta = 1$, each cluster formed by the remaining loops is a.s. bounded and,
$$
\mathbb{P} \Big( \{y\} \times \mathbb{S}^1  \stackrel{\LC}{\longleftrightarrow} \{x\} \times \mathbb{S}^1 \Big) \,\asymp\, (\log x)^{-1}.
$$
\item if $\theta > 1$, the remaining loops a.s. form a unique unbounded cluster.
\end{itemize}
\end{theoreme}
Theorem \ref{th:phase_transition_loopsoup} completes the picture for $\theta \le 1$. By establishing the exact logarithmic decay at criticality, it rigorously proves the absence of macroscopic percolation at $\theta = 1$, thereby demonstrating that this second phase transition is continuous. 
\medskip

Our third result concerns the super-critical regime $\theta > 1$. It is known that for $\theta > 1$, the separating loops (those disconnecting the two ends of the cylinder, see Section~\ref{sec:prelim_loops}) form a unique unbounded cluster. A very natural conjecture, raised by the authors in \cite[Section 1.3]{Jego_2025}, suggested that despite the presence of this giant unbounded cluster, the truncated loop soup would still contain \textit{infinitely many} bounded clusters. Our next theorem provides a surprising negative answer to this conjecture.

\begin{theoreme} \label{th:n_delta} 
For all $\theta > 1$ and any threshold $\tau > 0$, almost surely, the truncated loop soup $\LC$ contains only finitely many bounded connected components.
\end{theoreme}

At first glance, this result might seem counter-intuitive, which explains why the opposite was conjectured. Indeed, even for $\theta > 1$, the truncated loop soup $\LC$ only contains a finite number of loops in any compact slice of the cylinder, say $[n, n+1) \times \mathbb{S}^1$. Therefore, one might naturally expect to find an infinite sequence of unvisited microscopic holes (i.e., regions of the cylinder visited by no loop) escaping to infinity along the cylinder. Since a loop can have an arbitrarily small spatial diameter while still satisfying the duration constraint $T(\wp) \ge \tau$, one could imagine an infinite sequence of such small, isolated non-separating loops hiding in these unvisited holes, disconnected from the giant unbounded cluster. 

However, the loop soup on the semi-infinite cylinder lacks ergodicity. We prove that the unbounded cluster formed by the separating loops is \emph{asymptotically dense} in the following sense: almost surely, for all $\delta > 0$ the maximum number of unvisited disjoint $\delta$-discs on the entire infinite cylinder is finite (Proposition \ref{prop:N_log_finite}). As a consequence, far enough along the cylinder, there is simply no geometric space left for a small non-separating loop to exist without intersecting the giant separating cluster, and thus being absorbed by it.

\paragraph{Excursions from the boundary of a cluster.}
Consider the outermost cluster surrounding some point in a Brownian loop soup in the unit disc of intensity $\theta(\kappa)$ for some $\kappa \in (8/3,4]$, and let $O$ be the domain delimited by its outer boundary $\partial O$. The loops that touch $\partial O$ naturally decompose into a countable collection $\mathcal{E}_{\partial O}$ of excursions in $O$ with endpoints on $\partial O$. By conformal invariance \cite{Qian_2019}, mapping $O$ onto some simply connected domain $D$ via any conformal map $\psi$ yields a canonical collection of excursions $\mathcal{E}_{\partial D}^\kappa \,:=\: \psi(\mathcal{E}_{\partial O})$ whose law is independent of $\partial O$ and whose endpoints are attached to $\partial D$ (see Figure~\ref{fig:boundary_excursions_map}). However, only in the critical case $\kappa = 4$ (corresponding to $\theta = 1/2$) is $\mathcal{E}_{\partial D}^{\kappa}$ known to be a Poisson point process of Brownian excursions in $D$ \cite{Qian_2019}. For general subcritical intensities, the excursions interact and the full distribution of $\mathcal{E}_{\partial \mathbb{D}}^\kappa$ remains unknown. We can nevertheless consider the first-order intensity measure $\mu_{\partial D}^{\mathcal{E}, \kappa}$ defined as
\begin{equation} \label{eq:def_first_order_mu}
    \mathbb{E} \Big[ \sum_{e \in \mathcal{E}_{\partial D}^\kappa} F(e) \Big] \,:=\, \int F(e)\, \mu_{\partial D}^{\mathcal{E}, \kappa}(\dd e),
\end{equation}
for any nonnegative measurable path-functional $F$. Recall that the Brownian excursion measure $\mu_D^{\mathrm{exc}}$ is the canonical, conformally invariant, sigma-finite measure on boundary-to-boundary Brownian paths in $D$ (see Section~\ref{sec:bubbles} for a formal definition).

\begin{theoreme} \label{th:intensity} For all $\kappa  \in (8/3,4]$, and for all simply connected domains $D \subsetneq \mathbb{C}$, the first-order intensity measure of the excursions from the boundary of a Brownian loop soup cluster is a constant multiple of the excursion measure in $D$,
$$
\mu_{\partial D}^{\mathcal{E}, \kappa} (\dd e) \,=\, \beta(\kappa) \, \mu_D^{\mathrm{exc}}(\dd e). 
$$    
Furthermore, the intensity $\beta(\kappa)$ is explicitly given by,
\begin{equation} \label{eq:g_kappa_explicit}
\beta(\kappa)  \,=\, \frac{(3 \kappa-8)(6-\kappa)}{2 \kappa(4-\kappa)} \tan \Big(\frac{4 \pi}{\kappa}\Big).    
\end{equation}
\end{theoreme}

Each excursion in $D$ connects two boundary points $x, y \in \partial D$. By \cite{lawler2005conformally}, the measure on endpoints induced by $\mu_D^{\mathrm{exc}}$ is exactly given by $H_D(x,y) \dd x \dd y$, where $H_D$ is the boundary Poisson kernel (see Section~\ref{sec:bubbles}) and $\dd x, \dd y$ denote the arc-length measure on $\partial D$. This immediately identifies the first-order intensity of the excursion endpoints as $\beta(\kappa) H_D(x,y) \dd x \dd y$. The coefficient $\beta(\kappa)$ can be interpreted as the effective intensity of excursions touching the boundary: larger values of $\beta(\kappa)$ correspond to a higher density of excursions rooted at $\partial D$. A computation shows that $\beta$ extends continuously to the whole interval $[8/3, 4]$, with
$$
\lim_{\kappa \to 8/3^+} \beta(\kappa) = \frac{5}{2\pi} \quad \quad \text{and} \quad\quad \lim_{\kappa \to 4^-} \beta(\kappa) = \frac{\pi}{4}.
$$ 
The limit as $\kappa \to 4^-$ matches the critical case ($\theta = 1/2$), where boundary excursions form a Poisson point process of intensity $\pi/4$ \cite{Qian_2019}. Conversely, $\kappa \to 8/3^+$ corresponds to $\theta \to 0$, where loops are rare and each cluster reduces to a single Brownian loop. The limiting value $5/(2\pi)$ is consistent, as it heuristically recovers the expected area $\frac{\pi}{5}$ of a Brownian loop of duration $1$ \cite{Garban_2006} and the height gap $\frac{5}{\pi}$ of Brownian motion \cite{jego2026heightgapplanarbrownian} (see Remark~\ref{rem:brownian_area}).

More surprisingly, we prove that the intensity $\beta(\kappa)$ of excursions from the boundary of a cluster is not monotone in $\kappa$. We reveal the existence of a new special parameter $\kappa_0 \in (8/3, 4)$ such that $\beta$ is strictly increasing on $[8/3, \kappa_0)$ and strictly decreasing on $(\kappa_0, 4]$.

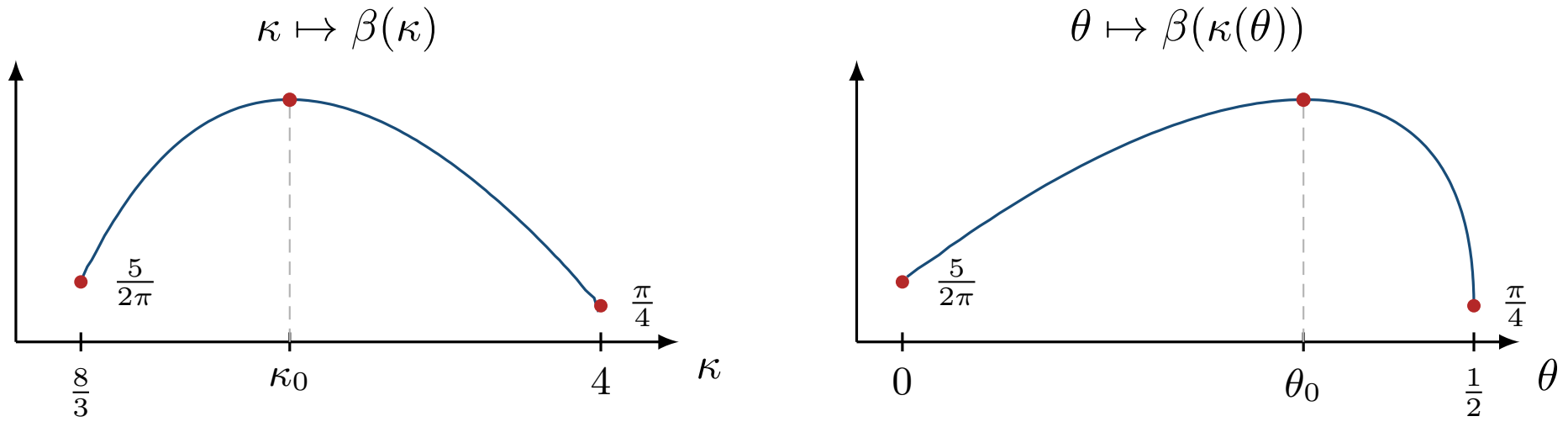
\begin{figure}[htbp]
\centering
\begin{tikzpicture}
  \definecolor{curveblue}{RGB}{24, 76, 120}
  \definecolor{ptred}{RGB}{180, 40, 40}

  \begin{axis}[
    name=ax1,
    width=6.8cm,
    height=3.8cm, 
    axis lines=left,
    axis line style={-latex, semithick},
    tick style={semithick, black},
    title={$\kappa \mapsto \beta(\kappa)$},
    title style={font=\small, yshift=-8pt}, 
    xlabel={$\kappa$},
    xlabel style={at={(ticklabel* cs:1)}, anchor=north west, font=\small},
    xmin=2.5, xmax=4.2,
    ymin=0.245, ymax=0.284,
    xtick={2.6667, 3.2022, 4.0},
    xticklabels={$\frac{8}{3}$, $\kappa_0$, $4$},
    ytick=\empty,
    tick label style={font=\footnotesize},
    clip=false
  ]
    \addplot[
      domain=2.672:3.994,
      samples=120,
      semithick,
      color=curveblue
    ] { ((3*x-8)*(6-x)) / (2*pi*x*(4-x)) * tan(720/x) };

    \draw[densely dashed, gray!60] (axis cs:3.2022, 0.245) -- (axis cs:3.2022, 0.2785);

    \filldraw[ptred] (axis cs:2.6667, 0.2533) circle (1.3pt)
      node[right=3pt, font=\scriptsize, text=black] {$\frac{5}{2\pi}$};

    \filldraw[ptred] (axis cs:3.2022, 0.2785) circle (1.4pt);

    \filldraw[ptred] (axis cs:4.0, 0.2500) circle (1.3pt)
      node[right=2pt, font=\scriptsize, text=black] {$\frac{\pi}{4}$};
  \end{axis}

  \begin{axis}[
    name=ax2,
    at={($(ax1.east)+(1.4cm,0)$)},
    anchor=west,
    width=6.8cm,
    height=3.8cm, 
    axis lines=left,
    axis line style={-latex, semithick},
    tick style={semithick, black},
    title={$\theta \mapsto \beta(\kappa(\theta))$},
    title style={font=\small, yshift=-8pt}, 
    xlabel={$\theta$},
    xlabel style={at={(ticklabel* cs:1)}, anchor=north west, font=\small},
    xmin=-0.04, xmax=0.54,
    ymin=0.245, ymax=0.284,
    xtick={0, 0.3509, 0.5},
    xticklabels={$0$, $\theta_0$, $\frac{1}{2}$},
    ytick=\empty,
    tick label style={font=\footnotesize},
    clip=false
  ]
    \addplot[
      domain=2.671:3.995,
      samples=130,
      semithick,
      color=curveblue
    ] (
      {(3*x-8)*(6-x) / (4*x)},
      {((3*x-8)*(6-x)) / (2*pi*x*(4-x)) * tan(720/x)}
    );

    \draw[densely dashed, gray!60] (axis cs:0.3509, 0.245) -- (axis cs:0.3509, 0.2785);

    \filldraw[ptred] (axis cs:0.0, 0.2533) circle (1.3pt)
      node[right=3pt, font=\scriptsize, text=black] {$\frac{5}{2\pi}$};

    \filldraw[ptred] (axis cs:0.3509, 0.2785) circle (1.4pt);

    \filldraw[ptred] (axis cs:0.5, 0.2500) circle (1.3pt)
      node[right=2pt, font=\scriptsize, text=black] {$\frac{\pi}{4}$};
  \end{axis}

\end{tikzpicture}
\caption{The functions $\kappa \mapsto \beta(\kappa)$ and $\theta \mapsto \beta(\kappa(\theta))$. In both representations, the unique maximum is attained at $\kappa_0 \approx 3.2022$ (corresponding to $\theta_0 \approx 0.3509$).}
\label{fig:g_kappa_theta}
\end{figure}

\noindent Furthermore, we show that this special parameter $\kappa_0$ is the unique solution to the equation 
\begin{equation}
\frac{3}{3 \kappa-8}-\frac{1}{6-\kappa}-\frac{1}{\kappa}+\frac{1}{4-\kappa}=\frac{8 \pi}{\kappa^2 \sin \big(\frac{8 \pi}{\kappa}\big)}.
\end{equation}

Our proof reveals that $\beta(\kappa)$ is in fact a constant multiple of  $\theta(\kappa) \mu(\kappa)$. Here, $\theta(\kappa)$ is the loop-soup intensity, and $\mu(\kappa)$ is the expected $|\log|$-conformal radius of the loop seen from the origin. This explains the non-monotonicity through two competing mechanisms. On the one hand, a higher intensity injects more Brownian loops into the domain, which increases the number of loops that can hit $\partial O$ and therefore generate boundary excursions. On the other hand, $\mu(\kappa)$ is strictly decreasing in $\kappa$. Indeed, as $\kappa$ (or, equivalently, $\theta$) increases, clusters merge and the boundary $\partial O$ expands outwards by swallowing neighboring clusters. Consequently, many loops that would have hit the boundary are now strictly within the core of the cluster, and no longer count as boundary excursions from $\partial O$. Near $\kappa = 8/3$, adding new loops dominates, so $\beta(\kappa)$ increases. Near $\kappa = 4$, the absorption of loops into the core of the cluster dominates, so $\beta(\kappa)$ decreases. The two effects balance each other at the critical parameter $\kappa_0$.

\subsection{Outline of the proofs}

The proofs of Theorem~\ref{th:crossing_exponent} and Theorem~\ref{th:phase_transition_loopsoup} are closely linked. We begin with the proof of Theorem~\ref{th:phase_transition_loopsoup}, which is based on a geometric reduction from 2D loop clusters to 1D random coverings and the theory of regenerative sets. We then turn to the proof of Theorem~\ref{th:crossing_exponent}, which relies on partitioning the non-truncated loop soup into small and large loops and on the estimates established for Theorem~\ref{th:phase_transition_loopsoup}. The proof of Theorem~\ref{th:n_delta} is independent and relies on two-dimensional Brownian bubble estimates and the asymptotic density of separating loops on the cylinder. Finally, the proof of Theorem~\ref{th:intensity} is completely independent  and relies on a mass-conservation argument for the occupation time.

\medskip

\noindent{\large \textbf{Phase transition and crossing probabilities for the truncated loop soup.}} Let $\bpi: \mathfrak{C} \to \mathbb{R}^+$ with $\bpi(x, e^{i\varphi}) := x$ be the longitudinal projection of the cylinder. For any loop $\wp$ on $\mathfrak{C}$, its longitudinal span is given by its minimum $a(\wp) := \min \bpi(\wp)$ and its maximum $b(\wp) := \max \bpi(\wp)$. We call the open interval $(a(\wp), b(\wp))$ the 1D \emph{projection} of $\wp$. The collection of all projections corresponding to loops in $\LC$ forms a one-dimensional random covering of $\mathbb{R}^+$. Clearly, if a cluster of 2D loops crosses from $\{y\} \times \mathbb{S}^1$ to $\{x\} \times \mathbb{S}^1$ on the cylinder (with $y < x$), their 1D projection intervals must form a connected component of intervals crossing $[y, x]$ on $(0,\infty)$. Consequently, the connectivity of the 1D projection provides an upper bound for the 2D loop connectivity. However, this is only an optimistic bound: two loops can easily have overlapping 1D projections without actually intersecting on the 2D cylinder. Therefore, upper-bounding the 2D connectivity by this 1D projection covering only provides a necessary condition for percolation, which is not sufficient to establish the exact crossing probabilities.

To complement this upper bound, we construct a second, more restrictive 1D covering that provides a sufficient condition for 2D intersection. To do so, we restrict our attention to separating loops, i.e., loops that separate the two ends of the cylinder, as they act as topological barriers. 
For any such separating loop $\wp$, we consider its inner boundary $\partial^{\mathrm{i}}\wp$, defined as the boundary of the connected component of $\mathfrak{C} \setminus \wp$ containing infinity. We then define $c(\wp) := \min \bpi(\partial^{\mathrm{i}}\wp)$ as the minimal longitudinal depth reached by this inner boundary. We call the interval $(a(\wp), c(\wp))$ the \emph{core interval} of $\wp$. The introduction of $c(\wp)$ relies on a crucial topological property: one can notice (see Lemma~\ref{lemma:coupling_coverings}) that if the core intervals $(a(\wp), c(\wp))$ and $(a(\wp'), c(\wp'))$ of two separating loops overlap in $\mathbb{R}^+$, then the loops $\wp$ and $\wp'$ \textit{must} intersect on the cylinder $\mathfrak{C}$. The collection of these core intervals thus provides a pessimistic 1D covering that strictly lower-bounds the 2D loop connectivity.

Thanks to these two coverings, the event of crossing a longitudinal distance $x$ on the cylinder is rigorously sandwiched between the connectivity events of two 1D random interval coverings. The initial 2D problem is therefore entirely reduced to the study of these two 1D models.  The core of our proof consists in demonstrating that both the optimistic and the pessimistic coverings belong to the same universality class. Specifically, we show that both 1D coverings undergo a phase transition at the exact same critical intensity $\theta=1$, and that their respective crossing probabilities exhibit the same asymptotic decay.

Once the problem is reduced to 1D Poisson point processes of intervals, the percolation behavior is completely governed by the sets of uncovered points, which form perfect regenerative sets on $\mathbb{R}^+$. Since the analysis of such regenerative sets is well-understood (see e.g \cite{Bertoin1999}, chapter 7), we simply apply the general theory to our case to find the phase transition and the exact asymptotics in the different regimes.
\medskip

\noindent{\large \textbf{Sharp crossing exponent.}} First, under the conformal map $f : re^{i\theta} \mapsto (|\log r|, e^{i\theta})$, crossing the annulus $\rho \mathbb{D} \setminus r\mathbb{D}$ in the disc is equivalent to crossing the slice $(y, x) \times \mathbb{S}^1$ on the cylinder with $y = -\log \rho$ and $x = -\log r$. The lower bound is therefore an immediate consequence of Theorem~\ref{th:phase_transition_loopsoup}, since the truncated loop soup $\LC$ is a subset of the full soup $\mathcal{L}_\mathfrak{C}^\theta$.

For the upper bound, the main challenge is to control the infinite density of small loops. Instead of a duration cutoff, we partition the soup using a spatial diameter threshold $\varepsilon > 0$, which is much better suited here because it provides strict spatial localisation: small loops in disjoint regions become independent.  When we project the large loops (those with a diameter larger than $\varepsilon$) onto the cylinder, they leave uncovered regions along the longitudinal axis, which we call \emph{projection gaps} (see Figure~\ref{fig:multiple_crossings}). By definition, no large loops can enter these gaps. Therefore, for a cluster to cross $(y,x) \times \mathbb{S}^1$, every single gap must be crossed by some cluster consisting only of small loops. 

We prove that the probability for a cluster of small loops to cross a large gap decays exponentially fast. Crucially, since removing large loops recovers some kind of translation invariance along the cylinder, this exponential bound is uniform and depends only on the length of the gap, not on its position on the cylinder. To prevent dependencies between gaps, we select a sub-collection of gaps that are both large and well-separated. This implies that crossing each additional gap with small loops costs a uniform penalty factor $q < 1$. 

The intuition is that configurations with many gaps are exponentially penalized by this factor $q$, and thus become completely negligible. The only configurations that can contribute to the crossing probability are those where the large loops leave very few gaps (or no gaps at all) between $y$ and $x$. By our analysis of the truncated loop soup in Theorem~\ref{th:phase_transition_loopsoup}, the probability that the large loops leave very few gaps across the slice $(y,x) \times \mathbb{S}^1$ decays  as $x^{-1+\theta}$, which yields the upper bound.

\medskip

\noindent{\large \textbf{Finiteness of bounded clusters in the super-critical regime.}}
For $\theta > 1$, we already know from Theorem~\ref{th:phase_transition_loopsoup} that the separating loops form a unique unbounded cluster (because the pessimistic covering, which is composed only of separating loops, undergoes a phase transition at $\theta = 1$). To prove that there are almost surely only finitely many bounded clusters overall, we must show that the remaining non-separating loops are almost entirely absorbed by this giant component. Our strategy is to divide these non-separating loops into two categories based on their inner spatial radius.

First, we rule out the existence of infinitely many microscopic loops. We show that almost surely, only a finite number of non-separating loops have an inner radius smaller than $(\log x)^{-1}$, where $x$ is their longitudinal coordinate. This is established by proving a small-ball probability estimate for the Brownian loop measure, which shows that loops with an extremely small inner radius are too rare to appear infinitely often as $x \to \infty$. Being finite in number, these microscopic loops can only create a finite number of isolated bounded clusters.

Second, we analyze the infinitely many remaining non-separating loops, which all have an inner radius of at least $(\log x)^{-1}$. If any of these loops forms an isolated bounded cluster, it must contain a spatial disc of radius $(\log x)^{-1}$ that is completely unvisited by the separating loops.  To conclude the proof, we demonstrate that the unbounded cluster is asymptotically dense: the total number of unvisited disjoint discs of radius $(\log x)^{-1}$ on the entire infinite cylinder is almost surely finite. 

Since there are only finitely many holes large enough to accommodate them, all but finitely many of the non-separating loops are geometrically forced to intersect the unbounded cluster. They thus coagulate with the giant component, which proves that the total number of bounded clusters is almost surely finite. \medskip

\noindent{\large \textbf{Excursions from the boundary of a cluster.}} We consider the outermost cluster surrounding the origin, with outer boundary $\partial O$. Conditionally on $O$, $\mathcal{L}_\mathbb{D}^\theta$ splits into a fresh loop soup $\mathcal{L}_O$ (strictly) inside $O$ and the remaining loops $\mathcal{L}_{\mathbb{D} \setminus O}$ \cite{Qian_2019}. The proof consists in calculating in two ways the expected occupation time density of $\mathcal{L}_{\mathbb{D} \setminus O}$ at $0$.  On the one hand, using Le Jan's framework, the density of the occupation time at the origin is given by $\frac{1}{\pi} \theta \mu(\theta)$, where $\mu(\theta)$ is the expected $|\log|$-conformal radius of $\partial O$ (which is explicit by \cite{Schramm_2009,sheffield2011conformalloopensemblesmarkovian}). On the other hand, only the boundary excursions $\mathcal{E}_{\partial O}$ can contribute to the occupation density. A direct calculation shows that an excursion point process of intensity $\beta \mu_O^{\mathrm{exc}}$ has a constant expected occupation density of $2 \beta$ inside $O$. Equating both expressions identifies  $2\pi \beta(\kappa) = \theta(\kappa) \mu(\kappa)$.

\subsection{Further results} \label{sec:further_results}

Our Theorem~\ref{th:crossing_exponent} implies further results on the conformally invariant fields out of Brownian loop soups $h_\theta$ introduced in \cite{jego2023conformallyinvariantfieldsbrownian}, on fractal properties of clusters and on polar sets of Brownian loop soup clusters. These are proved in Appendix~\ref{sec:proof_cors}. \newline

\noindent{\textbf{\large{Conformally invariant fields out of Brownian loop soup.}}}
In \cite{jego2023conformallyinvariantfieldsbrownian}, Jego, Lupu and Qian introduced and studied a new class of conformally invariant fields $h_\theta$ that generalize the Gaussian free field (corresponding to the critical intensity $\theta = 1/2$) to any subcritical $\theta \in (0, 1/2)$. Informally, these fields are constructed as the local time of the Brownian loop soup $\mathcal{L}_{\mathbb{D}}^\theta$ in the unit disc $\mathbb{D}$ to the power $1-\theta$ with independent random signs assigned to each cluster of loops. As already mentioned, in \cite{Jego_2025}, the authors proved that the crossing probability in Brownian loop soup decays as $|\log r|^{-1+\theta + o(1)}$ as $r \to 0$. Because of this $o(1)$ error in the exponent, several key properties of the field $h_\theta$ established in \cite{jego2023conformallyinvariantfieldsbrownian} could only be stated up to that $o(1)$ error. More precisely, the two-point covariance function $C_\theta(x, y) = \mathbb{E}[h_\theta(x)h_\theta(y)]$ exhibited a diagonal blow-up of the form $C_\theta(x, y) = |\log |x-y||^{2(1-\theta) + o(1)}$ as $x \to y$ (Theorem 1.3 in \cite{jego2023conformallyinvariantfieldsbrownian}), and the Minkowski content of the clusters could only be shown to exist in the implicit gauge function $r \mapsto r^2 |\log r|^{1-\theta + o(1)}$ (Theorem 1.9 in \cite{jego2023conformallyinvariantfieldsbrownian}).

Since Theorem~\ref{th:crossing_exponent} removes that $o(1)$ error from the crossing exponent, and since all proofs in \cite{jego2023conformallyinvariantfieldsbrownian} rely on the decay of these crossing probabilities, one can feed our sharper estimate into their framework. This yields the following improvements. First, we get the exact order of the diagonal singularity of the two-point function.

\begin{corollaire}
\label{cor:exact_covariance}
Let $D$ be a bounded simply connected domain. For $\theta \in (0, 1/2]$, denote by $C_\theta(x, y)$ the two-point covariance function of $h_\theta$ in $D$. Then, as $|x - y| \to 0$ while staying at positive distance from $\partial D$,
\begin{equation}
    C_\theta(x, y) \,\asymp\,  |\log |x - y||^{2(1-\theta)}.
\end{equation}
\end{corollaire}

Second, for some measurable and well-chosen enumeration $(\mathcal{C}_j)_{j \geq 1}$ of the clusters of $\mathcal{L}_D^\theta$, the authors proved that the conformally invariant fields $h_\theta$ can be described somewhat more geometrically as the signed sum $\sum_{j \geq 1} \sigma_{\mathcal{C}_j} \mu_j$, where $(\sigma_{\mathcal{C}_j})_{j \geq 1}$ are i.i.d. symmetric signs and $\mu_j$ is the Minkowski content measure of $\mathcal{C}_j$ in some common implicit gauge function $\phi(r) = r^2 |\log r|^{1-\theta+o(1)}$. Feeding our sharper crossing estimate into their construction yields the following improvement on the gauge.

\begin{corollaire} \label{cor:Minkowski} For all $\theta \in (0,1/2]$, the Minkowski gauge function $\phi$ of the clusters of $\mathcal{L}_D^\theta$ satisfies $\phi(r) \asymp r^2 |\log r|^{1-\theta}$ as $r \to 0$.
\end{corollaire}

Since Theorem~\ref{th:crossing_exponent} does not provide an exact asymptotic for the crossing probabilities, the Minkowski gauge function of Brownian loop soup clusters remains implicit, and a priori, we cannot rule out microscopic oscillations in $\phi(r)$. 
\newline

\noindent{\large{\textbf{Polar sets of loop soup clusters.}}}  Recall that a set $A \subset \mathbb{C}$ is said to be polar for planar Brownian motion if for all starting points $x$, we have $\mathbb{P}(\exists t > 0 : B_t \in A) = 0$. Since there are countably many loops in a Brownian loop soup, the polar sets of Brownian motion will, almost surely, not be visited by a cluster either. However, in \cite[Theorem 1.5]{Jego_2025}, the authors pointed out that the closures of the loop soup clusters can hit much thinner sets. They defined a set $A$ as being polar for the clusters of $\mathcal{L}_\mathbb{D}^\theta$ if, almost surely, no closure $\overline{\mathcal{C}}$ of any cluster of $\mathcal{L}_\mathbb{D}^\theta$ intersects $A$.  Following their framework, for any measurable function $K : \mathbb{C} \times \mathbb{C} \to [0, \infty]$, let 
$$
\operatorname{Cap}_K(A):=\left(\inf \left\{\int K(x, y) \mu(\mathrm{d} x) \mu(\mathrm{d} y): \mu \text { probability measure on } A\right\}\right)^{-1}
$$
be the capacity of $A$ with respect to the kernel $K$ and $\operatorname{Cap}_{\log^{1-\theta}}$ the capacity associated to the kernel $K_\theta(x, y) = |\log |x-y||^{1-\theta}$. Recall that Kakutani's theorem (see e.g \cite[Theorem 8.20]{morters2010brownian} and \cite{Kak44}) states that a set $A$ is polar for planar Brownian motion if and only if $\operatorname{Cap}_{\log^1}(A) = 0$.  In~\cite[Theorem~1.5]{Jego_2025}, the authors proved that for any closed set $A \subset \mathbb{D}$, if $A$ is polar for the clusters of $\mathcal{L}_\mathbb{D}^\theta$, then $\operatorname{Cap}_{\log^\alpha}(A) = 0$ for all $\alpha > 1-\theta$, whereas if $A$ is not polar, then $\operatorname{Cap}_{\log^\alpha}(A) > 0$ for all $\alpha < 1-\theta$. Using our sharper crossing exponent, we obtain the exact necessary condition for polarity at the critical exponent $\alpha = 1-\theta$.

\begin{corollaire} \label{cor:polar_sets} If $A$ is (closed and) polar for the clusters of $\mathcal{L}_\mathbb{D}^\theta$, then $\operatorname{Cap}_{\log^{1-\theta}}(A) = 0$.
\end{corollaire}

This generalizes one direction of Kakutani's theorem (which informally corresponds to the limiting case $\theta \to 0$) to Brownian loop soup clusters. However, this method does not yield the full equivalence. Whether $\operatorname{Cap}_{\log^{1-\theta}}(A) = 0$ is also sufficient for $A$ to be polar remains an open question (even in the critical case $\theta = 1/2$). \newline

\noindent {\large \textbf{Universality of the phase transition.}}  Although Theorem~\ref{th:phase_transition_loopsoup} is stated for the duration threshold $T(\wp) \ge \tau$, our proof method does not rely on this specific choice and should adapt to a wider class of truncated loop soups $\mathcal{L}_{\mathfrak{C}}^{\theta, F \ge \delta} := \{\wp \in \mathcal{L}_{\mathfrak{C}}^\theta : F(\wp) \ge \delta\}$. Indeed, the reduction to one-dimensional random coverings developed in Section~\ref{sec:phase_transition} only requires the functional $F$ to satisfy three natural properties that can be stated in terms of the Brownian bubble measure $\mu_{\mathfrak{C}}^{\mathrm{bub}}$ on the cylinder (see Section~\ref{sec:bubbles} for a definition of this measure):
\begin{itemize}[itemsep=3pt, topsep=5pt, parsep=0pt]
    \item[(1)] \emph{Longitudinal invariance:} $F(\tau_u(\wp)) = F(\wp)$ for all $u \ge 0$.
    \item[(2)] \emph{UV cutoff:} $\mu_{\mathfrak{C}}^{\mathrm{bub}}(F(\wp) \ge \delta) < \infty$, which filters out the infinite density of microscopic loops and prevents the real line from being trivially covered everywhere.
    \item[(3)] \emph{IR preservation:} $\mu_{\mathfrak{C}}^{\mathrm{bub}}(\ell(\wp) > x \text{ and } F(\wp) < \delta) = O(e^{-cx})$, which ensures that unusually large loops are statistically not affected by the truncation.
\end{itemize}

These conditions are satisfied by a wide range of functionals: metric cutoffs such as the spatial diameter $\mathrm{diam}(\wp)$, the longitudinal span $\ell(\wp)$, the inner radius $\mathrm{inrad}(\wp)$, or the area $\mathrm{Area}(\mathrm{Int}(\wp))$; topological cutoffs such as separating cutoffs $\mathbf{1}_{\wp \text{ is separating}}$ and the winding number $|\mathrm{wind}(\wp)|$, occupation quantities such as the maximum longitudinal local time $\sup_{z} \ell_z(\bpi (\wp))$ and conformal quantities such as the capacity. In each of these cases, the associated 1D projection and core coverings should satisfy Hypothesis~(H) from Section~\ref{sec:random_coverings}, yielding the exact same critical intensity $\theta = 1$ and the same crossing exponents.

Finally, this intensity $\theta = 1$ is also related to the spatial Markov property of the Brownian loop soup when the loops are viewed as oriented loops \cite{werner2015spatialmarkovpropertysoups}, and is known to be special due to its connection to the uniform spanning tree (see \cite{PW98,Wil96}), though to our knowledge no link with the phase transition has been established.
\newline

\noindent {\large \textbf{Organisation of the paper.}}   
In Section~\ref{sec:prelim}, we collect the necessary preliminaries and definitions concerning Brownian loops, bubbles, excursions and percolation properties of random coverings. 
Section~\ref{sec:phase_transition} is dedicated to the proof of Theorem~\ref{th:phase_transition_loopsoup} for the truncated loop soup on the cylinder: we detail the geometric reduction of 2D loop clusters to 1D random interval coverings and establish the crucial sandwiching property.
In Section~\ref{sec:crossing_exponent}, we prove Theorem~\ref{th:crossing_exponent}, and explain how clusters of small loops must cross every gap left by the large loops. 
Section~\ref{sec:finiteness_bounded} addresses the supercritical regime $\theta > 1$ and the proof of Theorem~\ref{th:n_delta}. We prove that the unbounded component is asymptotically dense by establishing the finiteness of unvisited microscopic gaps and ruling out the existence of infinitely many loops with tiny inner radii. In Section~\ref{sec:excursions_boundary_proof}, we prove Theorem~\ref{th:intensity}.

In Appendix~\ref{sec:ai}, we describe our use of artificial intelligence tools in the preparation of this paper. Appendix~\ref{sec:proof_cors} is dedicated to the proof of Corollaries~\ref{cor:exact_covariance}, \ref{cor:Minkowski} and \ref{cor:polar_sets}. Finally, Appendix~\ref{sec:appendix} gathers the proofs of several technical lemmas concerning Brownian excursions, the heat kernel on the circle, and complex analysis.

\paragraph{Acknowledgements} The author warmly thanks his PhD supervisor, Titus Lupu, for many insightful discussions, explanations, corrections, and for reviewing earlier draft versions of this paper. The author is supported by a PhD grant from Université Paris Cité through the École Doctorale de Sciences Mathématiques de Paris Centre (ED 386).

\section{Preliminaries} \label{sec:prelim}

In this section, we gather some definitions, intuitions and preliminary results about random coverings of $\mathbb{R}_+$, Brownian loops, Brownian bubbles, Brownian excursions and SLE$_\kappa$ bubbles in the disc and on the semi-infinite cylinder $(0,\infty) \times \mathbb{S}^1$.

\subsection{Brownian loop soup on the cylinder.} \label{sec:prelim_loops}

Consider a Brownian loop soup $\mathcal{L}_{\mathbb{D}}^\theta$ of intensity $\theta > 0$ in the unit disc $\mathbb{D}$. Because the punctured disc $\mathbb{D} \setminus \{0\}$ and the semi-infinite cylinder $\mathfrak{C} = (0, \infty) \times \mathbb{S}^1$ are conformally equivalent via the map $f : re^{i\theta} \mapsto (|\log r|, e^{i\theta})$ and because the Brownian loop soup $\mathcal{L}_{\mathbb{D}}^\theta$ is conformally invariant, we can work in $\mathfrak{C}$ instead of $\mathbb{D}$. The loop measure is then given by
$$
\mu_\mathfrak{C}^{\text{loop}}(d\wp) \,=\, \int_0^\infty dt \int_\mathfrak{C} \dd z \frac{p_\mathfrak{C}(t,z,z)}{t} \mathbb{P}_\mathfrak{C}^{t,z,z}(d\wp).
$$
One advantage of this representation is the following: on the cylinder, the two components of a Brownian loop $\wp$ of duration $t$ at $(x,e^{i\theta})$ are simply two independent Brownian bridges of duration $t$: the longitudinal one from $x$ to $x$ and the angular one from $e^{i \theta}$ to $e^{i \theta}$ on $\mathbb{S}^1$. In what follows, we are going to introduce topological definitions for loops on the cylinder. \newline

\noindent \textbf{Inner and outer boundaries of loops.} Notice that a loop $\wp'$ surrounds $0$ in the unit disc $\mathbb{D}$ if and only if its pushforward $\wp = f(\wp')$ on the cylinder wraps around the cylinder (in that case, we say that $\wp$ is \emph{separating}). For such a loop that surrounds the origin, we can define its inner and outer boundaries as follows: define $\partial^{\text{o}} \wp'$ as the boundary of the unique unbounded connected component of $\mathbb{C} \setminus \wp'$ and $\partial^{\text{i}} \wp'$ as the boundary of the unique connected component of $\mathbb{C} \setminus \wp'$ containing $0$. Similarly, we define $\partial^{\text{o}} \wp$ and $\partial^{\text{i}} \wp$ as their pushforward by $f$:
$$
\partial^{\text{o}} \wp := f(\partial^{\text{o}} \wp') \quad \quad \partial^{\text{i}} \wp := f(\partial^{\text{i}} \wp').
$$
Equivalently, the boundaries $\partial^{\text{o}} \wp$ and $\partial^{\text{i}} \wp$ can also be defined solely on the cylinder as the boundary of the unique connected component of $\mathfrak{C} \setminus \wp$ containing $\{0\} \times \mathbb{S}^1$ and as the boundary of the unique connected component of $\mathfrak{C} \setminus \wp$ containing $\{\infty\} \times \mathbb{S}^1$. See Figure~\ref{fig:loops} below. Note that while a non-zero winding number guarantees that the loop separates $\{0\} \times \mathbb{S}^1$ from $\{\infty\} \times \mathbb{S}^1$, the converse is false: a loop can easily disconnect the two ends of the cylinder while having zero winding number. Thus, being separating on the cylinder is not equivalent to being non-contractible (i.e., having a non-zero winding number).

\begin{figure}[h]
    \centering
    \includegraphics[width=12cm]{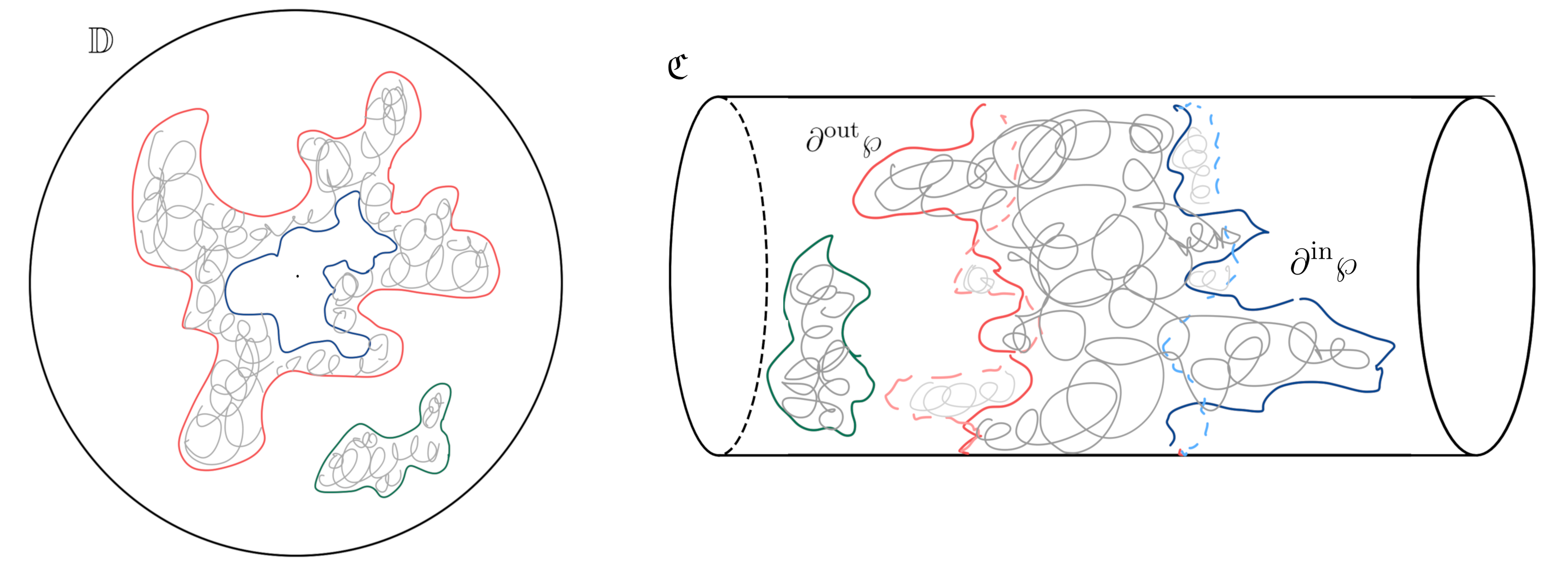}
    
    \caption{Two Brownian loops in $\mathbb{D} \setminus \{0\}$ conformally mapped to the semi-infinite cylinder $\mathfrak{C}$. The green one does not surround $0$ in the disc, so it is non-separating on the cylinder. The second one surrounds zero: its inner and outer boundaries are denoted $\partial^{\text{o}} \wp$ and $\partial^{\text{i}} \wp$.}
    \label{fig:loops}
\end{figure}

Let $\bpi : \mathfrak{C} \to \mathbb{R}_{+}$ be the longitudinal projection on the cylinder, defined by $\bpi(x, e^{i\varphi}) := x$. For any loop $\wp$ on the cylinder $\mathfrak{C}$, we define the following three topological quantities,
\begin{equation} \label{eq:def_a_b_c}
a(\wp) = \min \bpi(\wp) \quad\quad b(\wp) = \max \bpi(\wp) \quad\quad c(\wp) = \min \bpi(\partial^{\text{i}} \wp).
\end{equation}
The quantity $c(\wp)$ is only well-defined when $\wp$ is separating. The open interval $(a(\wp), b(\wp))$ represents the longitudinal span of the loop $\wp$ and $c(\wp)$ represents the minimal longitudinal position of $\wp$'s inner boundary on the cylinder. Alternatively, writing $\wp' = f^{-1}(\wp)$ on the unit disc, we have $c(\wp) = -\log \max \{|z|, z \in \partial^{\text{i}} \wp'\}$. Hence, in the disc $\mathbb{D}$, $c(\wp')$ represents (the negative logarithm of) the maximal distance of $\partial^{\text{i}} \wp'$ to the origin. \newline

\noindent{\textbf{Metric on the cylinder and diameter.}}  Throughout this paper, we will frequently use discs on the cylinder. For $(e^{i\theta}, e^{i\varphi}) \in \mathbb{S}^1 \times \mathbb{S}^1$, denote $d_{\mathbb{S}^1}(e^{i\theta}, e^{i\varphi}) = \min_{k \in \mathbb{Z} }| \theta - \varphi - 2\pi k|$ the geodesic distance on the unit circle. We define the canonical distance on the cylinder $\mathfrak{C}$ by
\begin{equation} \label{eq:def_metric_cyl}
d_\mathfrak{C}\big((x,e^{i\theta}), (y,  e^{i\varphi}) \big) := \sqrt{(x-y)^2 + d_{\mathbb{S}^1}(e^{i\theta},  e^{i\varphi})^2}.
\end{equation}
We denote by $D_\mathfrak{C}(z, r)$ the associated open disc of center $z \in \mathfrak{C}$ and radius $r > 0$. As long as $r < \pi$, these discs behave like standard planar discs: their area scales as $r^2$ and their Lebesgue measure is invariant under rotation and translation. For a loop $\wp$ on the cylinder, we define its diameter $\mathrm{diam}(\wp)$ as the diameter $2r$ of the smallest disc $D_\mathfrak{C}(z,r)$ that contains $\wp$.

\subsection{Brownian bubbles and SLE$_\kappa$-bubbles.} \label{sec:bubbles}

In this section, we introduce and recall some properties of the Brownian bubble measure, as well as the $\mathrm{SLE}_\kappa$ bubble measures in the unit disc. \medskip

\noindent \textbf{Boundary Poisson kernel, excursions, and Brownian bubbles.} 
For two distinct boundary points $a, b \in \partial \mathbb{D}$ and an interior point $x \in \mathbb{D}$, the \emph{boundary Poisson kernel} and the \emph{ Poisson kernel} in the unit disc are (respectively) given by
$$
H_\mathbb{D}(a,b) \,=\, \frac{1}{\pi |a-b|^2} \quad\quad\text{and}\quad\quad H_\mathbb{D}(x,a) \,=\, \frac{1}{2\pi} \frac{1-|x|^2}{|x-a|^2}.
$$
For any $a \neq b \in \partial \mathbb{D}$, we denote by $\mathbb{P}_\mathbb{D}^{a,b}$ the probability measure of a standard Brownian excursion in $\mathbb{D}$ linking $a$ to $b$, and define the weighted boundary-to-boundary excursion measure as $\mu_\mathbb{D}(a,b) := H_\mathbb{D}(a,b) \mathbb{P}_\mathbb{D}^{a,b}$. The canonical (unrooted) \emph{Brownian excursion measure} $\mu_\mathbb{D}^{\mathrm{exc}}$ in the unit disc is then obtained by integrating over all boundary endpoints:
\begin{equation} \label{eq:def_excur}
\mu_\mathbb{D}^{\mathrm{exc}} \,:=\, \iint_{(\partial \mathbb{D})^2} \mu_\mathbb{D}(a,b) \, \dd a \dd b,    
\end{equation}
where $\dd a$ and $\dd b$ denote the standard arc-length measure on $\partial \mathbb{D}$. When the two endpoints merge, this excursion measure gives rise to Brownian bubbles: for any $w \in \partial \mathbb{D}$, we define the \emph{rooted Brownian bubble measure} $\mu_\mathbb{D}^{\mathrm{bub}}(w)$ at $w$ as the weak limit
$$
\mu_\mathbb{D}^{\mathrm{bub}}(w) \,:=\, \lim_{b \to w} \mu_\mathbb{D}(w,b) \,=\, \lim_{b \to w} H_\mathbb{D}(w,b) \mathbb{P}_\mathbb{D}^{w,b}.
$$
Let us stress that our normalisation of the rooted Brownian bubble measure differs from the one in \cite{Lawler_2003, Jego_2025} (and most of the standard literature): if we denote $\widetilde{\mu}_\mathbb{D}^{\text{bub}}(w)$ the measure constructed in \cite{Lawler_2003} (Section 3.4.1) and conformally mapped to the unit disc by $\psi : \mathbb{H} \to \mathbb{D}$ sending $0 \in \partial \mathbb{H}$ to $w \in \partial \mathbb{D}$, we have 
$$
\mu_\mathbb{D}^{\text{bub}}(w) \,=\, \frac{1}{\pi}\widetilde{\mu}_\mathbb{D}^{\text{bub}}(w).
$$
The main reason for this choice is to avoid a prefactor $\frac{1}{\pi}$ in Equation~\eqref{eq:loop_bubble_decomp}. Since we will use this equation extensively, this convention avoids carrying an extra factor throughout the paper. The \emph{unrooted Brownian bubble measure} in the unit disc (or simply Brownian bubble measure) is then obtained by integrating over the roots $w \in \partial \mathbb{D}$:
\begin{equation} \label{eq:def_bubble_measure_disc}
\mu_\mathbb{D}^{\text{bub}} \,:=\, \int_{0}^{2\pi} \mu_\mathbb{D}^{\text{bub}}(e^{i \theta}) \dd \theta.
\end{equation}
As our normalisation of rooted bubbles is different from the one in \cite{Lawler_2003}, our normalisation of unrooted Brownian bubbles also differs by a factor $\frac{1}{\pi}$. Finally, by conformally mapping the bubbles in the unit disc to the cylinder $\mathfrak{C}$, we define the rooted and unrooted Brownian bubbles in $\mathfrak{C}$ as 
\begin{equation}
\mu_\mathfrak{C}^{\text{bub}}(e^{i\varphi}) \,:=\, f \circ \mu_{\mathbb{D}}^{\text{bub}}(w) \quad \text{ and } \quad  
\mu_\mathfrak{C}^{\text{bub}} \,:=\, f \circ \mu_{\mathbb{D}}^{\text{bub}},
\end{equation}
where $f$ conformally maps $\mathbb{D} \setminus \{0\}$ to the cylinder $\mathfrak{C}$ and sends $w \in \partial \mathbb{D}$ to $(0,e^{i\varphi}) \in \partial \mathfrak{C}$. This definition is valid since the Brownian bubble measures are conformally covariant and $f'$ has modulus $1$ on the boundary $\partial \mathbb{D}$. 
\medskip

\noindent \textbf{Brownian bubble decomposition of the loop measure.} By re-rooting each loop $\wp$ at its minimal longitudinal position $a(\wp) = \min \bpi (\wp)$, the shifted trajectory $\tau_{-a}(\wp)$ becomes a Brownian bubble rooted on the boundary $\partial \mathfrak{C} = \{0\} \times \mathbb{S}^1$, where $\tau_a(x, e^{i\theta}) := (x+a, e^{i\theta})$. Following Lawler and Werner, the loop measure on $\mathfrak{C}$ disintegrates over this minimum $a$ as
\begin{equation}
    \mu_{\mathfrak{C}}^{\mathrm{loop}} = \int_0^\infty \dd a \, \tau_a \circ \mu_{\mathfrak{C}}^{\mathrm{bub}}.
    \label{eq:loop_bubble_decomp}
\end{equation}  
This decomposition is similar to \cite[Proposition 8]{lawler2003brownianloopsoup} except that in our normalisation, the loop measure is twice as large as in \cite{lawler2003brownianloopsoup} (the critical intensity is $\theta = 1/2$ whereas it is $c = 2\theta = 1$ for them) and our normalisation of the bubble measure contains a factor $\frac{1}{\pi}$.

Under $\mu_{\mathfrak{C}}^{\mathrm{bub}}$, shifting the time origin of the loop to the hitting time of its longitudinal minimum decouples the trajectory $\wp = (X, \Theta)$ into an excursion and a bridge: the longitudinal component $X$ is a 1-dimensional excursion under Itô's measure $\mathbf{n}(\dd X)$ of lifetime $\zeta = \zeta(X)$, and, conditionally on $\zeta$, the angular component $\Theta$ is an independent Brownian bridge on $\mathbb{S}^1$ of duration $\zeta$ rooted at some $\theta \in \mathbb{S}^1$ (see \cite{Lupu_2018} Proposition 3.18). Taking into account the circular return density given by the heat kernel on $\mathbb{S}^1$  $p_{\mathbb{S}^1}(\zeta) := p_{\mathbb{S}^1}(\zeta, \theta, \theta)$ and integrating over all root positions $\theta \in \mathbb{S}^1$, the bubble measure on $\mathfrak{C}$ is explicitly given by
\begin{equation}
    \mu_{\mathfrak{C}}^{\mathrm{bub}}(\dd X, \dd \Theta) = p_{\mathbb{S}^1}(\zeta)  \int_{\mathbb{S}^1} \dd \theta \, \mathbf{n}(\dd  X) \otimes \mathbb{P}_{\mathbb{S}^1}^{\zeta, \theta, \theta}(\dd \Theta).
    \label{eq:bubble_explicit}
\end{equation} 

\noindent \textbf{SLE$_\kappa$ bubble measure and loop boundaries.} For any boundary root $w \in \partial\mathbb{D}$ and $\kappa \in (0, 8)$, the $\mathrm{SLE}_\kappa$ bubble measure $\nu_{\mathbb{D}}^{\mathrm{SLE}_\kappa} (w)$ rooted at $w$, first introduced in the conformal restriction framework of Lawler, Schramm, and Werner \cite{Lawler_2003} (see also \cite{lawler2005conformally}), is defined as the renormalized weak limit of chordal $\mathrm{SLE}_\kappa$ probability measures $\mathbb{P}_{\mathrm{SLE}_\kappa}^{\mathbb{D},w,w_\varepsilon}$ connecting $w$ to the point $w_\varepsilon = e^{i \varepsilon}w \in \partial\mathbb{D}$, namely 
\begin{equation} \label{eq:def_nu_SLE}
\nu_{\mathbb{D}}^{\mathrm{SLE}_\kappa} (w) := \lim_{\varepsilon \to 0} \varepsilon^{-(8/\kappa-1)}\mathbb{P}_{\mathrm{SLE}_\kappa}^{\mathbb{D},w,w_\varepsilon}
\end{equation}
In \cite[Section 10]{Lawler_2003} Lawler, Schramm, and Werner proved that the outer boundary $\partial^{\text{o}} \wp$ under $\mu_{\mathbb{D}}^{\mathrm{bub}}(w)$ coincides, up to a deterministic multiplicative constant $c > 0$, with the $\mathrm{SLE}_{8/3}$ bubble measure rooted at $w$: 
$$(\partial^{\text{o}})_* \mu_{\mathbb{D}}^{\mathrm{bub}}(w) = c \, \nu_{\mathbb{D}}^{\mathrm{SLE}_{8/3}} (w).$$
By integrating over the boundary root $w \in \partial\mathbb{D}$, this identity extends immediately to the total integrated bubble measures $\mu_{\mathbb{D}}^{\mathrm{bub}} := \int_{\partial\mathbb{D}} \mu_{\mathbb{D}}^{\mathrm{bub}}(w) |dw|$ and $\nu_{\mathbb{D}}^{\mathrm{SLE}_\kappa}  := \int_{\partial\mathbb{D}} \nu_{\mathbb{D}}^{\mathrm{SLE}_\kappa} (w) |dw|$. We now turn our attention to the following lemma, which is a direct corollary of Beffara's results.

\begin{lemma} \label{lem:SLE_bubble_hit}
For all $\kappa \in (0,8)$, the following holds true as $r \to 0$,
$$
\nu_{\mathbb{D}}^{\mathrm{SLE}_\kappa}  \big( \gamma \cap (r \mathbb{D}) \neq \emptyset \big) \,\asymp\, r^{1-\kappa/8}.
$$
\end{lemma}

\begin{proof}
By \cite[Corollary 5]{Beffara_2008}, for all $\kappa \in (0,8)$, $a,b \in \partial \mathbb{D}$, and all radii $r<1/2$, we have
\begin{equation} \label{eq:beffara}
\mathbb{P}_{\mathrm{SLE}_\kappa}^{\mathbb{D},a,b} \big( \gamma \cap (r \mathbb{D})\neq \emptyset \big) \,\asymp\, r^{1-\kappa/8} \left( \omega(ab) \wedge \omega(ba) \right)^{8/\kappa-1}.
\end{equation}
Here $\omega(ab)$ denotes the harmonic measure seen from the origin of the arc $ab$ on $\partial \mathbb{D}$. By Equation~\eqref{eq:def_nu_SLE}, let $w \in \partial \mathbb{D}$ and for $\varepsilon > 0$, let $w_\varepsilon = e^{\pm i \varepsilon} w$ be one of the two points on the circle $\partial \mathbb{D}$ at distance $\varepsilon$ from $w$. Then, we have $\omega(ww_\varepsilon) = \frac{\varepsilon}{2\pi}$ and $\omega(w_\varepsilon w) = 1-\frac{\varepsilon}{2\pi}$. We thus get $
\omega(ww_\varepsilon) \wedge \omega(w_\varepsilon w) \,\asymp\, \varepsilon$. Injecting this into Equation~\eqref{eq:beffara} yields,
$$
\mathbb{P}_{\mathrm{SLE}_\kappa}^{\mathbb{D},w,w_\varepsilon} \big( \gamma \cap (r \mathbb{D})\neq \emptyset \big) \,\asymp\, r^{1-\kappa/8} \,\varepsilon^{8/\kappa-1}.
$$
By definition of the $\mathrm{SLE}_\kappa$ bubble measure $\nu_{\mathbb{D}}^{\mathrm{SLE}_\kappa} (w)$ from Equation~\eqref{eq:def_nu_SLE}, we get
$$
\nu_{\mathbb{D}}^{\mathrm{SLE}_\kappa}(w)  \big( \gamma \cap (r \mathbb{D}) \neq \emptyset \big) \,=\, \lim_{\varepsilon \to 0} \varepsilon^{-(8/\kappa-1)}\mathbb{P}_{\mathrm{SLE}_\kappa}^{\mathbb{D},w,w_\varepsilon} \big( \gamma \cap (r \mathbb{D})\neq \emptyset \big) \,\asymp\, r^{1-\kappa/8}.
$$
Integrating over all roots $w \in \partial \mathbb{D}$ yields the desired result.
\end{proof}

\subsection{Excursions from the boundary of a cluster} \label{sec:excursions_boundary}

Let $\mathcal{L}_\mathbb{D}^\theta$ be a Brownian loop soup in the unit disc of intensity $\theta \in (0,\frac{1}{2}]$. We are interested in one of the outermost clusters in the soup, say, without loss of generality (by conformal invariance), the cluster $\mathcal{C}$ surrounding the origin. We denote by $\partial O$ the outer boundary of $\mathcal{C}$ and $O$ its interior. We know from \cite{sheffield2011conformalloopensemblesmarkovian} that the outer boundary $\partial O$ is a $\mathrm{CLE}_{\kappa}$ loop but are mainly interested in what happens inside the topological closure $\overline{O}$.

One can naturally split the loops $\wp \in \mathcal{L}_\mathbb{D}^\theta$ into two disjoint categories: those that stay strictly within $O$ without ever intersecting $\partial O$ (we call this set of loops $\mathcal{L}_O$) and those that touch the boundary $\partial \mathcal{O}$ at least once (which we call $\mathcal{L}_{\partial O}$). In \cite{Qian_2019}, Qian and Werner proved that, conditionally on $\partial O$, the loops $\mathcal{L}_O$ form a fresh loop soup in $O$ and that the law of $\mathcal{L}_{\partial O}$ is conformally invariant, in the sense that, for any $\partial O$-measurable conformal map $\psi$ from $O$ onto $\mathbb{D}$, the law of $\psi(\mathcal{L}_{\partial O})$ is independent of $\partial O$.

\begin{figure}[h]
    \centering
    \includegraphics[width=11cm]{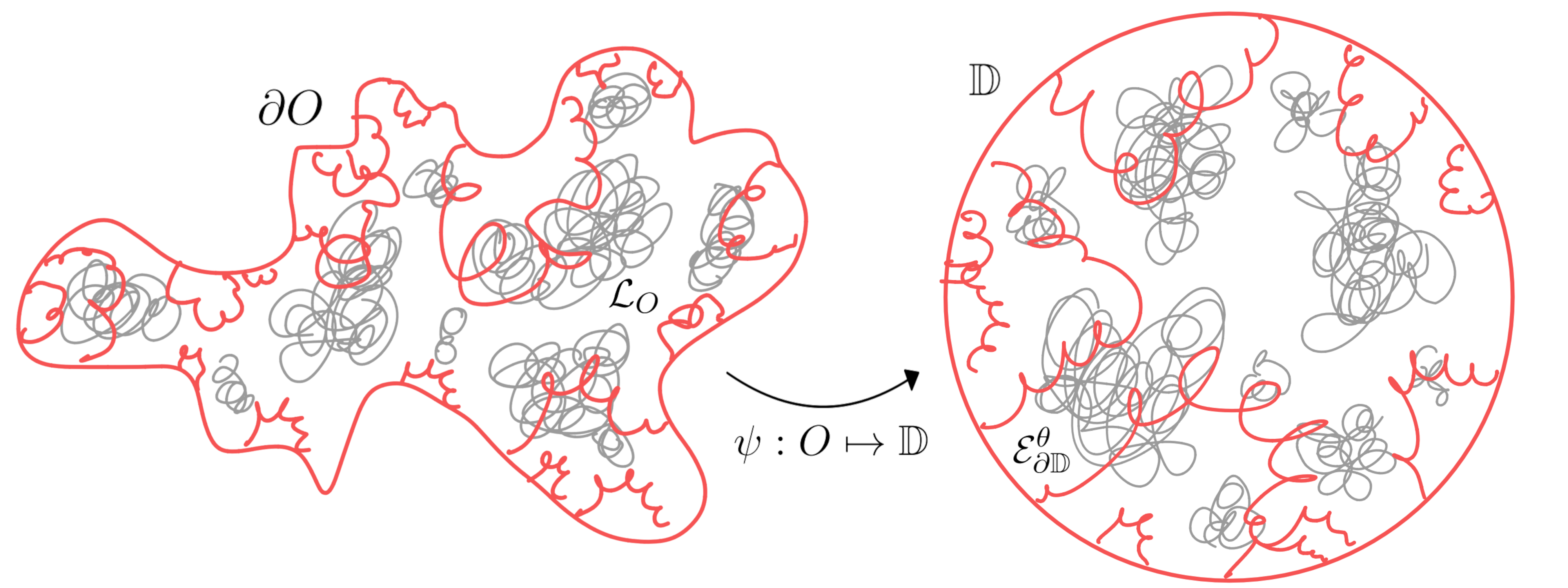}
    
    \caption{The outermost cluster $\mathcal{C}$ delimited by its outer boundary $\partial O$ (in red). Inside $O$, there is a fresh loop soup $\mathcal{L}_O$ (in gray) together with the excursions $\mathcal{E}_{\partial O}$ (in red). Right: the excursions $\mathcal{E}_{\partial \mathbb{D}}^\theta$ and the loop soup conformally mapped to the disc.}
    \label{fig:boundary_excursions_map}
\end{figure}

Each loop $\wp \in \mathcal{L}_{\partial O}$ hits the outer-boundary $\partial O$ a countable number of times. Between two such consecutive contact times, the Brownian trajectory  penetrates inside $O$ before returning to the boundary, forming an excursion in $O$ linking two contact points of the boundary. We denote by $\mathcal{E}_{\partial O}$ the set of these excursions in the domain $O$. By conformal invariance, we can also define the excursions from the boundary $\partial \mathbb{D}$ as 
$$
\mathcal{E}_{\partial \mathbb{D}}^\theta \,:=\: \psi(\mathcal{E}_{\partial O}).
$$
Here $\psi$ is any $\partial O$-measurable conformal map from $O$ onto $\mathbb{D}$ (see Figure~\ref{fig:boundary_excursions_map}). Due to the links between the critical Brownian loop soup and the GFF, the law of $\mathcal{E}_{\partial \mathbb{D}}^\theta$ is only explicit in the case $\theta = 1/2$ (in that case, its trace is distributed as the trace of a Poisson point process of excursions of intensity $\pi/4$, see \cite[Theorem 1]{Qian_2019}). However, for general intensities, such a Poissonian decomposition is not expected (see, again, \cite[Section 5]{Qian_2019}).

For each excursion $e \in \mathcal{E}_{\partial \mathbb{D}}^\theta$, we denote by $x(e) \in \partial \mathbb{D}$ and $y(e) \in \partial \mathbb{D}$ its starting point and endpoint, using the time-orientation inherited from its parent loop $\wp$. One important property, proved in \cite[Lemma 9]{qian2018lawpointprocessbrownian} is the following: conditionally on the endpoints, the trajectories $e \in \mathcal{E}_{\partial \mathbb{D}}^\theta$ are independent Brownian excursions in $\mathbb{D}$ linking $x(e)$ to $y(e)$.

\subsection{Hitting small discs with excursions and bubbles.}

The aim of this section is to prove estimates on the probability that a Brownian excursion in $\mathbb{D}$ linking endpoints $(a,b) \in \partial \mathbb{D}$ and a rooted Brownian bubble hit a small disc $r\mathbb{D}$.

\begin{lemma} \label{lemma:excur_hit_a_b}
Denote $\mathbb{P}_{\mathbb{D}}^{a,b}$ the law of the Brownian excursion in the unit disc with endpoints $a,b \in \partial \mathbb{D}$. We have, for $a \neq b$, as $r \to 0$,
$$
\mathbb{P}_{\mathbb{D}}^{a,b}\left(e \text{ hits } r\mathbb{D} \right) \,=\, \frac{|a-b|^2}{2|\log r|}\big( 1  + O(r) \big).
$$
Here the $O(r )$ error is uniform with respect to the endpoints $a,b \in \partial \mathbb{D}$.
\end{lemma}

\begin{proof} By rotational symmetry, we can set $a = 1$. Let $\mathbb{P}_x$ denote the law of a planar Brownian motion $(B_t)_{t \geq 0}$ started at $x \in \mathbb{D}$, and $\tau_{\partial \mathbb{D}}$ its first exit time from the unit disc. We recall that the density of the exit point $B_{\tau_{\partial \mathbb{D}}}$ on the circle $\partial \mathbb{D}$ exists and is given by the Poisson kernel: for all $\xi \in \partial \mathbb{D}$,
\begin{equation} \label{eq:def_poisson_kernel_int}
H_\mathbb{D}(x, \xi) \dd \xi \,:=\, \mathbb{P}_x\big( B_{\tau_{\partial \mathbb{D}}} \in \dd \xi \big) \,=\, \frac{1}{2\pi} \frac{1-|x|^2}{|x-\xi|^2} \dd \xi.
\end{equation}
Let $\tau_r$ be the first time $(B_t)_{t \geq 0}$ hits $B(0,r)$. By the strong Markov property at time $\tau_r$, we get, for all $x \in \mathbb{D} \setminus B(0,r)$,
$$
\mathbb{P}_x\left(\tau_r<\tau_{\partial \mathbb{D}} \text { and } B_{\tau_{\partial \mathbb{D}}} \in \dd b\right)=\mathbb{E}_x\left[\mathbf{1}_{\left\{\tau_r<\tau_{\partial \mathbb{D}}\right\}} H_{\mathbb{D}}\left(B_{\tau_r}, b\right)\right]\dd b.
$$
By conditioning on the event $\{B_{\tau_{\partial \mathbb{D}}} = b\}$ (which has density $H_\mathbb{D}(x,b)$ under $\mathbb{P}_x$), we get
\begin{equation} \label{eq:p_a_b_poisson}
\mathbb{P}_{\mathbb{D}}^{1,b}\left(e \cap B(0,r) \neq \emptyset\right) \,=\, \lim_{\varepsilon \to 0} \frac{\mathbb{E}_{1-\varepsilon}\big[\mathbf{1}_{\{\tau_r<\tau_{\partial \mathbb{D}}\}} H_{\mathbb{D}}(B_{\tau_r}, b)\big]}{H_{\mathbb{D}}(1-\varepsilon, b)}.
\end{equation}
This identity follows directly from the definition of the Brownian excursion in $\mathbb{D}$ from $a = 1$ to $b$, viewed as the limit as $\varepsilon \to 0$ of a planar Brownian motion started at $1-\varepsilon$ and conditioned to exit $\mathbb{D}$ at $b$. Now, at time $\tau_r$, we have $|B_{\tau_r}| = r$, hence, since $|b| = 1$, almost surely,
$$
|B_{\tau_r}-b|^2  \,=\,  1 - 2\Re(\bar{b} B_{\tau_r}) + r^2 \,=\, 1 + O(r).
$$
Here the $O(r)$ error is uniform with respect to the endpoint $b \in \partial \mathbb{D}$. By Equation~\eqref{eq:def_poisson_kernel_int}, we get, almost surely
$$
H_\mathbb{D}(B_{\tau_r},b) \,=\, \frac{1}{2\pi} \frac{1-r^2}{1+O(r)} \,=\, \frac{1}{2\pi} + O(r).
$$
Since this control is also uniform for all $z = B_{\tau_r} \in r \partial \mathbb{D}$, we can take it out of the expectation and get
$$
\mathbb{E}_{1-\varepsilon}\left[\mathbf{1}_{\left\{\tau_r<\tau_{\partial \mathbb{D}}\right\}} H_{\mathbb{D}}\left(B_{\tau_r}, b\right)\right]\,=\,\Big(\frac{1}{2 \pi}+O(r)\Big) \mathbb{P}_{1-\varepsilon}\left(\tau_r<\tau_{\partial \mathbb{D}}\right).
$$
Denoting $u(z) = \mathbb{P}_z(\tau_r < \tau_{\partial \mathbb{D}})$, the function $u(z)$ is harmonic on the annulus $\mathbb{D} \setminus \overline{r \mathbb{D}}$, its value is $1$ on $r\partial \mathbb{D}$, $0$ on $\partial \mathbb{D}$ and is rotationally invariant, hence $u(z) = \frac{-\log |z|}{|\log r|}$. From this we get
$$
\mathbb{P}_{1-\varepsilon}(\tau_r < \tau_{\partial \mathbb{D}}) \,=\, \frac{-\log (1-\varepsilon)}{|\log r|} \,=\, \frac{\varepsilon + O(\varepsilon^2)}{|\log r|}.
$$
Furthermore, since $u(z)$ does not depend on the endpoints, the error is independent of $b \in \partial \mathbb{D}$. It remains to study the denominator in Equation~\eqref{eq:p_a_b_poisson}. For $x = 1-\varepsilon$, we get
$$
1-|x|^2 = 2\varepsilon + O(\varepsilon^2) \quad \text{and} \quad |x-b|^2 = |1-b|^2 + O(\varepsilon).
$$
By Equation~\eqref{eq:def_poisson_kernel_int}, the denominator reduces to
$$
H_\mathbb{D}(1-\varepsilon, b) \,=\, \frac{1}{2\pi} \frac{2\varepsilon + O(\varepsilon^2)}{|1-b|^2 + O(\varepsilon)} \,=\, \frac{\varepsilon}{\pi|1-b|^2} + O(\varepsilon^2).
$$
Injecting the previous asymptotics into Equation~\eqref{eq:p_a_b_poisson}, we finally get
$$
\mathbb{P}_{\mathbb{D}}^{1, b}(e \cap B(0, r) \neq \emptyset)=\lim _{\varepsilon \rightarrow 0} \frac{\left(\frac{1}{2 \pi}+O(r)\right)\big(\frac{\varepsilon+O(\varepsilon^2)}{|\log r|}\big)}{\frac{\varepsilon}{\pi|1-b|^2}+O\left(\varepsilon^2\right)} \,=\, \frac{|1-b|^2}{2|\log r|}(1+O(r)).
$$
Moreover, as shown above, the $O(r)$ error is uniformly bounded in the endpoint $b \in  \partial \mathbb{D}$. 
\end{proof}

We will also need asymptotics on the mass of rooted Brownian bubbles that hit a small disc. We recall that our normalisation of the Brownian bubble measure is different from the one found in \cite{lawler2003brownianloopsoup} and differs by a factor $\frac{1}{\pi}$ (see Section~\ref{sec:bubbles}).

\begin{lemma} \label{lemma:bubble_hit_ball_disc}
As $r \to 0$, uniformly with respect to $\theta \in [0, 2\pi)$, 
$$\mu_\mathbb{D}^{\mathrm{bub}}(e^{i\theta}) \big(\big\{ \wp \text{ hits } r\mathbb{D} \big \} \big) \,=\, \frac{1}{2\pi} \frac{1}{|\log r|} \big(1 + O(r)).$$
\end{lemma}

\begin{proof}
Let $(a,b) \in \partial \mathbb{D}$, we recall the definition of the boundary Poisson kernel (see \cite{lawler2003brownianloopsoup}) as  $H_\mathbb{D}(a,b) \,=\, \frac{1}{\pi}|b-a|^{-2}$ and define the boundary-to-boundary excursion measure  $\mu_\mathbb{D}(a,b) \,=\, H_\mathbb{D}(a,b) \mathbb{P}_\mathbb{D}^{a,b}$. Then, the Brownian bubble measure rooted at $a \in \partial \mathbb{D}$ can be defined as 
$$
\mu_\mathbb{D}^{\text{bub}}(a) \,=\, \lim_{b \to a} \mu_\mathbb{D}(a,b).
$$
Again, we recall that our normalisation is different from the standard one in literature. For our event of interest, together with Lemma~\ref{lemma:excur_hit_a_b} with $a = e^{i\theta}$, we get
\begin{equation} \label{eq:hit_ball_bubble_error}
\mu_\mathbb{D}^{\mathrm{bub}}(e^{i\theta}) \big(\big\{ \wp \text{ hits } r\mathbb{D} \big \} \big) \,=\, \lim_{b \to a} H_\mathbb{D}(a,b) \mathbb{P}_\mathbb{D}^{a,b}\big(\big\{ \wp \text{ hits } r\mathbb{D} \big \} \big) \,=\, \frac{1}{2\pi} \frac{1}{|\log r|} \big(1 + O(r)).
\end{equation}
As the error is uniform over the endpoints $(a,b)\in \partial \mathbb{D}$ in Lemma~\ref{lemma:excur_hit_a_b}, the error in Equation~\eqref{eq:hit_ball_bubble_error} becomes uniform in $\theta \in [0,2\pi)$.
\end{proof}

\subsection{Random coverings of $(0,\infty)$.} \label{sec:random_coverings}

To study the longitudinal connectivity of the Brownian loop soup, we reduce the two-dimensional problem to two specific one-dimensional random coverings. Random coverings of the real line by intervals have been extensively studied since the seminal work of Shepp \cite{Shepp1972}. They are closely related to the theory of random cutouts, perfect regenerative sets, and subordinators (see e.g., \cite{Mandelbrot1972, Kingman1972, Bertoin1999}).
We consider $\theta > 0$ and a Poisson point process $\Lambda = \{a_i, \ell_i\}_{i \geq 1}$ on $\mathbb{R}_+ \times (0,\infty)$ with intensity measure
$$
\theta \, \dd a \otimes \mu(\dd\ell)
$$
We identify each pair $(a,\ell) \in \Lambda$ as the open interval $(a,a+\ell)$, and call $a$ the \emph{left endpoint} of the interval and $\ell$ its \emph{length}. We are interested in the connectivity properties of the random covered set $\bigcup_{(a,\ell) \in \Lambda} (a,a+\ell)$. In our context, the relevant random coverings are those satisfying the following Hypothesis. \medskip

\noindent\textbf{Hypothesis (H)}: \textit{The random covering $\Lambda$ is said to satisfy Hypothesis (H) if it is a Poisson point process of intensity $\theta \dd a \otimes \mu(\dd \ell)$, where $\mu$ is a finite measure on $(0,\infty)$, and, for some constant $c \in \mathbb{R}$, as $x \to \infty$,}
$$
\int_1^x \mu\big((v, \infty) \big) \dd v \,=\, \log x + c + o(1).
$$

The canonical example satisfying Hypothesis~(H) is the reference scale invariant measure with UV cutoff $\mu(\dd\ell) = \ell^{-2} \mathbf{1}_{\{\ell > 1\}} \dd\ell$ (which corresponds exactly to the endpoints of the 1D truncated Brownian loop soup, see \cite[Proposition 3.18]{Lupu_2018}). In general, Hypothesis~(H) means that $\mu$ behaves like this reference measure: the first condition filters out small intervals to avoid local covering by microscopic intervals, while the second condition ensures that $\mu$ has the same tail at infinity. The following proposition is a macroscopic universality result. It states that the percolation properties of the model (namely the critical threshold and the crossing probabilities) depend only on the macroscopic tail of $\mu$, and do not depend on its microscopic details.

\begin{proposition} \label{prop:random_coverings}  Let $y > 0$. Under Hypothesis (H), we have the following dichotomy
\begin{itemize}[itemsep=2pt, topsep=3pt, parsep=0pt]
\item if $\theta \in (0,1)$, every cluster is a.s. bounded and $\mathbb{P}([y,x] \text{ is covered by } \Lambda) \sim K(y) x^{-1+\theta}$.
\item if $\theta = 1$, every cluster is a.s. bounded and $\mathbb{P}([y,x] \text{ is covered by } \Lambda) \sim K^\prime(y) (\log x)^{-1}$.
\item if $\theta > 1$, almost surely, there exists a unique unbounded cluster.
\end{itemize}
All the asymptotics are taken when $x \rightarrow \infty$. The constants involved all depend on $y, \theta$ and $\mu$. Furthermore, for all $\delta > 0$, for all $y > \delta$, $K(y) \asymp y^{1-\theta}$ (with some constants depending only on $\delta, \theta, \mu$).
\end{proposition}

\begin{proof}
This proof is entirely based on \cite{Bertoin1999} (Chapter 7). Denote $\Lambda = \{(a_i, b_i)\}_{i \in I}$ and define the uncovered set $\mathcal{R}$ as 
\begin{equation} \label{eq:def_uncovered}
\mathcal{R} = [0,\infty) \setminus \bigcup_{i \in I} (a_i, b_i).    
\end{equation}
The first condition of Hypothesis (H) ensures that $\mu((v,\infty))$ is integrable around $v \rightarrow 0^+$. Denote $M = \int_0^1 \mu((v,\infty))\dd v < \infty$, we have, 
$$
\int_0^1 \exp \left\{ \theta \int_x^1 \mu((v, \infty)) \dd v \right\} \dd x \,\leq\, \int_0^1  e^{\theta M}  \dd x \,<\, \infty.
$$
By \cite{Bertoin1999} (Theorem 7.2), we have $\mathcal{R} \neq \{0\}$ almost surely, and $\mathcal{R}$ is a perfect regenerative set. Hence, by \cite{Bertoin1999} (Lemma 7.1), there exists a subordinator $(\sigma_t)_{t \geq 0}$ such that $\mathcal{R}$ is almost surely the range of $\sigma$. Mathematically, denoting $\zeta$ the lifetime of $\sigma$,
$$
\mathcal{R} = \overline{\{\sigma_t: 0 \leq t < \zeta \}} \quad \text{ a.s.}
$$
Further applying \cite{Bertoin1999} (Theorem 7.2), we also get that the renewal measure $U(\dd x)$ of $\sigma$ is absolutely continuous with respect to Lebesgue measure $\dd x$ with continuous density $u(x)$ given by
\begin{equation} \label{eq:renewal_measure}
U(\dd x) \,=\, u(x)\dd x \,=\, c \exp \left\{ \theta \int_x^1 \mu((v,\infty)) \dd v \right\} \dd x,
\end{equation}
for some positive constant $c > 0$. Note that since $\mu$ is a finite measure (by the first condition of (H)), we have $u(0^+) > 0$, it follows that the subordinator $\sigma$ has a positive drift $d = u(0^+)^{-1} > 0$ and a finite Lévy measure, i.e., $\sigma$ is a (killed) compound Poisson process with positive drift. Furthermore, the second condition of Hypothesis (H) is 
\begin{equation} \label{eq:asympt_int_mu}
\int_1^x \mu((v,\infty)) \dd v \,=\,  \log x + c^\prime + o(1).
\end{equation}
as $x \rightarrow \infty$. Now, by \cite{Bertoin1999} (Corollary 7.4), if $\int_1^\infty \exp \{ - \theta \int_1^x \mu((v,\infty)) \dd v\} \dd x = \infty$ then $\mathcal{R}$ is unbounded almost surely, otherwise $\mathcal{R}$ is bounded almost surely. In our case, by Equation~\eqref{eq:asympt_int_mu}, we have
$$
\exp \left\{ - \theta \int_1^x \mu((v,\infty)) \dd v \right\} \,=\, e^{-\theta \log x - \theta c^\prime + o(1)} \,\sim\, cx^{-\theta}.
$$
By noticing that $\Lambda$ contains an unbounded connected component if and only if $\mathcal{R}$ is bounded, we get that for all $\theta \in (0,1]$, almost surely all connected components of $\Lambda$ are bounded, and when $\theta > 1$, almost surely, there exists an unbounded connected component. \newline

\noindent We suppose here that $\theta \in (0,1]$. For any $y > 0$, denote $g_y = \sup \{ u \leq y : u \in \mathcal{R} \}$ the last passage time before $y$ and $D_y = \inf \{u > y : u \in \mathcal{R}\}$ the first passage time after $y$. By definition of $\mathcal{R}$, we get
$$
\{[y,x] \text{ is covered by } \Lambda\} \,=\, \{\mathcal{R} \cap (y,x) = \emptyset\} \,=\, \{D_y > x\}.$$
But by \cite{Bertoin1999} (Lemma 1.10), the joint law of $g_y$ and $D_y-g_y$ is given by
$$
\mathbb{P}(g_y \in \dd a, D_y-g_y \in \dd b) = \Pi(\dd b) U(\dd a).
$$
Where $\Pi$ is the Lévy measure of the subordinator $\sigma$ and $U$ is the renewal measure defined in Equation~\eqref{eq:renewal_measure}. By that same equation, we also have $U(\dd a) = u(a) \dd a$, hence
$$
\mathbb{P}\big(D_y > x\big) \,=\, \int_0^y \mathbb{P}\big(g_y \in \dd a, D_y-g_y > x-a\big) \,=\, \int_0^y \overline{\Pi}(x-a) u(a) \dd a,
$$
where $\overline{\Pi}(v) = \Pi((v, \infty))$. As we will see later in the proof, $\overline{\Pi}$ is a regularly varying function of index $1-\theta$ (and a slowly varying function when $\theta = 1$), this implies in particular that $\overline{\Pi}(x-a) \sim \overline{\Pi}(x)$ uniformly on all compact sets (see \cite{Bingham_Goldie_Teugels_1987}, Theorem 1.5.2). Hence, as $x \rightarrow \infty$,
\begin{equation} \label{eq:coverage_Pi}
\mathbb{P} \big( [y,x] \text{ is covered by } \Lambda \big) \,\sim\, \overline{\Pi}(x) U(y),
\end{equation}
where $U(y) = U((0,y)) < \infty$. It thus only remains to understand the behavior of $\overline{\Pi}(x)$ as $x \rightarrow \infty$. \newline

Suppose $\theta \in (0,1]$. By Equation~\eqref{eq:asympt_int_mu} and Equation~\eqref{eq:renewal_measure}, we have $u(x) \sim C_1 x^{-\theta}$ for some positive constant $C_1 > 0$. This implies that $U(v) = \int_0^v u(x) \dd x \sim C_2 v^{1-\theta}$ if $\theta \in (0,1)$ and $U(v) \sim C_5 \log v$ if $\theta = 1$, for some positive constants $C_2, C_5$. Hence $U$ is regularly varying as $v \rightarrow \infty$ of index $\alpha = 1-\theta > 0$ if $\theta \in (0,1)$ and slowly varying if $\theta = 1$. Denoting by $\Psi(\lambda)$ the Laplace exponent of $\sigma$, by the Tauberian theorem \cite{Bingham_Goldie_Teugels_1987} (Theorem 1.7.1), we get
\begin{align*}
\frac{1}{\Psi(\lambda)} \,\sim\, \Gamma(1+\alpha) U(1/\lambda) \,&\sim\, \Gamma(2-\theta) C_2 \lambda^{\theta-1} \quad\quad \text{if } \theta \in (0,1) \\
&\sim\, C_6 \log(1/\lambda) \quad\quad\quad\quad \text{if }  \theta =1,
\end{align*}
as $\lambda \to 0^+$. This implies that $\Psi$ is regularly varying of index $1-\theta$ (and slowly varying if $\theta = 1$) near $\lambda = 0^+$. Hence, by \cite{Bertoin1999} (Proposition 1.5), we get
\begin{align*}
\overline{\Pi}(x) \,\sim\, \frac{1}{\Gamma(1-\alpha)} \Psi(1/x) \,&\sim\, C_3 x^{-1 + \theta} \quad\quad\quad\quad \text{if } \theta \in (0,1) \\
&\sim\, C_7 (\log x)^{-1} \quad\quad\quad \text{if }  \theta =1.
\end{align*}
With $C_3 = (C_2 \Gamma(2-\theta) \Gamma(\theta))^{-1}$. Injecting this into Equation~\eqref{eq:coverage_Pi} yields the desired result. 

Finally, for $\theta \in (0,1)$, Equation~\eqref{eq:coverage_Pi} identifies the constant as $K(y) = C_3 U(y)$. By \eqref{eq:renewal_measure}, $U(y) \sim C_2 y^{1-\theta}$ as $y \to \infty$, it follows immediately that $K(y) \asymp y^{1-\theta}$ uniformly for all $y > \delta$.
\newline
\end{proof}

Suppose that $\Lambda$ satisfies Hypothesis (H). Since $\mu$ is a finite measure, the intensity measure $\theta da \otimes \mu(\mathrm{d}\ell)$ is finite on any compact subset of $(0, \infty) \times (0, \infty)$. Hence, almost surely, only finitely many intervals of $\Lambda$ intersect any given compact set and the uncovered set $\mathcal{R} = [0, \infty) \setminus \bigcup_{(a,\ell)\in\Lambda} (a, a+\ell)$ defined in~\eqref{eq:def_uncovered} consists of locally finite disjoint closed components. For any $0 < y < x$, the restriction $\mathcal{R} \cap (y,x)$ contains a finite number $N(y,x)$ of connected components. The proof of Theorem~\ref{th:crossing_exponent} relies on the following bound.

\begin{lemma} \label{lemma:generating_function_Nst}
    Let $\theta \in (0,1)$, let $q \in (0,1)$ and let $y > 0$ be fixed. For all $x > y$, we have
    $$
    \mathbb{E}[ q^{N(y,x)} ] \,\asymp\, x^{-1+\theta}.
    $$
\end{lemma}

\begin{proof}
We begin with the lower bound. Since $q \in (0,1)$, all terms in the expectation are non-negative. Keeping only the term $k = 0$, we have
\[
\mathbb{E}[q^{N(y,x)}] \,=\, \sum_{k=0}^\infty q^k \mathbb{P}(N(y,x) = k) \,\ge\, \mathbb{P}(N(y,x) = 0).
\]
The event $\{N(y,x) = 0\}$ means that $\mathcal{R} \cap (y,x) = \emptyset$, which implies that the interval $(y,x)$ is entirely covered by $\Lambda$. By Proposition~\ref{prop:random_coverings}, we obtain
\[
\mathbb{P}(N(y,x) = 0) \,=\, \mathbb{P}\big((y,x) \text{ is covered by } \Lambda\big) \,\asymp\, x^{-1+\theta},
\]
which yields the desired lower bound. We now turn to the upper bound. Since $\mu$ is a finite measure, the intervals of $\Lambda$ arrive at a finite rate $\lambda := \theta \mu((0,\infty)) < \infty$, so almost surely, we can write $\mathcal{R} = \bigcup_{i \geq 1} [s_i, t_i]$ with $0 = s_1 < t_1 < s_2 < t_2 < \dots$. In the subordinator representation of $\mathcal{R}$ from the proof of Proposition~\ref{prop:random_coverings}, the times $t_i - s_i$ are i.i.d.\ $\mathrm{Exp}(\lambda)$, while the covered clusters $s_{i+1} - t_i$ correspond to the jumps of the subordinator $\sigma$, which arrive at finite rate $\Pi((0,\infty)) \in (0,\infty)$ and satisfy $\overline{\Pi}(u) \sim C_3 u^{-1+\theta}$ as $u \to \infty$. Hence $\mathbb{P}(s_{i+1} - t_i > u) = \overline{\Pi}(u)/\Pi((0,\infty))$, and we get, for all $u > 0$,
\[
\mathbb{P}(s_2 - s_1 \ge u) \,\le\, C u^{-1+\theta},
\]
for some constant $C > 0$. Now let $D_y = \inf\{u > y : u \in \mathcal{R}\}$ be the first passage time after $y$, and set $K = \inf\{i \ge 1 : t_i > y\}$. Because $\mathcal{R} \cap [0,u]$ is measurable with respect to $\mathcal{F}_u := \sigma(\Lambda \cap ([0,u] \times (0,\infty)))$, both $D_y$ and $(s_{K+i})_{i \ge 1}$ are $\mathcal{F}_u$-stopping times. Define $L_{-1} := D_y - y$, $L_0 := s_{K+1} - D_y$, and $L_i := s_{K+i+1} - s_{K+i}$ for $i \ge 1$. Since almost surely $D_y = y$ when $y \in \mathcal{R}$ and $D_y = s_K$ when $y \notin \mathcal{R}$, the strong Markov property of $\Lambda$ at $D_y$ ensures that $(L_i)_{i \ge 0}$ are i.i.d.\ copies of $s_2 - s_1$. Moreover, by the renewal formula in the proof of Proposition~\ref{prop:random_coverings}, for all $u > 0$,
\[
\mathbb{P}(L_{-1} > u) \,=\, \mathbb{P}(D_y > y + u) \,=\, \int_0^y \overline{\Pi}(y + u - a)\,u(a)\,\mathrm{d}a \,\le\, U(y)\,\overline{\Pi}(u) \,\le\, C_y u^{-1+\theta}.
\]
Fix $k \ge 0$. On the event $\{N(y,x) = k\}$, almost surely, $(y,x)$ is covered by $(y, D_y) \cup [D_y, s_{K+1}) \cup \bigcup_{i=1}^{k-1} [s_{K+i}, s_{K+i+1})$ (just $(y, D_y)$ when $k = 0$), which gives $\sum_{i=-1}^{k-1} L_i \ge x - y$ almost surely. A union bound over these $k+1$ pieces then yields
\[
\mathbb{P}(N(y,x) = k) \,\le\, \sum_{i=-1}^{k-1} \mathbb{P}\!\left(L_i \ge \frac{x-y}{k+1}\right) \,\le\, C'' (k+1)^{2-\theta} (x - y)^{-1+\theta}.
\]
Summing over $k \ge 0$, we obtain
\[
\mathbb{E}[q^{N(y,x)}] \,=\, \sum_{k \geq 0} q^k \mathbb{P}(N(y,x) = k) \,\leq\, C'' (x-y)^{-1+\theta} \sum_{k \geq 0} (k+1)^{2-\theta} q^k.
\]
That sum on the right is a finite constant $S(\theta, q) < \infty$. Injecting this together with the fact that for all $x \geq 2y$ we have $x-y \geq x/2$, concludes the proof. The case $y<x < 2y$ is trivial: since $q \in (0,1)$, we have $\mathbb{E}[q^{N(y,x)}] \leq 1 \leq (2y)^{1-\theta} x^{-1+\theta}$.
\end{proof}

\section{Phase transition for the truncated loop soup.} \label{sec:phase_transition}

Let $\mathcal{L}_{\mathfrak{C}}^{\theta}$ be a Brownian loop soup on the cylinder $\mathfrak{C} = (0, \infty) \times \mathbb{S}^1$. Let $\tau > 0$ and define the truncated loop soup with a duration cutoff $\LC$ as the subset of $\mathcal{L}_{\mathfrak{C}}^{\theta}$ consisting of loops of duration (measured from the cylinder) at least $\tau$. In this section, we prove Theorem~\ref{th:phase_transition_loopsoup}.

\subsection{Sandwiching the loop soup between two random coverings.}

Let $\bpi: \mathfrak{C} \to \mathbb{R}_+$ be the longitudinal projection on the cylinder defined by $\bpi(x, e^{i\theta}) = x$. For any loop $\wp$ on $\mathfrak{C}$, its longitudinal span is the open interval $(a(\wp), b(\wp))$, where $a(\wp) := \min \bpi(\wp)$ and $b(\wp) := \max \bpi(\wp)$. We define the 1D \textit{projection} of $\wp$ as the open interval $\big(a(\wp), b(\wp)\big)$. We call $a(\wp)$ the left endpoint of $\wp$ and define the projection length $\ell(\wp) = b(\wp)-a(\wp)$. Following the framework of Section~\ref{sec:random_coverings}, we denote the collection of all projections as:
\begin{equation}
\Lp \,:=\, \left\{ a(\wp), \ell(\wp) \right\}_{\wp \in \LC}.
\end{equation}
For any $(a,\ell) \in \Lp$ we call $(a,a+\ell)$ a \emph{projection interval} of $\Lp$. By definition, any connected cluster of loops on the cylinder projects onto a connected component of projection intervals on the real line. Consequently, the absence of macroscopic percolation in the 1D model $\Lp$ implies the absence of percolation in the 2D loop soup $\LC$. However, this provides only a necessary condition for 2D connectivity, as two loops can easily have overlapping projections without intersecting on the cylinder. \medskip

To obtain a sufficient condition for 2D intersection, we introduce a second, more restrictive covering of $(0, \infty)$. We focus exclusively on separating loops, as they act as topological barriers around the cylinder. We recall (see Section~\ref{sec:prelim_loops}) that for any separating loop $\wp$ on the cylinder $\mathfrak{C}$, its inner boundary $\partial^{\text{i}}\wp$ is defined as the boundary of the connected component of $\mathfrak{C} \setminus \wp$ containing $\{\infty\} \times \mathbb{S}^1$. Topologically, $\partial^{\text{i}}\wp$ is almost surely a simple separating loop. We denote its minimal longitudinal depth by $c(\wp) := \min \bpi(\partial^{\text{i}}\wp)$ and define the 1D \emph{core interval} of $\wp$ as the open interval $(a(\wp), c(\wp))$. We still call $a(\wp)$ the left endpoint of $\wp$ and define the \emph{core length} as $w(\wp) = c(\wp)-a(\wp)$. We define the collection of \emph{core intervals} $\Lm$ as:
\begin{equation}
\Lm \,:=\, \left\{ a(\wp), w(\wp) \right\}_{\wp \in \LC : \wp \text{ is separating}}
\end{equation}
The introduction of $\Lm$ relies on a crucial topological property: if the core intervals of two separating loops overlap in $\mathbb{R}_+$, the loops must necessarily intersect on the cylinder $\mathfrak{C}$. The following deterministic lemma formalizes this coupling between the 2D loop clusters and our two 1D random coverings.

\begin{figure}[h]
    \centering
    \includegraphics[width=14cm]{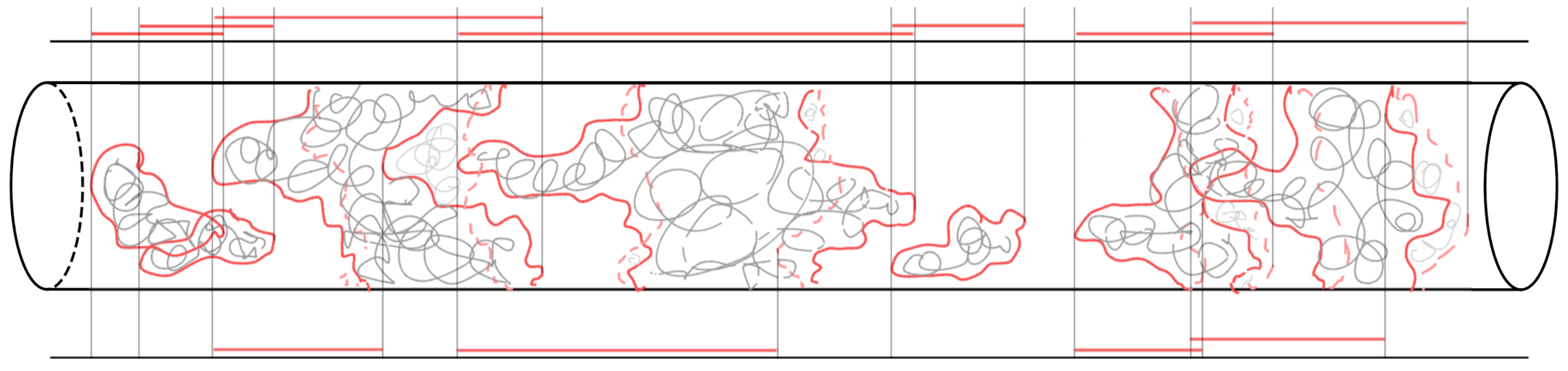}
    
    \caption{A Brownian loop soup $\LC$ composed of 7 loops (4 of which are separating) and 5 clusters on a slice of the cylinder $\mathfrak{C}$ together with (on top), the corresponding projection intervals $\Lp$ and (on the bottom), the core intervals $\Lm$. The projection intervals merge into 2 clusters, while the core intervals merge into 3. The non-separating loops of $\LC$ do not appear in $\Lm$.}
    \label{fig:coverings}
\end{figure}

\begin{lemma} \label{lemma:coupling_coverings} For any $0 \leq y<x \leq \infty$, covering the interval $(y, x)$ by core intervals of $\Lm$ is a sufficient condition for a cluster of $\LC$ to cross from $\{y\} \times \mathbb{S}^1$ to $\{x\} \times \mathbb{S}^1$ while covering $(y, x)$ by projection intervals of $\Lp$ is a necessary condition.
\end{lemma}

\begin{proof}
The necessary condition is immediate: if two loops $\wp, \wp^{\prime}$ intersect in $\mathfrak{C}$, their projections satisfy $\boldsymbol{\pi}(\wp) \cap \boldsymbol{\pi}\left(\wp^{\prime}\right) \neq \emptyset$, so their projection intervals in $\Lp$ overlap. For the sufficient condition, it is sufficient to prove that overlapping core intervals in $\Lm$ force their corresponding loops to intersect. We proceed by contraposition: let $\wp, \wp^{\prime} \in \mathcal{L}_{\mathfrak{C}}^{\theta, \tau}$ be two separating loops such that $\wp \cap \wp^{\prime}=\emptyset$. Without loss of generality, assume $\wp^{\prime}$ lies in the unbounded component of $\mathfrak{C} \backslash \wp$ containing $\{\infty\} \times \mathbb{S}^1$. Since the inner boundary $\partial^{\mathrm{i}} \wp$ separates this component from $\{0\} \times \mathbb{S}^1$, the loop $\wp^{\prime}$ lies strictly to the right of $\partial^{\mathrm{i}} \wp$. In particular,
$$
a(\wp^{\prime})=\min \boldsymbol{\pi}(\wp^{\prime})>\min \boldsymbol{\pi}\big(\partial^{\mathrm{i}} \wp\big)=c(\wp),
$$
which implies that $(a(\wp), c(\wp)) \cap\left(a\left(\wp^{\prime}\right), c\left(\wp^{\prime}\right)\right)=\emptyset$. By contraposition, overlapping core intervals in $\underline{\Lambda}^{\theta, \tau}$ force $\wp$ and $\wp^{\prime}$ to intersect. Hence, any chain of core intervals covering $(y, x)$ must generate a 2D connected cluster of separating loops crossing from $\{y\} \times \mathbb{S}^1$ to $\{x\} \times \mathbb{S}^1$.
\end{proof}

The previous lemma manages to sandwich the two-dimensional percolation on the cylinder between two random coverings of $(0,\infty)$. In order to prove Theorem~\ref{th:phase_transition_loopsoup}, it remains to prove that $\Lp$ and $\Lm$ belong to the same universality class, which, by Proposition~\ref{prop:random_coverings}, reduces to showing that both their intensity measures satisfy Hypothesis (H).

\begin{lemma} \label{lem:ppp_naive_projection}
$\Lp$ is a Poisson point process of intensity measure given by
$$
\theta \, \dd a \otimes \mu_\mathfrak{C}^{\mathrm{bub}, \tau}\big(\ell(\wp) \in \dd \ell\big).
$$
Furthermore, the length measure $\ell_* \mu_\mathfrak{C}^{\mathrm{bub}, \tau}$ satisfies Hypothesis (H) from Section~\ref{sec:random_coverings}.
\end{lemma}

\begin{proof} Denote $(a\otimes \ell) : \wp \mapsto (a(\wp), \ell(\wp))$ the map that associates to each loop $\wp$ its left endpoint and its projection length. Since $\Lp = (a \otimes \ell)(\LC)$ and $a \otimes \ell$ is a measurable functional, the mapping theorem ensures that $\Lp$ is a Poisson point process of intensity $\theta (a \otimes \ell)_* \mu_{\mathfrak{C}}^{\text{loop}, \tau}$. Furthermore, by Section~\ref{sec:prelim_loops} and Equation~\eqref{eq:loop_bubble_decomp}, we have
$$\mu_{\mathfrak{C}}^{\mathrm{loop}, \tau} =  \int_0^\infty \dd a \, \tau_a \circ \mu_{\mathfrak{C}}^{\mathrm{bub}, \tau}.$$
Here $\mu_{\mathfrak{C}}^{\mathrm{bub}, \tau}$ denotes the measure $\mu_{\mathfrak{C}}^{\mathrm{bub}}$ restricted to the loops of duration at least $\tau$. Since each loop $\wp$ under $\tau_a \circ \mu_{\mathfrak{C}}^{\mathrm{bub}, \tau}$, has its minimal longitudinal position $a(\wp) = \min \bpi(\wp)$ at exactly $a$, and $\ell = b-a$, we have
$$
(a \otimes \ell)_* \mu_{\mathfrak{C}}^{\text{loop}, \tau}(\dd a, \dd \ell) \,=\, \dd a \, \mu_{\mathfrak{C}}^{\text{bub}, \tau}\big( \ell(\wp) \in \dd\ell \big).
$$
This proves the first identity of the lemma. Denoting by $(X_t)_{t \geq 0}$ the longitudinal component of $\wp$ and $\zeta(X)$ its lifetime and integrating the Brownian bubble measure decomposition from Equation~\eqref{eq:bubble_explicit}, we get
$$
\mu_\mathfrak{C}^{\text{bub},\tau}(\dd X) = 2\pi p_{\mathbb{S}^1}(\zeta) \mathbf{1}_{\{\zeta > \tau\}} \mathbf{n}(\dd X).
$$
Since $\ell(\wp)$ is exactly $\sup_{t \geq 0} X_t$, we get 
\begin{equation} \label{eq:mu_bub_ell_n}
\mu_{\mathfrak{C}}^{\text{bub}, \tau}\big( \ell(\wp) \in \dd\ell \big) \,=\, \mathbf{n}\left(2\pi p_{\mathbb{S}^1}(\zeta) \mathbf{1}_{\{\zeta > \tau\}}  \mathbf{1}_{\{\sup_{t \geq 0} X_t \in \dd\ell\}}\right).
\end{equation}
For the rest of this proof, for the sake of readability, we will denote $\mu(\dd \ell) = \mu_{\mathfrak{C}}^{\text{bub}, \tau}\big( \ell(\wp) \in \dd\ell \big)$. To conclude the proof of the Lemma, it remains to prove that $\mu$ satisfies Hypothesis (H). We start with the first condition. We have
$$
\mu((0,\infty)) \,=\, \mathbf{n} \left( 2\pi p_{\mathbb{S}^1}(\zeta) \mathbf{1}_{\{\zeta > \tau\}} \right) \,\leq\, C_\tau \mathbf{n} \big( \zeta > \tau \big) \,<\, \infty.
$$
We have used the fact that the heat kernel $2\pi p_{\mathbb{S}^1}(t) \leq C_\tau < \infty$ on the circle $\mathbb{S}^1$ is uniformly bounded for $t > \tau > 0$, and the fact that $\mathbf{n}$ is locally finite (see for example \cite{revuz1999continuous}). Hence $\mu$ is a finite measure, which concludes the proof of the first condition. Let us now prove the second condition. Denote $M = \sup_{t \geq 0} X_t$. Under Itô's excursion measure, it is well known that $\mathbf{n}(M > u) = \frac{1}{u}$ for all $u > 0$. In what follows, we will establish the following asymptotic
\begin{equation} \label{eq:pointwise_mu_tail}
\mu((u, \infty)) \,=\, \frac{1}{u} + O(e^{-cu}) \quad \text{as } u \to \infty,
\end{equation}
for some constant $c > 0$. By Equation~\eqref{eq:mu_bub_ell_n} and the triangle inequality, we have for all $u \geq 1$:
\begin{equation} \label{eq:decomp_error_mu}
\Big|\mu((u, \infty)) - \frac{1}{u}\Big| \leq \mathbf{n}\left(\left|2 \pi p_{\mathbb{S}^1}(\zeta)-1\right| \mathbf{1}_{\{\zeta>\tau\}} \mathbf{1}_{\{M>u\}}\right)+\mathbf{n}\left(\mathbf{1}_{\{\zeta \leq \tau\}} \mathbf{1}_{\{M>u\}}\right) .
\end{equation}
By Lemma~\ref{lemma:ito_max_tau}, the second term satisfies $\mathbf{n}(M > u \text{ and } \zeta \leq \tau) \leq C'_\tau e^{-cu^2} = O(e^{-cu})$. It remains to bound the first term of Equation~\eqref{eq:decomp_error_mu}. Intuitively, as the lifetime $\zeta$ grows large, the Brownian motion on the cylinder mixes rapidly, and the heat kernel converges to the uniform measure $\frac{1}{2\pi} \dd \varphi$ on $\mathbb{S}^1$. More precisely, by Lemma~\ref{lemma:_heat_kernel_bound}, there exist positive constants $C', c'$ such that for all $\zeta > \tau$,
$$
|2\pi p_{\mathbb{S}^1}(\zeta)-1| \,\leq\, C'e^{-c' \zeta}.
$$
By disintegrating Itô's measure $\mathbf{n}$ with respect to the lifetime $\zeta$, we get
$$
\mathbf{n}\big(e^{-c' \zeta} \mathbf{1}_{\{\zeta>\tau\}} \mathbf{1}_{\{M>u\}}\big) \,=\, \int_\tau^{\infty} \frac{\dd s}{\sqrt{2 \pi s^3}} e^{-c' s} \, \mathbf{n}_s(M>u),
$$
where $\mathbf{n}_s$ denotes the probability law of a Brownian excursion conditioned on $\{\zeta(X) = s\}$. To bound this integral, we split the domain into $s > u$ and $\tau < s \leq u$. On the one hand, for $s > u$, we simply bound the probability $\mathbf{n}_s(M > u)$ by $1$:
$$
\int_u^\infty \frac{\dd s}{\sqrt{2\pi s^3}} e^{-c' s} \,\leq\, \frac{1}{\sqrt{2\pi u^3}} \int_u^\infty e^{-c' s} \dd s \,=\, O(e^{-c' u}).
$$
On the other hand, for $\tau < s \leq u$, by Lemma~\ref{lemma:excur_max_gaussian}, together with Brownian scaling, we have $\mathbf{n}_s(M > u) \leq C_1 e^{-u^2/s}$. Since $u/s \geq 1$, we have $u^2/s \geq u$, which yields $\mathbf{n}_s(M > u) \leq C_1 e^{-u}$. Hence,
$$
\int_\tau^u \frac{\dd s}{\sqrt{2\pi s^3}} e^{-c' s} \, \mathbf{n}_s(M>u) \,\leq\, C_1 e^{-u} \int_\tau^\infty \frac{\dd s}{\sqrt{2\pi s^3}} \,=\, O(e^{-u}).
$$
Combining both regimes, we obtain the estimate~\eqref{eq:pointwise_mu_tail}. Finally, integrating this estimate from $u = 1$ to $x$, we conclude that, as $x \to \infty$,
$$
\int_1^x \mu((u, \infty)) \dd u \,=\, \int_1^x \Big( \frac{1}{u} + O(e^{-cu}) \Big) \dd u \,=\, \log x + c'' + o(1)
$$
where $c'' := \int_1^\infty (\mu((u,\infty)) - \frac{1}{u}) \dd u < \infty$, which completes the proof.
\end{proof}

Lemma~\ref{lem:ppp_naive_projection} together with Proposition~\ref{prop:random_coverings} implies in particular that $\LC$ does not contain any unbounded connected components of loops on the cylinder when $\theta = 1$. Otherwise, their projections $\Lp$ would overlap and form an unbounded connected component on $(0, \infty)$. Hence, the truncated loop soup on the cylinder does not percolate at the critical intensity, and the phase transition is continuous.

\subsection{The core intervals and the boundary span.}

In this section, we prove that the core intervals $\Lm$ satisfy Hypothesis (H), and deduce from this Theorem~\ref{th:phase_transition_loopsoup}. We first recall that $\Lm$ is only defined on separating loops and also recall the definition of the core length $w(\wp)$ of any separating loop $\wp$ as $w(\wp) = c(\wp)-a(\wp)$. An important quantity in the following analysis is the \emph{boundary span} $\delta(\wp)$ that we define as 
\begin{equation} \label{eq:def_delta}
\delta(\wp) \,=\, \ell(\wp) - w(\wp)
\end{equation}
The term boundary span comes from the fact that $\delta(\wp)$ is also equal to $\max \bpi(\partial^{\text{i}} \wp) - \min \bpi(\partial^{\text{i}} \wp)$. The following lemma states that, under the bubble measure $\mu_\mathfrak{C}^{\text{bub}, \tau}$, the boundary span has exponentially light tails.

\begin{lemma} \label{lemma:delta_exp} There exists a positive constant $C$ (depending on $\tau$) such that, for all $u > 0$,
$$
\mu_\mathfrak{C}^{\mathrm{bub}, \tau}\left( \big \{ \delta(\wp) > u \text{ and } \wp \text{ is separating} \big\} \right) \,\leq\, Ce^{- \frac{2}{3} u}.
$$
\end{lemma}

\begin{proof}
Let $\wp$ be a separating loop anchored at $x = 0$. By flipping the cylinder at the exact longitudinal point where $\wp$ achieves its maximum (i.e., applying the map $\psi_\wp : (x,\varphi) \mapsto (\ell(\wp)-x, \varphi)$), we get that the outer boundary $\partial^{\text{o}}\wp$ is equal in distribution to the reflected inner boundary $\psi_\wp(\partial^{\text{i}} \wp)$ (this comes from Vervaat's transform \cite{Ver79}, Equation~\eqref{eq:bubble_explicit} and the fact that, for an excursion $e$, re-rooting the time at the first time $e$ achieves its maximum and flipping the space by $x \to (\sup e - x)$ preserves the law of $e$). This implies in turn that, under $\mu_\mathfrak{C}^{\text{bub}, \tau}$, and conditionally on the fact that $\wp$ is separating, the variable $\delta(\wp)$ is equal in distribution to $\max \bpi (\partial^{\text{o}} \wp)$. See Figure~\ref{fig:mirror}. 

\begin{figure}[h]
    \centering
    \includegraphics[width=12.5cm]{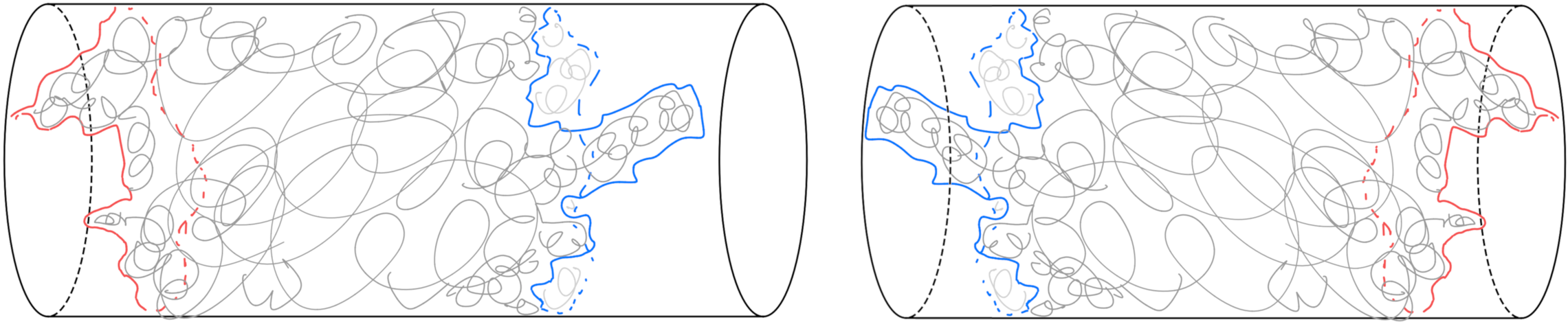}
    
    \caption{Left: the Brownian bubble $\wp$, right: the flipped bubble $\psi_\wp(\wp)$. The map $\psi_\wp$ is a reflection across $\{\ell/2\} \times \mathbb{S}^1$, the inner boundary (in blue) becomes the outer boundary (in red) and the law of $\wp$ remains unchanged.}
    \label{fig:mirror}
\end{figure}
\noindent Secondly, denoting by $f : re^{i\theta} \mapsto (|\log r|, e^{i\theta})$ the conformal map from $\mathbb{D} \setminus \{0\}$ to the cylinder, we know that $\mu_\mathfrak{C}^{\text{bub}} = f_* \mu_\mathbb{D}^{\text{bub}}$, hence
\begin{align*}
&\mu_\mathfrak{C}^{\mathrm{bub}}\left( \big \{ \delta(\wp) > u \text{ and } \wp \text{ is separating} \big\} \right) \\
&\quad\quad\quad\quad\quad\quad=\, \mu_\mathfrak{C}^{\text{bub}}\big( \{ \max \bpi(\partial^{\text{o}} \wp) > u \text{ and } \wp \text{ is separating} \} \big) \\
&\quad\quad\quad\quad\quad\quad=\, \mu_\mathbb{D}^{\text{bub}} \big( \{ \partial^{\text{o}} \wp \cap (e^{-u} \mathbb{D}) \neq \emptyset \text{ and } \wp \text{ surrounds }0 \}\big) \\
&\quad\quad\quad\quad\quad\quad\leq\, \mu_\mathbb{D}^{\text{bub}} \big(  \partial^{\text{o}} \wp \cap (e^{-u} \mathbb{D}) \neq \emptyset  \big).
\end{align*}
Now, by \cite{Lawler_2003}, we know that the outer boundary of a Brownian bubble measure in the disc is (up to scaling) identical to the $\mathrm{SLE}_{8/3}$ bubble measure. Combining this fact with Lemma~\ref{lem:SLE_bubble_hit}, the exponent simplifies to $1-\kappa/8 = 2/3$ and yields
$$
\mu_\mathbb{D}^{\text{bub}} \big(  \partial^{\text{o}} \wp \cap (e^{-u} \mathbb{D}) \neq \emptyset  \big) \,\asymp\,  \nu_\mathbb{D}^{\mathrm{SLE}_{8/3}} \big(\gamma \cap (e^{-u} \mathbb{D}) \neq \emptyset \big) \,\asymp\, (e^{-u})^{2/3}.
$$
Since $\mu_\mathfrak{C}^{\text{bub}, \tau}$ is the bubble measure restricted to the loops of duration at least $\tau$, and is a finite measure, the previous equation yields an upper bound for the desired mass.
\end{proof}

The previous lemma is key to conclude. Indeed, we know from Lemma~\ref{lem:ppp_naive_projection} that under $\mu_\mathfrak{C}^{\text{bub}}$, the projection length $\ell(\wp)$ has heavy tails (in fact, the tail is of order $\ell^{-2} \dd \ell$). In addition, the previous lemma shows that the difference between the core length $w(\wp)$ and the total projection length $\ell(\wp)$ is exactly $\delta(\wp)$ and has exponential tails. These light tails are too weak to change the heavy-tailed behaviour $\ell$ and change the universality class of the covering. Before proving that $\Lm$ satisfies Hypothesis (H), we need to prove that it is exponentially unlikely that a large loop is non-separating on the cylinder.

\begin{lemma} \label{lemma:l_contractible}
There exists a positive constant $C$ such that, for all $v$ sufficiently large,
$$
\mu_\mathfrak{C}^{\mathrm{bub}, \tau}\left( \big \{ \ell(\wp) > v \text{ and } \wp \text{ is non-separating} \big\} \right) \,\leq\, Ce^{- \frac{2}{3} v}.
$$
\end{lemma}

\begin{proof}
Let $f : re^{i\theta} \to (|\log r|, e^{i\theta})$ be the conformal map from $\mathbb{D} \setminus \{0\}$ to the cylinder. We have
\begin{align*}
&\mu_\mathfrak{C}^{\mathrm{bub}, \tau}\left( \big \{ \ell(\wp) > v \text{ and } \wp \text{ is non-separating} \big\} \right) \\
&\quad\quad\quad\quad\quad\quad=\, \mu_\mathbb{D}^{\mathrm{bub}, \tau} \left( \big \{ \wp \cap (e^{-v} \mathbb{D}) \neq \emptyset \text{ and } \wp \text{ does not surround } 0 \big\} \right).
\end{align*}
Now, notice that a loop $\wp$ does not surround $0$ if and only if its outer boundary $\partial^{\text{o}} \wp$ does not surround $0$. In that case, if in addition, the loop $\wp$ intersects the disc $e^{-v} \mathbb{D}$, then Jordan's theorem ensures that the outer boundary also has to intersect the disc $e^{-v} \mathbb{D}$ (see Figure~\ref{fig:outer_wp_hit}).

\begin{figure}[h]
    \centering
    \includegraphics[width=9cm]{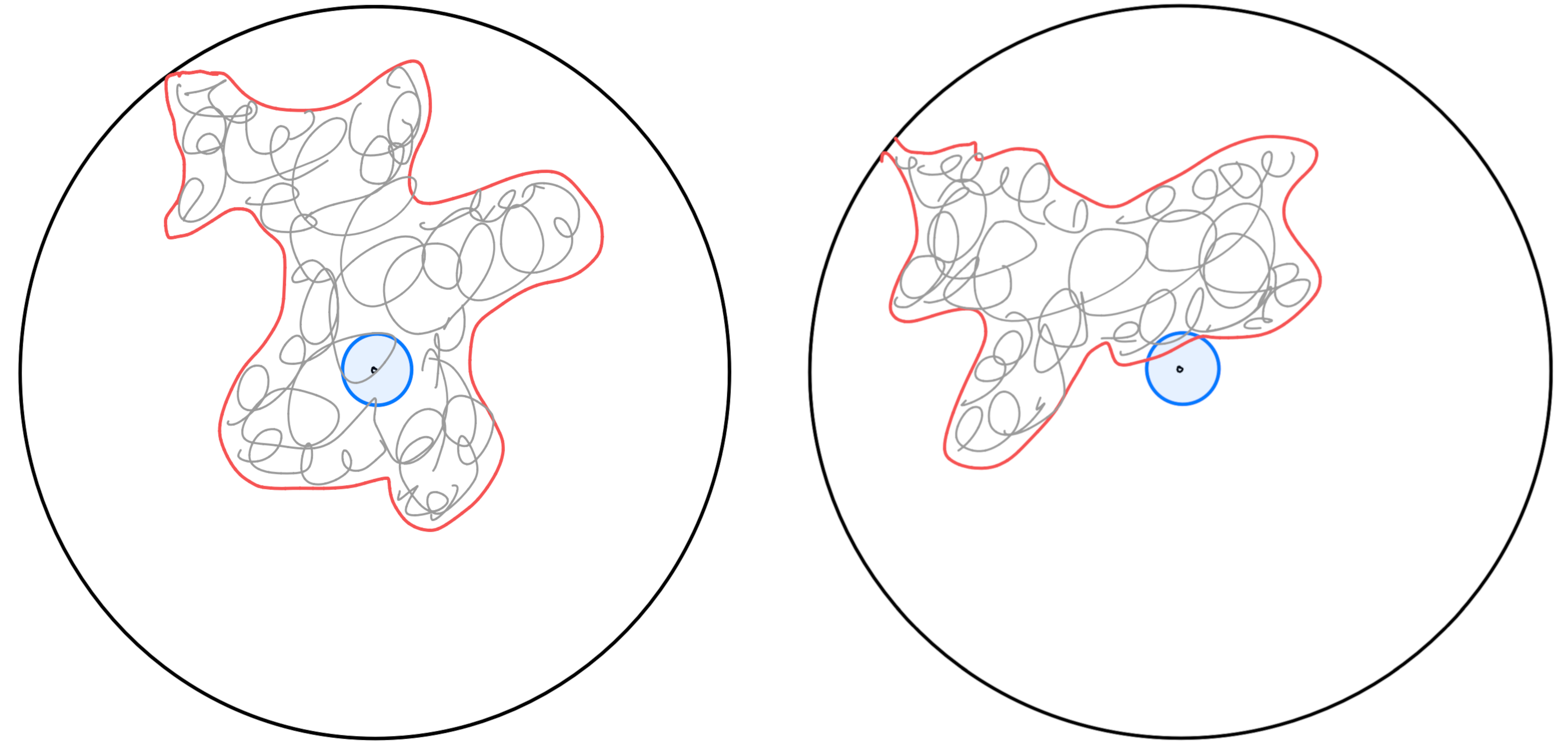}
    
    \caption{For a loop $\wp$ (in red and gray) in $\mathbb{D}$ that does not surround $0$ but hits $e^{-v} \mathbb{D}$ (in blue), its outer-boundary $\partial^{\text{o}} \wp$ (in red) has to hit $e^{-v} \mathbb{D}$. Right: the loop hits the disc but surrounds $0$, left: the loop hits the disc and does not surround $0$.}
    \label{fig:outer_wp_hit}
\end{figure}

\noindent Now, by \cite{Lawler_2003}, the outer boundary of a Brownian bubble measure in the disc is (up to scaling) identical to the $\mathrm{SLE}_{8/3}$ bubble measure. Combining this fact with Lemma~\ref{lem:SLE_bubble_hit}, the exponent simplifies to $1-\kappa/8 = 2/3$ and yields
\begin{align*}
\mu_\mathfrak{C}^{\mathrm{bub}, \tau}\left( \big \{ \ell(\wp) > v \text{ and } \wp \text{ is non-separating} \big\} \right) \,\leq\, c\nu_\mathbb{D}^{\mathrm{SLE}_{8/3}} \left( \gamma \cap (e^{-v} \mathbb{D}) \neq \emptyset \right) \,\asymp\, (e^{-v})^{2/3}.
\end{align*}
The proof is complete.
\end{proof}

\begin{lemma} \label{lem:lm_hypothesis_h} The random covering $\Lm$ satisfies Hypothesis (H).
\end{lemma}

\begin{proof}
Denoting $(a \otimes w) : \wp \to (a(\wp), w(\wp))$, the random set $\Lm$ is almost surely equal to $(a\otimes w) (\LC \cap \{\wp \text{ is separating}\})$. By the same arguments as in the proof of Lemma~\ref{lem:ppp_naive_projection}, we get that $\Lm$ is a Poisson point process of intensity measure
$$
\theta \,\dd a \otimes \mu_\mathfrak{C}^{\text{bub}, \tau}\big( w(\wp) \in \dd w \text{ and } \wp \text{ is separating}\big).
$$
For the rest of this proof, we will denote $\nu(\dd w) = \mu_\mathfrak{C}^{\text{bub}, \tau}( w(\wp) \in \dd w \text{ and } \wp \text{ is separating})$. To conclude the proof, it only remains to prove that the core length measure $\nu$ satisfies Hypothesis (H). We start with the first condition. Since, for $\mu_\mathfrak{C}^{\text{bub}, \tau}$-almost every bubble $\wp$, we have $w(\wp) \leq \ell(\wp)$ (by definition), we get 
$$
\nu((0,\infty))  \,\leq\, \mu_\mathfrak{C}^{\text{bub}, \tau}(\ell(\wp) > 0) \,=\, \mu((0,\infty)) \,<\, \infty.
$$
The last term is finite by the proof of Lemma~\ref{lem:ppp_naive_projection}.
We now prove the second condition. Recall that $w = \ell-\delta$ where $\ell$ is the total projection length of $\wp$ and $\delta$ is the boundary span defined in \eqref{eq:def_delta}. We proceed by upper and lower bounding $\nu((u, \infty))$. First, since we have seen that $w \leq \ell$, for all $u > 0$
\begin{equation}
\nu((u, \infty)) \,\leq\,\mu_\mathfrak{C}^{\text{bub}, \tau} \big(\{ \ell > u \text{ and } \wp \text{ is separating} \} \big) \,\leq\, \mu_\mathfrak{C}^{\text{bub}, \tau} \big( \ell > u \big). \label{eq:lower_mu_ell}
\end{equation}
For the sake of readability, in this proof, we simply denote $\ell, \delta, w$ for $\ell(\wp), \delta(\wp), w(\wp)$. We now turn our attention to the lower bound. Let $M > 0$ be some deterministic constant, note that if $\ell> u + M$ and $\delta \leq M$, then $\ell-\delta > u$. Hence, by subadditivity of the measure $\mu_\mathfrak{C}^{\text{bub}, \tau}$, we get
\begin{align}
\nu((u,\infty)) &= \mu_\mathfrak{C}^{\text{bub}, \tau} \big( \{\ell-\delta > u \text{ and } \wp \text{ is separating} \}\big)  \notag \\
&\geq\, \mu_\mathfrak{C}^{\text{bub}, \tau} \big( \{ \ell> u+M \text{ and } \delta \leq M \text{ and } \wp \text{ is separating} \} \big) \notag  \\
&\geq\, \mu_\mathfrak{C}^{\text{bub}, \tau} \big( \{ \ell> u+M  \text{ and } \wp \text{ is separating} \} \big)  \label{eq:upper_mu_M} \\
&\quad\quad\quad\quad\quad-  \mu_\mathfrak{C}^{\text{bub}, \tau} \big( \{ \delta > M  \text{ and } \wp \text{ is separating} \} \big). \notag 
\end{align}
In the last line, we have used the fact that for any measure $\mu$ we have $\mu(A\cap B \cap C) \geq \mu(A \cap C) - \mu(B^c \cap C)$. Now, by Equation~\eqref{eq:pointwise_mu_tail} of the proof of Lemma~\ref{lem:ppp_naive_projection}, we have
\begin{equation} \label{eq:bound_mu_ell}
\mu_\mathfrak{C}^{\text{bub}, \tau} \big( \ell > u \big) \,=\,  \frac{1}{u} + O(e^{-c u}).
\end{equation}
Consequently, for all $M \geq 0$,
\begin{align*}
&\mu_\mathfrak{C}^{\text{bub}, \tau} \big( \{ \ell> u+M  \text{ and } \wp \text{ is separating} \} \big) \\
&\quad\quad\quad\quad\quad\quad =\, \mu_\mathfrak{C}^{\text{bub}, \tau} \big( \ell > u+M \big) - \mu_\mathfrak{C}^{\text{bub}, \tau} \big( \ell > u+M \text{ and } \wp \text{ is non-separating} \big) \\
&\quad\quad\quad\quad\quad\quad =\, \frac{1}{u+M} + O(e^{-c(u+M)}) - O(e^{-2u/3}) \\
&\quad\quad\quad\quad\quad\quad =\, \frac{1}{u} - \frac{M}{u(u+M)} + O(e^{-c'u}). 
\end{align*}
Here $c' = c \wedge (2/3)$. In the third line, we have used the result from Lemma~\ref{lemma:l_contractible}. On the other hand, by Lemma~\ref{lemma:delta_exp}, we have, for all $M \geq 0$,
$$\mu_\mathfrak{C}^{\text{bub}, \tau} \big( \{ \delta > M  \text{ and } \wp \text{ is separating} \} \big) \,\leq\, Ce^{-\frac{2}{3}M}.$$
We now set $M = M(u) = 3\log u$. Injecting that specific choice of $M(u)$ and the previous bounds into Equation~\eqref{eq:upper_mu_M} yields
\begin{align*}
\nu((u,\infty)) \,&\geq\, \frac{1}{u} - \frac{3\log u}{u(u+3\log u)} + O(e^{-c' u}) -  Ce^{-\frac{2}{3} \cdot 3 \log u} \\
&\geq \frac{1}{u} - 3 \frac{\log u}{u^2} - \frac{C}{u^2} + O(e^{-c' u})
\end{align*}
By the bounds of Equation~\eqref{eq:upper_mu_M} and Equation~\eqref{eq:lower_mu_ell}, we get
$$
\frac{1}{u} - O\Big(\frac{\log u}{u^2} \Big) \,\leq\, \nu((u,\infty)) \,\leq\, \frac{1}{u} + O(e^{-c' u}),
$$
which of course implies that $\nu((u,\infty)) = u^{-1} + O(u^{-2} \log u)$ as $u \to \infty$. Finally, we get
$$
\int_1^x \nu((u,\infty)) \dd u \,=\, \int_1^x \Big( \frac{1}{u} + O\Big( \frac{\log u}{u^2} \Big) \Big) \dd u \,=\, \log x + c'' + o(1),
$$
as $x \to \infty$. This is the desired result.
\end{proof}

\noindent \textbf{Proof of Theorem~\ref{th:phase_transition_loopsoup}.} By Lemma~\ref{lem:ppp_naive_projection} and Lemma~\ref{lem:lm_hypothesis_h}, both the random covering of projections $\Lp$ and the random covering of core intervals $\Lm$ satisfy Hypothesis (H). By Proposition~\ref{prop:random_coverings}, if $\theta \in (0,1]$, almost surely, both the 1D random coverings have no unbounded connected component. On the other hand, if $\theta > 1$, almost surely, both $\Lp$ and $\Lm$ have a unique unbounded connected component. The sandwiching property of Lemma~\ref{lemma:coupling_coverings} readily implies that for $\theta \leq 1$, almost surely, there is no unbounded cluster in the loop soup $\LC$, whereas if $\theta > 1$, almost surely, there exists (at least one) unbounded cluster. The sandwiching Lemma~\ref{lemma:coupling_coverings} also implies that for all $0<y<x$,
\begin{align*}
\mathbb{P} \Big([y,x] \text{ is covered by } \Lm\Big)
 \,\leq\, \mathbb{P} \Big( \{y\} \times \mathbb{S}^1  \stackrel{\LC}{\longleftrightarrow} \{x\} \times \mathbb{S}^1 \Big) \,\leq\, \mathbb{P}\Big([y,x] \text{ is covered by } \Lp\Big)
\end{align*}
Further applying Proposition~\ref{prop:random_coverings} to both the random coverings $\Lm$ and $\Lp$, we get, as $x \rightarrow \infty$ and $y > 0$ is fixed,
\begin{align*}
&\mathbb{P} \Big( \{y\} \times \mathbb{S}^1  \stackrel{\LC}{\longleftrightarrow} \{x\} \times \mathbb{S}^1 \Big) \,\asymp\, x^{-1 + \theta} \quad\quad\, \text{if } \theta \in (0,1) \\
&\mathbb{P} \Big( \{y\} \times \mathbb{S}^1  \stackrel{\LC}{\longleftrightarrow} \{x\} \times \mathbb{S}^1 \Big) \,\asymp\, (\log x)^{-1} \quad \text{if } \theta = 1
\end{align*}
It remains to prove that if $\theta > 1$, the unbounded connected component is unique. For this, notice that the loops $\wp\in \LC$ that are separating form an unbounded cluster (because their core intervals $\Lm$ form an unbounded cluster on $\mathbb{R}_+$). Since each loop is bounded, we can find a sequence $(\wp_n)_{n \geq 1}$ of separating loops such that $a(\wp_n) \to \infty$, which prevents any other cluster from being unbounded.
$\hfill \square$

\section{Sharp crossing exponent in the Brownian loop soup} \label{sec:crossing_exponent}

In this section, we prove Theorem~\ref{th:crossing_exponent}. We consider the non-truncated Brownian loop soup $\mathcal{L}_\mathbb{D}^\theta$ of intensity $\theta \in (0,1/2]$ in the unit disc. Under the conformal map $re^{i\theta} \mapsto (|\log r|, e^{i\theta})$, setting $y = -\log \rho$ and $x = -\log r$, the cylinder slice $(y,x) \times \mathbb{S}^1$ corresponds to the annulus $\rho \mathbb{D} \setminus r \mathbb{D}$ in the unit disc, and the crossing event translates to the cylinder $\mathfrak{C} = (0, \infty) \times \mathbb{S}^1$ as
$$
\mathbb{P} \Big( \rho \partial \mathbb{D}  \stackrel{\mathcal{L}_\mathbb{D}^{\theta}}{\longleftrightarrow} r \partial \mathbb{D} \Big) \,=\, \mathbb{P} \Big( \{y\} \times \mathbb{S}^1  \stackrel{\mathcal{L}_\mathfrak{C}^{\theta}}{\longleftrightarrow} \{x\} \times \mathbb{S}^1 \Big).
$$
First, the lower bound is an immediate consequence of our previous results. Indeed, since the truncated loop soup satisfies $\mathcal{L}_\mathfrak{C}^{\theta, \tau} \subset \mathcal{L}_\mathfrak{C}^\theta$ for any fixed threshold $\tau > 0$, and crossing events are increasing, Theorem~\ref{th:phase_transition_loopsoup} directly implies
\begin{equation} \label{eq:crossing_lower_bound}
    \mathbb{P} \Big( \{y\} \times \mathbb{S}^1  \stackrel{\mathcal{L}_\mathfrak{C}^{\theta}}{\longleftrightarrow} \{x\} \times \mathbb{S}^1 \Big) \,\geq\, \mathbb{P} \Big( \{y\} \times \mathbb{S}^1  \stackrel{\mathcal{L}_\mathfrak{C}^{\theta, \tau}}{\longleftrightarrow} \{x\} \times \mathbb{S}^1 \Big) \,\asymp\, x^{-1+\theta}.
\end{equation}
Switching back to the unit disc, we get the  $|\log r|^{-1+\theta}$ lower bound. The entire goal of this section is therefore to establish the upper bound. To this end, we completely abandon the duration cutoff $\tau$ used in Section~\ref{sec:phase_transition} in favor of a \emph{diameter cutoff}: For any loop $\wp$ on $\mathfrak{C}$, we define its diameter $\mathrm{diam}(\wp)$ as the diameter of the smallest disc containing $\wp$ (see Equation~\eqref{eq:def_metric_cyl} for a definition of the metric on the cylinder). For a fixed small parameter $\varepsilon \in (0, 1)$, we partition the soup into two independent Poisson point processes:
$$
\LCl \,:=\, \mathcal{L}_\mathfrak{C}^\theta \cap \{\wp : \mathrm{diam}(\wp) \le \varepsilon\} \quad \text{and} \quad \LCg \,:=\, \mathcal{L}_\mathfrak{C}^\theta \cap \{\wp : \mathrm{diam}(\wp) > \varepsilon\},
$$
which we refer to as the \emph{small loops} and the \emph{large loops}, respectively. A diameter cutoff is much better here because it provides strict spatial localisation, which is essential to ensure independence between disjoint regions of the cylinder (Because of the Poissonian nature of the loop soup, the sets of small loops that intersect domains at distance at least $2\varepsilon$ from one another are independent). 

Note that we cannot use the same strategy as in Section~\ref{sec:phase_transition} and project the small loops $\LCl$ onto $(0, \infty)$. Indeed, because of the infinite density of small loops, the projection length measure of $\LCl$ satisfies $\mu^{\leq \varepsilon}(\dd \ell) \asymp  \mathbf{n}( p_{\mathbb{S}^1}(\zeta) \mathbf{1}_{\{ \sup_{t \geq 0} X_t \in \dd \ell\}}) \asymp \ell^{-3} \dd \ell$ as $\ell \to 0$. Therefore, it does not satisfy Hypothesis~(H). In fact, the small loops cover the entire half-line:
$$
\bpi(\LCl) \,=\, (0, \infty) \quad \text{almost surely.}
$$
Before proving Theorem~\ref{th:crossing_exponent}, we first need to check that this new diameter cutoff does not break our previous results. Indeed, all the estimates in Section~\ref{sec:phase_transition} were established using a duration cutoff $T(\wp) \ge \tau$.  

\subsection{Projection intervals with a diameter cutoff.} \label{sec:diameter_projection}

In this section, we introduce the direct analogue of the covering $\Lp$ from Section~\ref{sec:phase_transition}, using a diameter cutoff instead of a duration cutoff:
$$
\overline{\Lambda}{}^{\theta, > \varepsilon} \,:=\, \big\{ (a(\wp), \ell(\wp)) \big\}_{\wp \in \LCg}.
$$
Our goal is to show that this change of cutoff does not affect the behavior of the covering. Specifically, we prove that $\overline{\Lambda}{}^{\theta, > \varepsilon}$ still satisfies Hypothesis~(H), so that all the 1D covering and renewal results from Section~\ref{sec:random_coverings} remain valid.

\begin{lemma} \label{lemma:diameter_hypothesis_H}
    The random covering $\overline{\Lambda}{}^{\theta, > \varepsilon}$ satisfies Hypothesis (H).
\end{lemma}

\begin{proof}
    Since $(a\otimes \ell)$ is a measurable functional, the same arguments as in the proof of Lemma~\ref{lem:ppp_naive_projection} yield that the projection covering of the large loops $\overline{\Lambda}{}^{\theta, > \varepsilon}$ is a Poisson point process of intensity $\theta \dd a \otimes \mu^{> \varepsilon}(\dd \ell)$ where $\mu^{> \varepsilon}(\dd \ell) = \mu_\mathfrak{C}^{\mathrm{bub}}(\ell(\wp) \in \dd \ell \text{ and } \mathrm{diam}(\wp) > \varepsilon)$ is the Brownian bubble's length measure on the cylinder restricted to the loops of diameter $\mathrm{diam}(\wp) > \varepsilon$. First, by \cite{lawler2003brownianloopsoup}, the Brownian bubble measure is locally finite, hence $\mu^{> \varepsilon}((0, \infty)) < \infty$, which implies that the length measure of $\overline{\Lambda}{}^{\theta, > \varepsilon}$ is finite and concludes the proof of the first condition. Now, note that if a loop $\wp$ satisfies $\ell(\wp) > \varepsilon$, then, $\mathrm{diam}(\wp) > \varepsilon$. Hence for all $v \geq \varepsilon$,
    $$
    \mu^{> \varepsilon}((v, \infty)) \,=\, \mu_\mathfrak{C}^{\text{bub}}(\ell(\wp) > v).
    $$
    Introducing an arbitrary time threshold $\tau = 1$ (in order to apply our previous results), we get
    \begin{align*}
    \mu^{> \varepsilon}((v,\infty)) \,&=\, \mu_\mathfrak{C}^{\text{bub}, \tau}(\ell(\wp) > v) + \mu_\mathfrak{C}^{\text{bub}}(\ell(\wp) > v \text{ and } T(\wp) \leq \tau) \\        
    &=\, v^{-1} + O(e^{-c v}) + O(e^{-c' v^2}).
    \end{align*}
    In the last line, we have used $\mu_\mathfrak{C}^{\text{bub}, \tau}(\ell(\wp) > v) = \frac{1}{v} + O(e^{-cv})$ (Proof of Lemma~\ref{lem:ppp_naive_projection}, Equation~\eqref{eq:pointwise_mu_tail}) together with the proof of Lemma~\ref{lemma:hits_and_tau} (Equation~\eqref{eq:mu_bub_bound}). Integrating this relation  implies that $\int_1^x \mu^{> \varepsilon}((v,\infty)) \dd v = \log x + c'' + o(1)$.
\end{proof}

\subsection{Filling large gaps with clusters of small loops.}

By definition, loops in $\LCl$ have diameter at most $\varepsilon$. Therefore, the boundary condition at $\partial \mathfrak{C} = \{0\} \times \mathbb{S}^1$ only affects the slice $(0, \varepsilon) \times \mathbb{S}^1$. For any $x \geq \varepsilon$, the soup $\LCl$ inside the shifted cylinder $\mathfrak{C}_x := (x,\infty) \times \mathbb{S}^1$ thus has the same law as on the infinite cylinder $\mathbb{R} \times \mathbb{S}^1$. In particular, it is invariant under the longitudinal translations that send $\mathfrak{C}_{x-\varepsilon}$ to a subset of $\mathfrak{C}$.

\begin{lemma} \label{lemma:translation_invariant}
Let $x \ge \varepsilon$, and let $E$ be an event measurable with respect to the loops contained in $\mathfrak{C}_{x-\varepsilon}$. Then, for all $a \in \mathbb{R}$ such that $x + a \ge \varepsilon$, we have $\mathbb{P}(\tau_a E) = \mathbb{P}(E)$.
\end{lemma}

\begin{proof}
On the infinite cylinder $\mathbb{R} \times \mathbb{S}^1$, the Brownian small-loop soup $\mathcal{L}_{\mathbb{R}\times \mathbb{S}^1}^{\theta, \leq \varepsilon}$ is invariant under the longitudinal translations $(\tau_a)_{a \in \mathbb{R}}$. Since $x \geq \varepsilon$, any loop $\wp$ with diameter at most $\varepsilon$ contained in $\mathfrak{C}_{x-\varepsilon}$ is also contained in $\mathfrak{C}$. Hence, the small loops of $\LCl$ contained in $\mathfrak{C}_{x-\varepsilon}$ are exactly those of the loop soup on the full infinite cylinder $\mathbb{R} \times \mathbb{S}^1$ in that same region. Since $x + a \ge \varepsilon$, the same holds for the shifted region $\mathfrak{C}_{x+a-\varepsilon}$, and the result follows immediately from the translation invariance on $\mathbb{R} \times \mathbb{S}^1$.
\end{proof}

In Theorem~\ref{th:phase_transition_loopsoup}, we have shown that the crossing probability for the large loops of $\LCg$ decays polynomially as $x^{-1+\theta}$. In contrast, the following lemma shows that crossings made by clusters of small loops decay exponentially fast. Most importantly, by the translation invariance of Lemma~\ref{lemma:translation_invariant}, this bound is uniform along the cylinder and only depends on the crossing distance.

\begin{lemma} \label{lemma:exp_crossing} Let $\theta \in (0,\frac{1}{2}]$, there exist positive constants $C,c > 0$, such that, for any $0 < u < v$,
    $$\mathbb{P} \Big( \{u\} \times \mathbb{S}^1  \stackrel{\LCl}{\longleftrightarrow} \{v\} \times \mathbb{S}^1 \Big) \,\leq\, Ce^{-c(v-u)}.
$$
\end{lemma}

\begin{proof}
Let $0 < u < v$ and set $K = \lfloor \frac{v-u}{4\varepsilon} \rfloor$. For each $0 \le k < K$, we define the slice $\mathfrak{C}_{k, \varepsilon} := (u + 4\varepsilon k + \varepsilon, \, u + 4\varepsilon k + 3\varepsilon) \times \mathbb{S}^1$, and consider the event
\[
E_k := \big\{ \text{no cluster of } \LCl \text{ crosses } \mathfrak{C}_{k, \varepsilon} \big\}.
\]
Because every loop in $\LCl$ has diameter at most $\varepsilon$, for all $k \geq 0$, the event $E_k$ is measurable with respect to the loops contained in $(u + 4\varepsilon k, \, u + 4\varepsilon(k+1)) \times \mathbb{S}^1$.  Since these slices are all disjoint, the spatial independence property of Poisson point processes ensures that the collections of loops contained in each slice are mutually independent. Hence, the events $(E_k)_{0 \le k < K}$ are mutually independent. Moreover, since each $E_k$ is measurable with respect to the loops contained in $\mathfrak{C}_{u+4k\varepsilon}$, Lemma~\ref{lemma:translation_invariant} with $a = -u-4k\varepsilon$ yields
\begin{equation} \label{eq:def_p_gamma}
p \,:=\, \mathbb{P}(E_k) \,=\, \mathbb{P}(\tau_{a} E_k) \,=\, \mathbb{P} \big( \text{no cluster of } \LCl \text{ crosses } (\varepsilon, 3\varepsilon) \times \mathbb{S}^1 \big),
\end{equation}
which implies that $p := \mathbb{P}(E_k)$ depends neither on $k$ nor on $u$. For a crossing to occur from $\{u\} \times \mathbb{S}^1$ to $\{v\} \times \mathbb{S}^1$, each event $E_k$ must fail (because each subslice must of course also be crossed), which yields
\[
\mathbb{P} \Big( \{u\} \times \mathbb{S}^1 \stackrel{\mathcal{L}_\mathfrak{C}^{\theta,\le \varepsilon}}{\longleftrightarrow} \{v\} \times \mathbb{S}^1 \Big) \,\leq\, \mathbb{P}\left( E_0^c \cap \cdots \cap E_{K-1}^c \right) \,=\, (1 - p)^K \,\leq\, C e^{-c(v-u)},
\]
with $c = -(4\varepsilon)^{-1}\log(1 - p) > 0$. It remains to check that $p > 0$. The event defining $p$ in \eqref{eq:def_p_gamma} is decreasing (adding more loops can only make a crossing easier). Since $\LCl \subset \mathcal{L}_\mathfrak{C}^\theta$, by Equation~\eqref{eq:def_p_gamma} we get
\begin{align*}
p \,&\geq\, \mathbb{P}\big(\text{no cluster of } \mathcal{L}_\mathfrak{C}^\theta \text{ crosses } (\varepsilon, 3\varepsilon) \times \mathbb{S}^1\big)  \\
&=\, \mathbb{P}\big(\text{no cluster of } \mathcal{L}_\mathbb{D}^\theta \text{ crosses } e^{-\varepsilon} \mathbb{D} \setminus e^{-3\varepsilon} \overline{\mathbb{D}} \big).
\end{align*}
In the last line, we have applied the conformal map $(x, e^{i\varphi}) \mapsto e^{-x+i\varphi}$: the cylinder slice is mapped onto the annulus $e^{-\varepsilon}\mathbb{D} \setminus e^{-3\varepsilon}\mathbb{D}$. Now, denote $\mathcal{L}_{e^{-\varepsilon} \mathbb{D}}$ the loops of $\mathcal{L}_\mathbb{D}^\theta$ contained in $e^{-\varepsilon} \mathbb{D}$ and $\mathcal{L}_{\mathbb{D} \setminus e^{-\varepsilon} \mathbb{D}}$ the remaining loops. Denote $\mathcal{C}$ the outermost cluster of $\mathcal{L}_{e^{-\varepsilon} \mathbb{D}}$ surrounding the origin. By \cite[Lemma 2.5]{Jego_2025}, there is a positive probability $p_1$ that some loop $\wp \in \mathcal{L}_{e^{-\varepsilon} \mathbb{D}}$ surrounds $e^{-3\varepsilon} \mathbb{D}$, on this event, the outer boundary of $\mathcal{C}$ has to be contained in the annulus $e^{-\varepsilon} \mathbb{D} \setminus e^{-3\varepsilon} \mathbb{D}$. One can note that if no loop of $\mathcal{L}_{\mathbb{D} \setminus e^{-\varepsilon} \mathbb{D}}$ hits $\mathcal{C}$ (which happens with positive probability $p_2$, since $\mathcal{L}_{\mathbb{D} \setminus e^{-\varepsilon} \mathbb{D}}$ is Poissonian, locally finite in $e^{-\varepsilon} \mathbb{D}$ and $\mathcal{C}$ is almost surely at positive distance from $e^{-\varepsilon}\partial \mathbb{D}$), then no cluster of $\mathcal{L}_\mathbb{D}^\theta$ crosses $e^{-\varepsilon} \mathbb{D} \setminus e^{-3\varepsilon} \mathbb{D}$. Hence, since $\mathcal{L}_{e^{-\varepsilon} \mathbb{D}}$ and $\mathcal{L}_{\mathbb{D} \setminus e^{-\varepsilon} \mathbb{D}}$ are independent, we get $p \geq p_1p_2 > 0$.
\end{proof}

Fix $0 < y < x$. Let $\overline{\Lambda}{}^{\theta,> \varepsilon}$ denote the set of projection intervals of the large loops as in Section~\ref{sec:diameter_projection}, and let $\mathcal{R}^{\theta,> \varepsilon} = (0,\infty) \setminus \bigcup_{I \in \overline{\Lambda}{}^{\theta,> \varepsilon}} I$ be the uncovered projection set. Because the diameter threshold $\varepsilon > 0$ is strictly positive, the number of large loops in $\mathcal{L}_\mathfrak{C}^{\theta,> \varepsilon}$ intersecting the cylinder slice $(y,x) \times \mathbb{S}^1$ is almost surely finite. Consequently, the uncovered set restricted to this slice, $\mathcal{R}^{\theta,> \varepsilon} \cap (y,x)$, consists of a finite union of disjoint intervals, which we refer to as the \emph{projection gaps} of the large loops and denote by $\mathcal{G}(y,x)$. Finally, we denote by $N(y,x) = \# \mathcal{G}(y,x)$ the number of projection gaps.
\medskip

\begin{figure}[h]
    \centering
    \includegraphics[width=13.5cm]{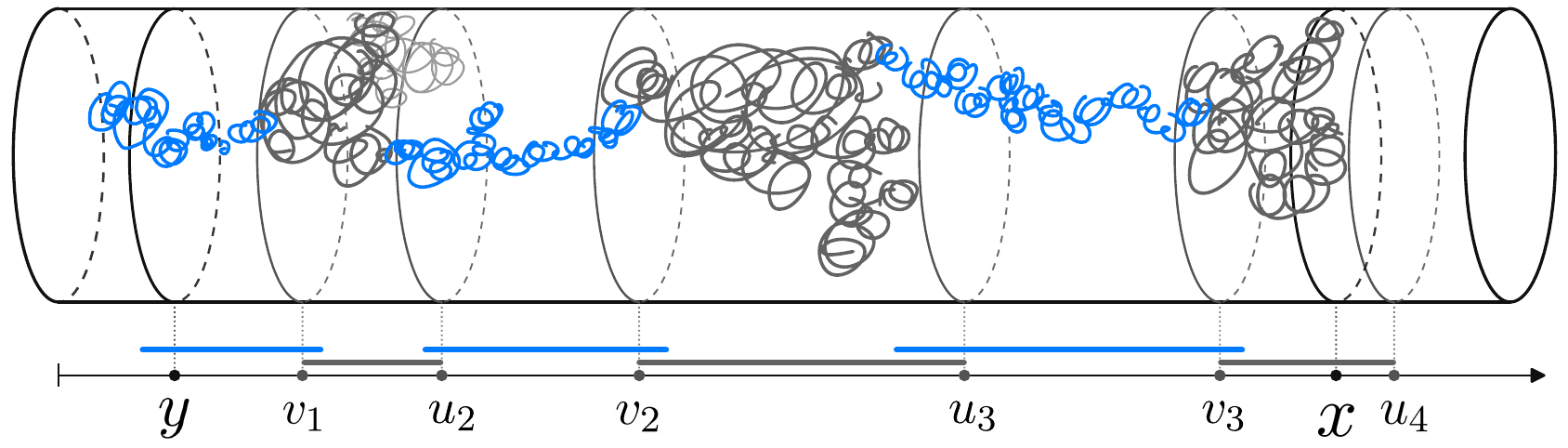}

    \caption{A cluster $\mathcal{C}$ in $\mathcal{L}_\mathfrak{C}^\theta$ crossing from $\{y\} \times \mathbb{S}^1$ to $\{x\} \times \mathbb{S}^1$, composed of large loops (dark gray) and small loops (blue). Below the cylinder, the 1D projection intervals of large loops $\overline{\Lambda}{}^{\theta, > \varepsilon} = \bigcup_{i} (v_i,u_{i+1})$ leave uncovered gaps $\mathcal{R}^{\theta, > \varepsilon}$. For $\mathcal{C}$ to make the crossing from $y$ to $x$, each gap left by the large loops must be crossed by some cluster consisting only of small loops. Here $N(y,x) = 3$ and $\mathcal{G}(y,x) = \{(y,v_1), (u_2,v_2), (u_3,v_3)\}$.}
    \label{fig:multiple_crossings}
\end{figure}

Suppose now that there is a crossing from $\{y\} \times \mathbb{S}^1$ to $\{x\} \times \mathbb{S}^1$ in $\mathcal{L}_\mathfrak{C}^\theta$. This means that there exists a cluster $\mathcal{C}$ of the full non-truncated loop soup $\mathcal{L}_\mathfrak{C}^\theta$ that intersects both $\{y\} \times \mathbb{S}^1$ and $\{x\} \times \mathbb{S}^1$, which implies that its longitudinal projection $\bpi(\mathcal{C})$ contains the entire interval $(y,x)$. For any projection gap $(u,v) \in \mathcal{G}(y,x)$, the definition of $\mathcal{R}^{\theta,> \varepsilon}$ ensures that no loop of diameter at least $\varepsilon$ has a projection intersecting $(u,v)$. Hence, no large loop can intersect the cylinder slice $(u,v) \times \mathbb{S}^1$. As a consequence, the portion of $\mathcal{C}$ crossing this slice is composed exclusively of small loops from $\LCl$. Since $\bpi(\mathcal{C})$ covers $(u,v)$, there must exist a connected component consisting only of small loops that crosses from $\{u\} \times \mathbb{S}^1$ to $\{v\} \times \mathbb{S}^1$. Because this condition must be satisfied for all the gaps in $\mathcal{G}(y,x)$, we obtain the following necessary condition:
\begin{equation}
\label{eq:crossing_inclusion}
\Big\{ \{y\} \times \mathbb{S}^1 \stackrel{\mathcal{L}_\mathfrak{C}^\theta}{\longleftrightarrow} \{x\} \times \mathbb{S}^1 \Big\} \,\subseteq\, \Big\{ \forall (u,v) \in \mathcal{G}(y,x) : \, \{u\} \times \mathbb{S}^1 \stackrel{\LCl}{\longleftrightarrow} \{v\} \times \mathbb{S}^1 \Big\}.
\end{equation}
Note that this condition covers all cases: if $(y,x)$ is completely covered by $\overline{\Lambda}{}^{\theta,> \varepsilon}$, the set $\mathcal{G}(y,x)$ is empty and the condition holds trivially. On the other hand, if no large loops appear in $(y,x) \times \mathbb{S}^1$, then $\mathcal{G}(y,x) = \{(y,x)\}$, the small loops must cross the entire slice $(y,x) \times \mathbb{S}^1$ by themselves. A schematic representation of \eqref{eq:crossing_inclusion} is available in Figure~\ref{fig:multiple_crossings}.
\medskip

Directly bounding the right-hand event in~\eqref{eq:crossing_inclusion} is difficult. Indeed, the gaps in $\mathcal{G}(y,x)$ can be arbitrarily small and arbitrarily close to each other. If a gap is too small, the bound $C e^{-c(v-u)}$ from Lemma~\ref{lemma:exp_crossing} is larger than $1$, so it does not provide any penalty. In addition, gaps that are too close create unwanted correlations between crossings. However, crossing all gaps trivially implies crossing any sub-collection of them. We can therefore focus only on gaps that are both large and well-separated. Hence, we fix some large parameter $K > 0$, and we select a largest subset $\mathcal{G}^*(y,x) \subset \mathcal{G}(y,x)$ such that every gap $(u,v)$ satisfies $v - u \geq K$, and any two gaps in $\mathcal{G}^*(y,x)$ are at a distance of at least $3\varepsilon$ from one another (to make this selection measurable and well defined, one can simply use a greedy algorithm from left to right).  Denoting by $N^*(y,x) := \#\mathcal{G}^*(y,x)$ the number of well-behaved gaps and restricting~\eqref{eq:crossing_inclusion} to $\mathcal{G}^*(y,x)$ yields the simpler inclusion

\begin{equation}
\label{eq:crossing_inclusion_good}
\Big\{ \{y\} \times \mathbb{S}^1 \stackrel{\mathcal{L}_\mathfrak{C}^\theta}{\longleftrightarrow} \{x\} \times \mathbb{S}^1 \Big\} \,\subseteq\, \Big\{ \forall (u,v) \in \mathcal{G}^*(y,x) : \, \{u\} \times \mathbb{S}^1 \stackrel{\LCl}{\longleftrightarrow} \{v\} \times \mathbb{S}^1 \Big\}.
\end{equation}
The advantage of~\eqref{eq:crossing_inclusion_good} is that every well-behaved gap from $\mathcal{G}^*(y,x)$ has a length bounded from below by $K$, ensuring an exponential penalty of each gap by Lemma~\ref{lemma:exp_crossing}. Furthermore, the distance $3\varepsilon$ between consecutive gaps ensures independence of crossings by the small loops. The following lemma will be key to the proof.

\begin{lemma} \label{lemma:filling_gaps_BK}
    For $K > 0$ sufficiently large, there exists $q(K) \in (0,1)$ such that for all $0 < u_1 < v_1 < u_2 < v_2 < \dots < u_n < v_n$, satisfying $v_i-u_i \geq K$ and $u_{i+1}-v_i \geq 3\varepsilon$ we have
$$
\mathbb{P} \Big( \forall i \in \llbracket 1,n \rrbracket: \, \{u_i\} \times \mathbb{S}^1 \stackrel{\mathcal{L}_\mathfrak{C}^{\theta,\le \varepsilon}}{\longleftrightarrow} \{v_i\} \times \mathbb{S}^1 \Big) \,\leq\, q(K)^n.
$$
\end{lemma}

\begin{proof} For $i = 1, \dots, n$, define the events $E_i = \{ \exists \mathcal{C} \text{ cluster of } \LCl \text{ crossing } (u_i,v_i) \times \mathbb{S}^1 \}$. Since each loop $\wp \in \LCl$ has a diameter at most $\varepsilon$, the events $E_i$ are measurable with respect to the loops contained in $(u_i-\varepsilon, v_i+\varepsilon) \times \mathbb{S}^1$. Moreover, since $u_{i+1}-v_i \geq 3\varepsilon$, these regions are disjoint. By spatial independence of Poisson point processes, the events $(E_i)_{1\leq i\leq n}$ are independent. Combining this with Lemma~\ref{lemma:exp_crossing} and the fact that $v_i-u_i \geq K$ for all $i$, we get
\begin{align*}
    \mathbb{P} \Big( \forall i \in \llbracket 1,n \rrbracket: \, \{u_i\} \times \mathbb{S}^1 \stackrel{\mathcal{L}_\mathfrak{C}^{\theta,\le \varepsilon}}{\longleftrightarrow} \{v_i\} \times \mathbb{S}^1 \Big) \,&=\, \prod_{i=1}^n \mathbb{P} \Big( \{u_i\} \times \mathbb{S}^1 \stackrel{\mathcal{L}_\mathfrak{C}^{\theta,\le \varepsilon}}{\longleftrightarrow} \{v_i\} \times \mathbb{S}^1 \Big) \\    
    &\leq\, \big(Ce^{-cK} \big)^n.
\end{align*}
Choosing $K$ sufficiently large, we get $Ce^{-cK} < 1$, which concludes the proof.
\end{proof}

We fix $K > 0$ (that we choose so that $K > 3\varepsilon$) for the remainder of the proof, and denote  $q := q(K) \in (0,1)$.  Partitioning the Brownian loop soup according to the loop diameter splits the loop soup into two independent Poisson point processes: $\mathcal{L}_\mathfrak{C}^\theta = \LCl \cup \LCg$. By construction, the collection of well-behaved gaps $\mathcal{G}^*(y,x)$ and its number of elements $N^*(y,x)$ depend only on the configuration of large loops $\LCg$, and are therefore measurable with respect to the set of large loops. By Equation~\eqref{eq:crossing_inclusion_good}, conditioning on the set of large loops and using the law of total expectation, we obtain
\begin{equation*} \label{eq:cond_exp_crossing}
\mathbb{P} \Big( \{y\} \times \mathbb{S}^1 \stackrel{\mathcal{L}_\mathfrak{C}^\theta}{\longleftrightarrow} \{x\} \times \mathbb{S}^1 \Big) 
\,\leq\, \mathbb{E} \left[ \mathbb{P} \Big( \forall (u,v) \in \mathcal{G}^*(y,x) : \, \{u\} \times \mathbb{S}^1 \stackrel{\LCl}{\longleftrightarrow} \{v\} \times \mathbb{S}^1 \,\Big|\, \LCg \Big) \right].
\end{equation*}
Conditionally on $\LCg$, the set of well-behaved gaps $\mathcal{G}^*(y,x)$ is deterministic. Moreover, by definition, each interval $(u_i,v_i) \in \mathcal{G}^*(y,x)$, satisfies $v_i - u_i \geq K$ and $u_{i+1} - v_i \geq 3\varepsilon$. Since the set of small loops $\LCl$ is independent of $\LCg$, Lemma~\ref{lemma:filling_gaps_BK} applies conditionally almost surely, yielding that the inner probability of filling all well-behaved gaps by clusters of small loops is almost surely bounded by $q^{N^*(y,x)}$. Plugging this bound back, we get,
\begin{align} \label{eq:_bound_crossing_Nstar}
\mathbb{P} \Big( \{y\} \times \mathbb{S}^1 \stackrel{\mathcal{L}_\mathfrak{C}^\theta}{\longleftrightarrow} \{x\} \times \mathbb{S}^1 \Big)  \,&\leq\, \mathbb{E} \big[ q^{N^*(y,x)} \big].
\end{align}
To prove Theorem~\ref{th:crossing_exponent}, it only remains to show that $\mathbb{E} \big[ q^{N^*(y,x)} \big]$ is of order $x^{-1+\theta}$. Note that we cannot apply Lemma~\ref{lemma:generating_function_Nst} here, because it only applies to the total number of gaps $N(y,x)$ and not to a strict subset of those. \newline

\noindent \textbf{Proof of Theorem~\ref{th:crossing_exponent}.} 
Recall that $N^*(y,x) = \# \mathcal{G}^*(y,x)$ denotes the number of well-behaved gaps in $(y,x)$ (of length at least $K$ and separated by at least $3\varepsilon$). To study these large gaps, we introduce the widened random covering of $(0,\infty)$ defined by 
$$
\Lambda^{(K)} \,:=\, \big\{(a_i, \ell_i + K) \big\}_{(a_i, \ell_i) \in \overline{\Lambda}{}^{\theta, > \varepsilon}},
$$
and we denote by $N_K(y,x)$ the number of gaps of $\Lambda^{(K)}$ in $(y,x)$. Extending each interval by $K$ to the right completely closes any gap of $\overline{\Lambda}{}^{\theta, > \varepsilon}$ of length smaller than $K$, while any gap $(u,v)$ of length at least $K$ survives as a smaller gap $(u+K, v)$. Therefore, projection intervals of large loops $\overline{\Lambda}{}^{\theta, > \varepsilon}$ contain at least $N_K(y,x)$ disjoint gaps of length at least $K$ between $y$ and $x$.

From these gaps, we can extract the well-behaved $\mathcal{G}^*(y,x)$ using a greedy algorithm. Order the gaps of $\overline{\Lambda}{}^{\theta, > \varepsilon}$ of length at least $K$ from left to right: $\{ I_1, I_2, \dots \}$ (with $|I_i| \geq K$). We select $I_1$, remove any other gap intersecting $[\sup I_1, \sup I_1 + 3\varepsilon]$, and repeat on the remaining gaps. Since each gap has length at least $K$ and $3\varepsilon < K$, the interval $[\sup I_1, \sup I_1 + 3\varepsilon]$ has a length of $3\varepsilon < K$ strictly smaller than $K$, so it can intersect at most one other gap. Thus, each selection eliminates at most one other gap from the list $\{I_j, j \geq 1\}$, which proves that almost surely,
$N^*(y,x) \,\geq\, \frac{1}{2} (N_K(y,x)-1)$ (the offset $N_K-1$ comes from the fact that the gap containing $\{y\}$ may have a length $< K$ in $\mathcal{G}(y,x)$ while being counted in $N_K(y,x)$). Since $q \in (0, 1)$, this implies
\begin{equation} \label{eq:q_Nstar_NK}
\mathbb{E} \big[ q^{N^*(y,x)} \big] \,\leq\, q^{-1/2}\mathbb{E} \big[ (\sqrt{q})^{N_K(y,x)} \big].    
\end{equation}
Furthermore, denoting the map $p : (a,\ell) \mapsto (a, \ell+K)$, we have $\Lambda^{(K)} = p(\overline{\Lambda}{}^{\theta, > \varepsilon})$. By the mapping theorem, $\Lambda^{(K)}$ is a Poisson point process. Recall from Lemma~\ref{lemma:diameter_hypothesis_H} that the intensity measure of $\overline{\Lambda}{}^{\theta, > \varepsilon}$ is given by $\theta \dd a \otimes \mu^{> \varepsilon}(\dd \ell)$, where $\mu^{> \varepsilon}$ was introduced as the bubble length measure restricted to loops of diameter at least $\varepsilon$. Therefore, the intensity measure of $\Lambda^{(K)}$ is $\theta \dd a \otimes \mu^{(K)}(\dd \ell)$, where $\mu^{(K)}$ is the pushforward of $\mu^{> \varepsilon}$ under the shift $\ell \mapsto \ell + K$. By the tail estimates established in the proof of Lemma~\ref{lemma:diameter_hypothesis_H}, we have, as $v \to \infty$,
$$
\mu^{(K)}((v,\infty)) \,=\, \mu^{> \varepsilon}((v-K, \infty)) \,=\, \frac{1}{v-K} + O(e^{-c(v-K)}) \,=\, \frac{1}{v} + O(v^{-2}),
$$
which gives $\int_1^x \mu^{(K)}((v,\infty)) \dd v = \log x + c' + o(1)$. Hence, $\Lambda^{(K)}$ satisfies Hypothesis~(H). Notice that $N_K(y,x)$ counts \emph{all} the uncovered gaps of $\Lambda^{(K)}$ in $(y,x)$, without any selection. This is the reason why we introduced this new covering: it reduces the problem of estimating a subset of large gaps to counting every gap of $\Lambda^{(K)}$. Since $\sqrt{q} \in (0,1)$, Lemma~\ref{lemma:generating_function_Nst} now applies directly to $N_K(y,x)$ and yields
$$
\mathbb{E} \big[ (\sqrt{q})^{N_K(y,x)} \big] \,\leq\, C x^{-1+\theta}.
$$
Injecting this into \eqref{eq:q_Nstar_NK} and \eqref{eq:_bound_crossing_Nstar} gives the upper bound
$$
\mathbb{P} \left( \{y\} \times \mathbb{S}^1 \overset{\mathcal{L}_\mathfrak{C}^\theta}{\longleftrightarrow} \{x\} \times \mathbb{S}^1 \right) \,\leq\, Cq^{-1/2}\, x^{-1+\theta}.
$$
Combined with the lower bound \eqref{eq:crossing_lower_bound}, this completes the proof of Theorem~\ref{th:crossing_exponent}. \qed

\section{Bounded clusters at super-critical intensities} \label{sec:finiteness_bounded}
In this section, we prove Theorem~\ref{th:n_delta}. We focus only on the super-critical regime $\theta > 1$. We completely abandon the diameter cutoff $\varepsilon$, and we keep denoting $\LC$ as the set of loops that have a duration at least $\tau$, just as in Section~\ref{sec:phase_transition}. The proof strategy relies on splitting the loop soup into separating and non-separating loops. Define, for all $x > 0$, 
\begin{equation} \label{eq:def_r_x}
r(x) := (\log (e+x))^{-1}.    
\end{equation}
Note that $r(x) \in (0,1)$ and that $r(x) \to 0$ as $x \to \infty$. In Section~\ref{sec:finiteness_unvisited}, we show that separating loops form an unbounded connected component which is dense enough to leave almost no empty space.  More precisely, we prove that there are only finitely many unvisited discs of the form $D_\mathfrak{C}((x, e^{i\varphi}), r(x))$, where $x$ is the longitudinal coordinate.
Then, in Section~\ref{sec:tiny_loops}, we show that only finitely many non-separating loops have an inner radius smaller than $(\log (e+x))^{-1}$. This ensures that all but finitely many non-separating loops are too large to fit in the remaining unvisited discs, forcing them to intersect the unbounded cluster.

\subsection{Finiteness of unvisited microscopic discs.} \label{sec:finiteness_unvisited}

Define $\mathcal{N}^{\mathrm{disc}}$ as the maximum number of non-intersecting discs of the form $D_\mathfrak{C}((x,e^{i\varphi}), r(x))$ with $x > 0$ and $\varphi \in [0,2\pi)$ that do not intersect any separating loop from $\LC$ (see Figure~\ref{fig:dense_N} below). By rotational symmetry, and for the sake of simplicity, we will frequently omit the angular component and denote $D_\mathfrak{C}(x,r) := D_\mathfrak{C}((x,e^{i\varphi}), r)$ for discs around the cylinder. This section is dedicated to the proof of the following proposition.

\begin{figure}[h]
    \centering
    \includegraphics[width=11.5cm]{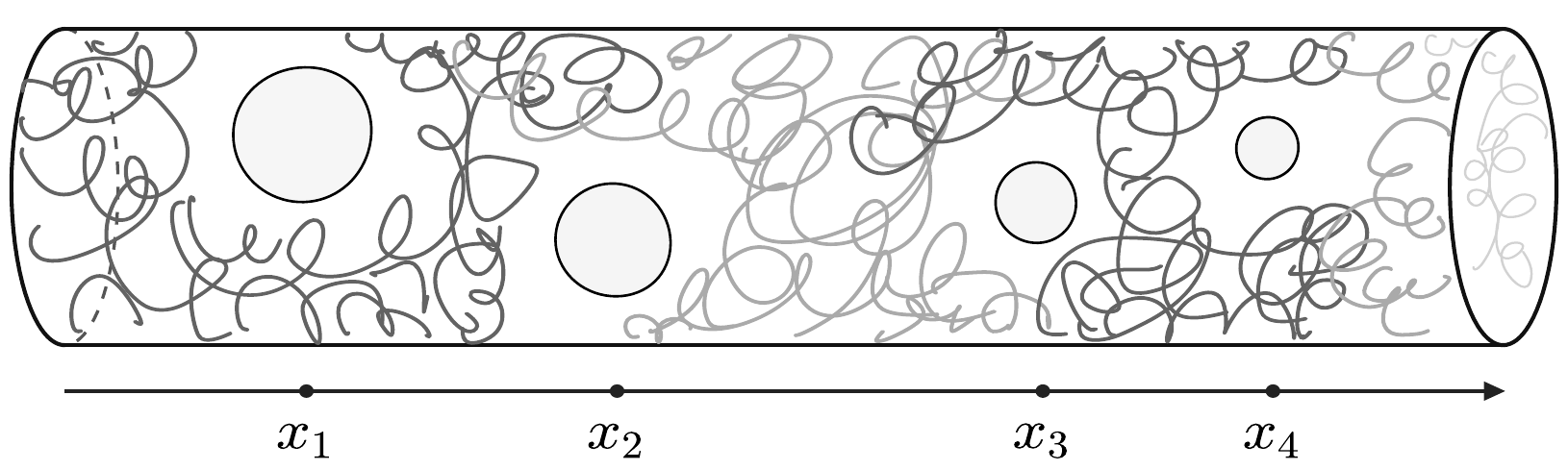}
\caption{A slice far along the cylinder ($x \gg 1$) in the supercritical regime $\theta > 1$. Only the separating loops are represented (in various shades of gray); they all merge into the unique giant unbounded cluster. This giant component leaves empty holes unvisited by any loop. $\mathcal{N}^{\mathrm{disc}}$ denotes the maximal number of disjoint unvisited discs of radius $r(x) = (\log (e+x))^{-1}$, which become smaller as the longitudinal coordinate $x$ increases (here, four such discs centered at $x_1 < x_2 < x_3 < x_4$ are shown).}
    \label{fig:dense_N}
\end{figure}

\begin{proposition} \label{prop:N_log_finite} For all $\theta > 1$, almost surely, $\mathcal{N}^{\mathrm{disc}}$ is finite.
\end{proposition}

Of course, since $r(x)$ is decreasing to zero as $x \rightarrow \infty$, Proposition~\ref{prop:N_log_finite} implies in particular that, almost surely, for all $\delta > 0$, the maximum number of non-intersecting $\delta$-discs on the cylinder that do not intersect any separating loop is finite. Heuristically, the unbounded connected component formed by the separating loops is dense enough to leave almost no empty space. The proof of Proposition~\ref{prop:N_log_finite} relies heavily on estimating the mass of separating loops that hit a given disc on the cylinder.

\begin{lemma} \label{lemma:mu_loop_log_x_hits} We have, 
as $x \rightarrow \infty$:
$$
\mu_{\mathfrak{C}}^{\mathrm{loop}, \tau} \left( \big\{ \wp \text{ hits } D_\mathfrak{C} \big (x,{\textstyle \frac{1}{4}} r(x) \big) \text{ and } \wp \text{ is separating}\big\} \right) \,\geq\, \log x - \log \log \log x + O \big( 1 \big).
$$
\end{lemma}

The factor $1/4$ in the radius is purely technical and does not affect the asymptotic expansion: it is only needed for the discretisation grid in the proof of Proposition~\ref{prop:N_log_finite}. Several intermediate lemmas are needed to prove Lemma~\ref{lemma:mu_loop_log_x_hits}. We recall that the Brownian bubble measures on the disc and on the cylinder have been defined in Section~\ref{sec:bubbles}.

\begin{lemma} \label{lemma:hits_surrounds} There exists a positive constant $C$, such that for all $z=(x,e^{i\theta}) \in \mathfrak{C}$ and $r < \pi$,
$$
\mu_{\mathfrak{C}}^{\mathrm{bub}, \tau} \big( \big\{ \wp \text{ is non-separating and } \wp \text { hits } D_\mathfrak{C}(z, r) \big\} \big) \,\leq\, C e^{-2x/3}.
$$
\end{lemma}

\begin{proof}
Recall that $\bpi : (x,e^{i\theta}) \to x$ denotes the longitudinal projection from the cylinder onto $(0, \infty)$. By definition of the metric $d_\mathfrak{C}$ on the cylinder from Equation~\eqref{eq:def_metric_cyl}, we have $\bpi(D_\mathfrak{C}(z, r)) = (x-r, x+r)$. This implies that any bubble $\wp$ that intersects $D_\mathfrak{C}(z,r)$ has to intersect the cylinder's slice $(x-r, \infty) \times \mathbb{S}^1$. But this implies in turn that $\ell(\wp) > x-r$. 
\begin{align*}
&\mu_{\mathfrak{C}}^{\mathrm{bub}, \tau} \big( \big\{ \wp \text{ is non-separating and } \wp \text { hits } D_\mathfrak{C}(z, r) \big\} \big) \\
&\quad\quad\quad\quad\quad\quad\leq\, \mu_{\mathfrak{C}}^{\mathrm{bub}, \tau} \big( \big\{ \wp \text{ is non-separating and } \ell(\wp) > x-r \big\} \big) \\
&\quad\quad\quad\quad\quad\quad\leq\, Ce^{-\frac{2}{3}(x-r)} \,<\, Ce^{2\pi/3} \, e^{-2x/3}.
\end{align*}
In the last line, we have used Lemma~\ref{lemma:l_contractible} together with the assumption $r < \pi$.
\end{proof}

\begin{lemma} \label{lemma:hits_and_tau}
There exist some positive constants $C,c$, depending on $\tau$, such that for all $z = (x,e^{i\varphi}) \in \mathfrak{C}$ with $x \geq 2$ and $r < 1$,
$$
\mu_\mathfrak{C}^{\mathrm{bub}} \big( \big\{ \wp \text{ hits } D _\mathfrak{C}(z,r) \text{ and } T(\wp) \leq \tau \big\} \big) \,\leq\, C e^{- c x^2}.
$$
\end{lemma}

\begin{proof}
Just as in the proof of the previous lemma, if a bubble $\wp$ hits $D_\mathfrak{C}(z,r)$, then $\ell(\wp) > x-r$. Denoting $(X_t)_{t \geq 0} = \bpi(\wp_t)$ the longitudinal coordinate of $\wp$ and $\zeta$ its lifetime, we get 
\begin{equation} \label{eq:sup_lifetime_projection}
\mu_\mathfrak{C}^{\mathrm{bub}} \big( \big\{ \wp \text{ hits } D_\mathfrak{C}(z,r) \text{ and } T(\wp) \leq \tau \big\} \big) \,\leq\, \mu_\mathfrak{C}^{\text{bub}} \big( \sup_{t \geq 0} X_t > x-r \text{ and } \zeta \leq \tau).
\end{equation}
Integrating the Brownian bubble decomposition in longitudinal and angular components from Section~\ref{sec:bubbles}, Equation~\eqref{eq:bubble_explicit}, we get 
$$
\mu_\mathfrak{C}^{\text{bub}}(\dd X) \,=\, 2\pi p_{\mathbb{S}^1}(\zeta) \mathbf{n}(\dd X),
$$
where $\zeta(X)$ is the lifetime of $X$, $p_{\mathbb{S}^1}$ is the heat kernel on the circle and $\mathbf{n}$ is Itô's excursion measure. Denoting $u = x-r$ and disintegrating Itô's measure on the lifetime $\zeta$, we get
\begin{align} \label{eq:mu_bub_bound}
\mu_\mathfrak{C}^{\text{bub}} \big( \sup_{t \geq 0} X_t > u \text{ and } \zeta \leq \tau) \,=\, \frac{2\pi}{\sqrt{2\pi}} \int_0^\tau \frac{\dd t}{t^{3/2}} p_{\mathbb{S}^1}(t) \mathbf{n}_t \big( \sup_{s \geq 0} X_s > u \big),
\end{align}
where $\mathbf{n}_t$ denotes the law of a Brownian excursion of duration $t$. Now, by Lemma~\ref{lemma:_heat_kernel_bound} there exists a constant $C_\tau > 0$, such that for all $t \leq \tau$, we have $p_{\mathbb{S}^1}(t) \leq C_\tau/\sqrt{t}$. Furthermore, denoting $M = \sup_{s \geq 0} e_s$ where $e$ is an excursion of duration $1$, by Brownian scaling we get
$$
\mu_\mathfrak{C}^{\text{bub}} \big( \sup_{t \geq 0} X_t > u \text{ and } \zeta \leq \tau) \,\leq\, C_\tau' \int_0^\tau \frac{\dd t}{t^2} \mathbb{P}\big(M > \textstyle{\frac{u}{\sqrt{t}}} \big).
$$
Now, by Lemma~\ref{lemma:excur_max_gaussian}, there exists a constant $C' > 0$ such that $\mathbb{P}(M > u) \leq C' e^{-u^2}$. Injecting this bound and swallowing the constants into $C_\tau''$ yields
$$
C_\tau' \int_0^\tau \frac{\dd t}{t^2} \mathbb{P}\big(M > {\textstyle\frac{u}{\sqrt{t}}} \big) \,\leq\, C''_\tau \int_0^\tau \frac{\dd t}{t^2} e^{-\frac{u^2}{t}} \,=\, \frac{C''_\tau}{u^2} \exp \Big( -\frac{u^2}{\tau} \Big).
$$
Injecting this into Equation~\eqref{eq:sup_lifetime_projection} with $u = x-r \geq 1$ (because $x \geq 2$ and $r \leq 1$) yields the desired result.
\end{proof}

We recall (see Section~\ref{sec:bubbles}) that we use a different normalisation of the Brownian bubble measure than the standard one that can be found in \cite{Lawler_2003, Jego_2025}. The next lemma estimates the mass of Brownian bubbles that hit a given disc $B(y,r)$.

\begin{lemma} \label{lemma:hits_D_y_r} Let $y \in (0,1)$ and $r > 0$ such that $B(y,r) \subset \mathbb{D}$. Define the parameters $K = \frac{1+y^2-r^2}{2y}$, $J = \frac{1+r^2-y^2}{2r}$, $a = K - \sqrt{K^2-1}$ and $\rho = J  - \sqrt{J^2 -1}$. Then, as $\rho \to 0$, we have
$$
\mu_{\mathbb{D}}^{\mathrm {bub}}(\{\wp \text { hits } B(y, r)\})=\Big(\frac{1+a^2}{1-a^2}\Big) \frac{1}{|\log \rho|} \big(1 + O(\rho) \big).
$$
\end{lemma}

\begin{proof}
With our normalisation of the Brownian bubble measure, by Lemma~\ref{lemma:bubble_hit_ball_disc}, as $R \to 0$ and uniformly in $\theta \in [0, 2\pi)$, we have
\begin{equation} \label{eq:bubble_rooted_hit_ball}
\mu_{\mathbb{D}}^{\text {bub }}(e^{i \theta})(\{\wp \text { hits } R\mathbb{D} \})= \frac{1+O(R)}{2\pi |\log R|}.
\end{equation}
With the definitions of $a$ and $\rho$ from the lemma's statement, we prove in Appendix~\ref{sec:complex_analysis}, Lemma~\ref{lemma:conformal_circles} that the map $g : \mathbb{D} \to \mathbb{D}$ defined by
\begin{equation} \label{eq:def_g}
g(z) \,:=\, \frac{z-a}{1-az},
\end{equation}
is a conformal automorphism of $\mathbb{D}$ that sends $B(y,r)$ onto $\rho \mathbb{D}$. Hence, by conformal covariance of the Brownian bubble measure and by Equation~\eqref{eq:bubble_rooted_hit_ball}, we get, as $\rho \to 0$ (which is equivalent to $r \to 0$ because a direct calculation shows that $r \sim (1-|y|^2) \rho$ as $\rho \to 0$):
\begin{align}
\mu_{\mathbb{D}}^{\text {bub}}\big(\{\wp \text { hits } B(y, r)\} \big) \,&=\, \int_0^{2\pi} \dd \theta \, |g'(e^{i \theta})|^2\,\mu_{\mathbb{D}}^{\text {bub}}(e^{i \theta})\big(\{\wp \text { hits } \rho \mathbb{D} \} \big) \notag \\
&=\, \frac{1+O(\rho)}{2\pi |\log \rho|} \int_0^{2\pi} |g'(e^{i\theta})|^2 \dd \theta. \label{eq:mu_bub_integrated}
\end{align}
We removed the $(1+O(\rho))$ from the integral as the error is uniform in $\theta$. It now remains to compute the integral $\int_0^{2\pi} |g'(e^{i\theta})|^2 \dd \theta$. By differentiating Equation~\eqref{eq:def_g}, we get $g'(z) = \frac{1-a^2}{(1-az)^2}$. When $z = e^{i\theta}$, identifying the expression with the Poisson kernel in the unit disc yields
$$
|g'(e^{i\theta})| \,=\, \frac{1-a^2}{|1-ae^{i\theta}|^2} \,=\, 2\pi H_\mathbb{D}(a, e^{i\theta}).
$$
On the other hand (see e.g. \cite{stein2003fourier} Section 5.4), we have the following Fourier decomposition $H_\mathbb{D}(a,e^{i\theta}) = \frac{1}{2\pi} \sum_{n \in \mathbb{Z}} a^{|n|} e^{i n \theta}$. By Parseval's identity, we get
$$
4\pi^2 \int_0^{2\pi} H_\mathbb{D}(a, e^{i\theta})^2 \dd \theta \,=\, 2\pi \sum_{n \in \mathbb{Z}} a^{2|n|} \,=\, 2\pi \frac{1+a^2}{1-a^2}.
$$
Injecting the aforementioned facts into Equation~\eqref{eq:mu_bub_integrated} yields the desired result.
\end{proof}

The previous Lemma gives an expression for the mass of Brownian bubbles in $\mathbb{D}$ that hit a given disc. As we are mainly interested in Brownian bubbles on the cylinder, we need a way to compare discs $D_\mathfrak{C}(z,r)$ on the cylinder with their pushback $f^{-1}(D_\mathfrak{C}(z,r))$ on the disc $\mathbb{D}$.

\begin{lemma} \label{lemma:discs_cylinder_disc} Let $f$ be the conformal map from $\mathbb{D} \setminus \{0\}$ onto $\mathfrak{C}$ defined by $re^{i\theta} \mapsto (|\log r|, e^{i\theta})$. Let $z = (x,e^{i\theta}) \in \mathfrak{C}$ and $r < \min(\pi, x)$. We have the following set inclusion
$$
B\big(f^{-1}(z), e^{-x}(1-e^{-r}) \big) \,\subset\, f^{-1}\big( D_\mathfrak{C}(z,r) \big).
$$
\end{lemma}

\begin{proof} See Appendix~\ref{sec:complex_analysis} for a proof of this lemma.
\end{proof}

We can now bound the mass of separating Brownian bubbles on $\mathfrak{C}$ that hit a given disc on the cylinder. We recall that $r(x) = (\log (e+x))^{-1}$ has been defined in Equation~\eqref{eq:def_r_x}.

\begin{lemma} \label{lemma:mu_bub_hits_log} There exists some positive constant $c > 0$ such that, for any $v \geq 2$ and $x > 0$, 
$$
\mu_\mathfrak{C}^{\mathrm{bub}, \tau} \big( \big\{ \wp \text{ hits } D_\mathfrak{C} \big (v,{ \textstyle \frac{1}{4}}r(x) \big) \text{ and } \wp \text{ is separating} \big\} \big) \,\geq\, \frac{1}{v+\log \log (e+x) + c} - O(e^{-2v/3}).
$$
\end{lemma}

\begin{proof}
First, since for any measure $\mu$ and measurable sets $A,B$, we always have $\mu(A\cap B) = \mu(A) - \mu(A\setminus B)$, we get
\begin{align*}
&\mu_\mathfrak{C}^{\mathrm{bub}, \tau} \big( \big\{ \wp \text{ hits } D_\mathfrak{C} \big (v,{ \textstyle \frac{1}{4}}r(x) \big) \text{ and } \wp \text{ is separating} \big\} \big) \\
&\quad\quad\quad=\, \mu_\mathfrak{C}^{\mathrm{bub}, \tau} \big( \wp \text{ hits } D_\mathfrak{C} \big (v,{ \textstyle \frac{1}{4}}r(x) \big)  \big) - \mu_\mathfrak{C}^{\mathrm{bub}, \tau} \big( \big\{ \wp \text{ hits } D_\mathfrak{C} \big (v,{ \textstyle \frac{1}{4}}r(x) \big) \text{ and } \wp \text{ is non-separating} \big\} \big) \\
&\quad\quad\quad\geq\, \mu_\mathfrak{C}^{\mathrm{bub}, \tau} \big(  \wp \text{ hits } D_\mathfrak{C} \big (v,{ \textstyle \frac{1}{4}}r(x) \big) \big) - O(e^{-2v/3}).
\end{align*}
In the last line, we have used Lemma~\ref{lemma:hits_surrounds}. Second, since $\mu_\mathfrak{C}^{\mathrm{bub}, \tau}$ is exactly $\mu_\mathfrak{C}^{\mathrm{bub}}$ restricted to the loops of duration at least $\tau$, that same trick yields
\begin{align*}
&\mu_\mathfrak{C}^{\mathrm{bub}, \tau} \big(  \wp \text{ hits } D_\mathfrak{C} \big (v,{ \textstyle \frac{1}{4}}r(x) \big) \big) \\
&\quad\quad\quad=\, \mu_\mathfrak{C}^{\mathrm{bub}} \big(  \wp \text{ hits } D_\mathfrak{C} \big (v,{ \textstyle \frac{1}{4}}r(x) \big) \big)- \mu_\mathfrak{C}^{\mathrm{bub}} \big(  \big\{ \wp \text{ hits } D_\mathfrak{C} \big (v,{ \textstyle \frac{1}{4}}r(x) \big) \text{ and } T(\wp) \leq \tau \big\}\big) \\
&\quad\quad\quad\geq\, \mu_\mathfrak{C}^{\mathrm{bub}} \big(  \wp \text{ hits } D_\mathfrak{C} \big (v,{ \textstyle \frac{1}{4}}r(x) \big) \big) - O(e^{-cv^2}).
\end{align*}
In the last line, we have used Lemma~\ref{lemma:hits_and_tau}. Putting together the previous bounds, we get
\begin{align}
&\mu_\mathfrak{C}^{\mathrm{bub}, \tau} \big( \big\{ \wp \text{ hits } D_\mathfrak{C} \big (v,{ \textstyle \frac{1}{4}}r(x) \big) \text{ and } \wp \text{ is separating} \big\} \big) \notag \\
&\quad\quad\quad\quad\quad\quad\geq\, \mu_\mathfrak{C}^{\mathrm{bub}} \big(  \wp \text{ hits } D_\mathfrak{C} \big (v,{ \textstyle \frac{1}{4}}r(x) \big) \big) - O(e^{-2v/3}). \label{eq:bub_tau_clean}
\end{align}
Note that Equation~\eqref{eq:bub_tau_clean} makes physical sense: a loop $\wp$ that hits the disc $D_\mathfrak{C}(v,{ \textstyle \frac{1}{4}}r(x))$ has to penetrate far inside the cylinder, making it (exponentially) unlikely that $\wp$ is non-separating or that $\wp$ has a small duration. It only remains to study the $\mu_\mathfrak{C}^{\mathrm{bub}}$-mass (without the cutoff $\tau$) of the loops that intersect the given disc. Let $f$ be the conformal map from $\mathbb{D} \setminus \{0\}$ onto $\mathfrak{C}$ given by $re^{i\theta} \mapsto (|\log r|, e^{i\theta})$. By Lemma~\ref{lemma:discs_cylinder_disc}, we have 
$$
B\big(e^{-v}, e^{-v}(1-e^{- {\scriptstyle \frac{1}{4}} r(x)}) \big) \,\subset\, f^{-1}\big( D_\mathfrak{C}(v,{ \textstyle \frac{1}{4}}r(x)) \big).
$$
Let $\alpha(x) = 1-e^{-\frac{1}{4}r(x)} = 1-e^{-(4\log (e+x))^{-1}}$. Since, by definition, we have $\mu_\mathfrak{C}^{\text{bub}} = f \circ \mu_\mathbb{D}^{\text{bub}}$, we get
\begin{align}
\mu_\mathfrak{C}^{\mathrm{bub}} \big(  \wp \text{ hits } D_\mathfrak{C} \big (v,{ \textstyle \frac{1}{4}}r(x) \big) \big) \,&\geq\, \mu_\mathbb{D}^{\mathrm{bub}} \big(  \wp \text{ hits } B \big (e^{-v}, \alpha(x)e^{-v} \big) \big) \notag \\
&=\, \Big( \frac{1+a^2}{1-a^2} \Big) \frac{1}{|\log \rho|} \big( 1 + O(\rho) \big), \label{eq:bound_bub_disc}
\end{align}
as $\rho \to 0$. In the last line, we have applied Lemma~\ref{lemma:hits_D_y_r} with $y = e^{-v}$, $r = \alpha(x) e^{-v}$ (hence, $a$ and $\rho$ both depend on $v$ and $x$ according to the Lemma's statement). Now, since $a = K-\sqrt{K^2-1} \in (0,1)$, we get, by Equation~\eqref{eq:bound_bub_disc}
\begin{align*}
\mu_\mathfrak{C}^{\mathrm{bub}} \big(  \wp \text{ hits } D_\mathfrak{C} \big (v,{ \textstyle \frac{1}{4}}r(x) \big) \big) \,\geq\, \Big( \frac{1+a^2}{1-a^2} \Big) \frac{1+O(\rho)}{|\log \rho|} \,\geq\, \frac{1}{|\log \rho|} \big( 1 + O(\rho) \big). 
\end{align*}
Denoting the error by $\varepsilon(\rho) = O(\rho)$, we have $|\varepsilon(\rho)| \leq C \rho$ for some positive constant $C$. We get
\begin{align}
\mu_\mathfrak{C}^{\mathrm{bub}} \big(  \wp \text{ hits } D_\mathfrak{C} \big (v,{ \textstyle \frac{1}{4}}r(x) \big) \big) \,\geq\, \frac{1}{|\log \rho|} \big( 1 + \varepsilon(\rho) \big) \,\geq\, \frac{1}{|\log \rho|} - C\frac{\rho}{|\log \rho|}. \label{eq:bound_bub_disc_clean}
\end{align}
By the definition of $J$ from Lemma~\ref{lemma:hits_D_y_r}, some algebraic manipulations yield
$$
2J \,=\, \frac{1 + (\alpha(x)e^{-v})^2 - (e^{-v})^2}{\alpha(x) e^{-v}} \,=\, \frac{e^v-e^{-v}}{\alpha(x)}+\alpha(x) e^{-v} \,\leq\, \frac{e^v}{\alpha(x)} + 1 \,\leq\, \frac{2e^v}{\alpha(x)}.
$$
Before the last bound, we have used the fact that $e^{-v} \alpha(x) \leq 1$, and in the last bound we have used the fact that $1 \leq e^{v}/\alpha(x)$ because $v \geq 0$ and $\alpha(x) \leq 1$. On the other hand, we have
\begin{equation} \label{eq:bound_rho}
\rho \,=\, J-\sqrt{J^2-1} \,=\, \frac{1}{J + \sqrt{J^2-1}} \,\geq\, \frac{1}{2J} \,\geq\, \frac{\alpha(x)}{2e^{v}}
\end{equation}
Now, recalling that $\alpha(x) = 1-e^{-1/(4\log (e+x))}$ and using the inequality $1-e^{-u} \geq \frac{u}{2}$ (for all $u \leq 1$), we get $\alpha(x) \geq \frac{1}{8} \log (e+x)^{-1}$ for all $x > 0$. Hence,
$$
|\log \alpha(x)| \,\leq\, \log \log (e+x) + 3\log 2.
$$
Injecting this into Equation~\eqref{eq:bound_rho} we get,
\begin{equation} \label{eq:bound_log_rho}
|\log \rho| \,\leq\, \log \Big( \frac{2e^v}{\alpha(x)} \Big) \,\leq\, v + \log \log (e+x) + 4\log 2. 
\end{equation}
This concludes the lower bound for $\rho$. We also need an upper bound for the second term in Equation~\eqref{eq:bound_bub_disc_clean}. We have
\begin{equation} \label{eq:upper_bound_rho}
\rho \,=\, J-\sqrt{J^2-1} \,=\, \frac{1}{J + \sqrt{J^2-1}} \,\leq\, \frac{1}{J}
\end{equation}
Now, by definition of $J$, we have
$$
2J \,=\, \frac{e^v-e^{-v}}{\alpha(x)} + \alpha(x)e^{-v} \,\geq\, \frac{e^v-e^{-v}}{\alpha(x)}. 
$$
Since we have supposed that $v > 1$, we have $e^{v}-e^{-v} \geq (1-e^{-2}) e^v$. Furthermore, since $\alpha(x) = 1-e^{-1/(4\log (e+x))} \leq (4\log (e+x))^{-1}$ we get
$$
2J \,\geq\, 4(1-e^{-2})e^v \log (e+x).
$$
Injecting this into Equation~\eqref{eq:upper_bound_rho}, since $x > 0$ we get
$$
\rho \,\leq\, \frac{1}{J} \,\leq\, \frac{2}{1-e^{-2}} \frac{e^{-v}}{4 \log (e+x)} \,\leq\, C' e^{-v}.
$$
Here, we have denoted $C' = \frac{1}{2}(1-e^{-2})^{-1}$. We finally get the following upper bound for the error term,
$$
C \frac{\rho}{|\log \rho|} \,\leq\, C \rho \,\leq\, C C' \frac{e^{-v}}{\log (e+x)} \,\leq\, C_1 e^{-v}
$$
Injecting this together with the bound of Equation~\eqref{eq:bound_log_rho} into Equation~\eqref{eq:bound_bub_disc_clean} finally yields, for all $x > 0$ and $v \geq 2$,
\begin{align*} \label{eq:exact_bound_mu_hits_log}
\mu_\mathfrak{C}^{\mathrm{bub}} \big(  \wp \text{ hits } D_\mathfrak{C} \big (v,{ \textstyle \frac{1}{4}}r(x) \big) \big)  \,&\geq\, \frac{1}{|\log \rho|} - C \frac{\rho}{|\log \rho|} \\
&\geq\, \frac{1}{v+\log \log (e+x) + c} - O(e^{-v}),
\end{align*}
where we have set $c = 4\log 2 > 0$.  Injecting this into Equation~\eqref{eq:bub_tau_clean} concludes the proof.
\end{proof}

\noindent \textbf{Proof of Lemma~\ref{lemma:mu_loop_log_x_hits}.} Denoting $\tau_u : \mathfrak{C} \mapsto \mathfrak{C}$ given by $\tau_u((x,e^{i\theta})) = (x+u, e^{i\theta})$, we recall (see Equation~\eqref{eq:loop_bubble_decomp}) that the Brownian loop measure on the cylinder disintegrates with respect to the left-endpoint of the loops as
\begin{equation} \label{eq:rooted_soup}
\mu_\mathfrak{C}^{\text{loop}} \,=\, \int_0^\infty \dd u \, \tau_u \circ \mu_{\mathfrak{C}}^{\text{bub}}.
\end{equation}
Let $x > 2$, since $\mu_{\mathfrak{C}}^{\text{bub}}$ is a positive measure, we can stop the integral at $x-2$ to get a lower bound. Translating the events by $-u$ yields,
\begin{align} \label{eq:rooted_loop_log_x}
&\mu_{\mathfrak{C}}^{\mathrm{loop}, \tau} \left( \big\{ \wp \text{ hits } D_\mathfrak{C} \big (x,{ \textstyle \frac{1}{4}}r(x) \big) \text{ and } \wp \text{ is separating}\big\} \right) \\
&\qquad\quad \geq \, \int_0^{x-2} \dd u\, \mu_\mathfrak{C}^{\text{bub}, \tau} \big( \big\{ \wp \text{ hits } D_\mathfrak{C} \big (x-u,{ \textstyle \frac{1}{4}}r(x) \big) \text{ and } \wp \text{ is separating}\big\} \big) \notag.
\end{align}
Now, injecting the resulting lower bound of Lemma~\ref{lemma:mu_bub_hits_log} into Equation~\eqref{eq:rooted_loop_log_x} (with $v = x-u > 2$), the integral of $O(e^{-2v/3})$ gives a $O(1)$, and
\begin{align*}
&\mu_{\mathfrak{C}}^{\mathrm{loop}, \tau} \left( \big\{ \wp \text{ hits } D_\mathfrak{C} \big (x,{ \textstyle \frac{1}{4}}r(x) \big) \text{ and } \wp \text{ is separating}\big\} \right) \\
&\qquad\qquad\qquad\geq\, \int_2^{x} \frac{\dd v}{v + \log \log (e+x) + c} - O(1) \\
&\qquad\qquad\qquad=\, \log \left(x+\log \log (e+x)+c\right)-\log \left(2+\log \log (e+x)+c\right) - O(1) \\
&\qquad\qquad\qquad=\, \log x - \log \log \log x + O(1).
\end{align*}
The proof is complete.
$\hfill \square$ \newline

For the sake of readability, define $\LCtc$ as the set of loops $\wp \in \LC$ that are separating. Recall that $\mathcal{N}^{\mathrm{disc}}$ is defined as the maximum number of non-intersecting discs of the form $D_\mathfrak{C}((x,e^{i\varphi}), r(x))$ with $x > 0$ and $\varphi \in [0, 2\pi)$ that do not intersect any loop from $\LCtc$. \newline

\noindent \textbf{Proof of Proposition~\ref{prop:N_log_finite}}
Let $B_n$ be the event that there is a hole of radius $r(x)$ in the soup $\LCtc$ between longitudinal coordinates $n$ and $n+1$, i.e.,
$$
B_n \,=\, \left\{ \exists y=(x,e^{i \varphi}) \in [n,n+1) \times \mathbb{S}^1 :\, D_\mathfrak{C}(y,r(x)) \cap \LCtc = \emptyset \right\}
$$
We know that $\mathcal{N}^{\mathrm{disc}}$ is finite if and only if $B_n$ occurs only finitely many times (because, there can be at most $\asymp (\log n)^2$ non intersecting such discs inside each slice $[n,n+1) \times \mathbb{S}^1$). Now, discretising the cylinder with smaller balls of radius $\delta_n = \frac{1}{4} r(n)$ (for any disc $D_\mathfrak{C}(z, r(x))$ centered at $z = (x, e^{i\varphi}) \in [n, n+1) \times \mathbb{S}^1$, there exist grid points $i, j$ such that $D_\mathfrak{C}\big((n + \delta_n i, \delta_n j), \delta_n\big) \subset D_\mathfrak{C}(z, r(x))$), we get
$$
\mathbb{P}(B_n) \,\leq\, \mathbb{P} \left( \exists i \in \{0, \dots \lfloor \delta_n^{-1} \rfloor \}, \exists j \in \{0, \dots \lfloor 2\pi \delta_n^{-1} \rfloor \} : D_\mathfrak{C}\big((n+\delta_n i, \delta_n j), \delta_n \big) \cap \LCtc = \emptyset \right).
$$
By a union bound over all grid points $i \in \{0, \dots, \lfloor \delta_n^{-1} \rfloor \}$ and $j \in \{0, \dots, \lfloor 2\pi \delta_n^{-1} \rfloor \}$, we get,
$$\mathbb{P}(B_n) \,\leq\, \sum_{i=0}^{\lfloor \delta_n^{-1}\rfloor} \sum_{j=0}^{\lfloor 2\pi \delta_n^{-1}\rfloor} \mathbb{P}(A_{i,j}^{n}),$$
where $A_{i,j}^{n}:= \{D_\mathfrak{C}\big((n+\delta_n i, \delta_n j), \delta_n \big) \cap \LCtc = \emptyset\}$. By rotational symmetry, $\mathbb{P}(A_{i, j}^n)$ does not depend on the discretised angle $j$. Furthermore, this probability is non-increasing in $i$ (this comes directly by the fact that $\mathcal{L}_\mathfrak{C}^\theta$ is Poissonian and the Bubble decomposition of $\mu_\mathfrak{C}^{\mathrm{loop}}$ in \eqref{eq:loop_bubble_decomp}). Hence, recalling that $\delta_n = \frac{1}{4}r(n) = \frac{1}{4} (\log (e+n))^{-1}$, we get $2\pi \delta_n^{-2} \leq C (\log n)^2$, hence
\begin{equation} \label{eq:bound_P_Bn}
\mathbb{P}(B_n) \,\leq\, C (\log n)^2 \,\mathbb{P} \big( D_\mathfrak{C}(n,\delta_n) \cap \LCtc = \emptyset \big).    
\end{equation}
On the other hand, by Lemma~\ref{lemma:mu_loop_log_x_hits},
\begin{align*}
\mathbb{P} \big( D_\mathfrak{C}(n,\delta_n) \cap \LCtc = \emptyset \big) \,&=\, e^{-\theta \mu_{\mathfrak{C}}^{\mathrm{loop}, \tau} \big( \{ \wp \text{ hits } D_\mathfrak{C}(n,\frac{1}{4}r(n)) \text{ and } \wp \text{ is separating}\} \big)} \\
&\leq e^{- \theta \big(\log n - \log\log\log n + O(1)\big)} \,=\, n^{-\theta} (\log\log n)^\theta e^{O(1)}.
\end{align*} 
Injecting this into Equation~\eqref{eq:bound_P_Bn}, we get
$$
\mathbb{P}\left(B_n\right) \,\leq\, O\big(n^{-\theta}(\log n)^2(\log \log n)^\theta\big) .
$$
Applying the Borel-Cantelli lemma concludes the proof. 
$\hfill \square$

\begin{remarque} \label{rem:fastest_decay}
While Proposition~\ref{prop:N_log_finite} is stated with the logarithmic radius $r(x) \asymp (\log x)^{-1}$ for convenience, the proof adapts to much faster-decaying functions $r(x)$. Let us briefly track the dependence on an arbitrary radius $r(x) \to 0$. First, replacing $(\log(e+x))^{-1}$ by $r(x)$ in Equation~\eqref{eq:bound_log_rho} gives
$$
|\log \rho| \,\leq\, v + \log\big({\textstyle \frac{1}{r(x)}}\big) + c.
$$
Integrating this bound over $v \in [2, x]$ as in Lemma~\ref{lemma:mu_loop_log_x_hits}, the mass of hitting loops becomes
$$
\mu_{\mathfrak{C}}^{\mathrm{loop}, \tau}\big(\big\{\wp \text{ hits } D_{\mathfrak{C}}(x, {\textstyle \frac{1}{4}} r(x)) \text{ and is separating}\big\}\big) \,\geq\, \log x - \log \log \big({\textstyle \frac{1}{r(x)}}\big) + O(1).
$$
Next, in the proof of Proposition~\ref{prop:N_log_finite}, covering each slice $[n,n+1) \times \mathbb{S}^1$ requires $O(r(n)^{-2})$ grid balls of radius $\delta_n \asymp r(n)$. The probability that a single ball is avoided by all separating loops satisfies
$$
\mathbb{P} \big( D_\mathfrak{C}(n,\delta_n) \cap \LCtc = \emptyset \big) 
\,\leq\, e^{- \theta \big( \log n - \log \log ( \frac{1}{r(n)} ) + O(1) \big)} 
\,\leq\, C n^{-\theta} \left( \log \Big( \frac{1}{r(n)} \Big) \right)^\theta.
$$
Summing over all grid points, the Borel--Cantelli lemma applies as soon as
\begin{equation} \label{eq:borel_cantelli_condition}
\sum_{n \geq 1} \frac{n^{-\theta}}{r(n)^2} \left( \log \Big( \frac{1}{r(n)} \Big) \right)^\theta \,<\, \infty.    
\end{equation}
This condition reveals that the giant cluster formed by separating loops is extremely dense:  It easily absorbs any logarithmically  decaying radii: \eqref{eq:borel_cantelli_condition} converges for $r(x) = (\log x)^{-\beta}$ for any $\beta > 0$. More remarkably, it even rules out polynomially shrinking radii: taking $r(x) = x^{-\alpha}$, the general term is of order $n^{2\alpha - \theta}(\log n)^\theta$, which is summable for any power $\alpha < \frac{\theta-1}{2}$. Finally, up to a first order logarithmic correction, the optimal decay given by this method is
$$
r_{\mathrm{opt}}(x) \,\asymp\, \frac{(\log x)^{\frac{\theta+1}{2} + \delta}}{x^{\frac{\theta-1}{2}}},
$$
for any $\delta > 0$. In other words, as $x \to \infty$, the giant cluster leaves no empty space: even microscopic holes that shrink polynomially fast along the cylinder are almost surely entirely filled. (See Figure~\ref{fig:dense_N})
\end{remarque}

\subsection{Ruling out loops with tiny inner radius.} \label{sec:tiny_loops}

In this section, we prove Theorem~\ref{th:n_delta}. The proof relies on bounding the number of non-separating loops that have a very small inner radius. First, we define the inner radius: Let $\eta : [0,t] \to \mathfrak{C}$ be a continuous path on the cylinder. We define $\mathrm{Int}(\eta)$ as the complement of the unique unbounded connected component of $\mathfrak{C} \setminus \eta([0,t])$. Then,  the inner radius of $\eta$ is defined as the radius $r$ of the largest disc $D_{\mathfrak{C}}(x,r)$ that is contained in $\mathrm{Int}(\eta)$, in other words
$$
\mathrm{inrad}(\eta) \,:=\, \sup_{x \in \mathrm{Int}(\eta)} \mathrm{d}_{\mathfrak{C}}\big(x, \partial \mathrm{Int}(\eta)\big).
$$
Where $\mathrm{d}_{\mathfrak{C}}$ is the canonical metric on the cylinder defined in Equation~\eqref{eq:def_metric_cyl}.  If $\eta$ is non-separating, its inner radius is smaller than $\pi$. Furthermore, it is clear that if $\mathrm{inrad}(\eta) > \delta$, then $\eta$ cannot fit inside any disc of radius $r \leq \delta$.

\begin{lemma} \label{lemma:inrad} Denote $\mathbb{P}^{t, z, z}$ the law of a Brownian loop of duration $t$ in the complex plane $\mathbb{C}$ started at $z$. There exist positive constants $C$ and $c$ such that, for all $t \geq 0$, $z \in \mathbb{C}$ and $\delta > 0$,
$$\mathbb{P}^{t, z, z}\big(\mathrm{inrad}(\wp) \leq \delta \big) \,\leq\, Ce^{-ct \delta^{-2}}.
$$
\end{lemma}

\begin{proof} By translation symmetry, we can set $z = 0$. Let $\wp$ be a Brownian loop of duration $t$ started at $0 \in \mathbb{C}$. On the time interval $[0, \frac{t}{2}]$, the law of $(\wp_s)_{s \in [0,t/2]}$ is absolutely continuous with respect to the law of a planar standard Brownian motion $(B_s)_{s \in [0,t/2]}$ with the Radon-Nikodym derivative $D_{t/2}$ given by Bayes' formula (see e.g., \cite{revuz1999continuous} Chapter XI),
$$
D_{t/2} \,=\, \frac{p_{\mathbb{C}}(t/2; B_{t/2}, 0)}{p_{\mathbb{C}}(t,0,0)} \,=\, 2 \exp \Big(-\frac{\|B_{t / 2}\|^2}{t}\Big) \,\leq\, 2 \quad \text{a.s}.
$$
By monotonicity, the inner-radius of $\wp([0,s])$ is non-decreasing in $s \geq 0$. Hence,
\begin{equation} \label{eq:loop_radon_bm}
\mathbb{P}^{t, z, z}\big(\mathrm{inrad}(\wp) \leq \delta \big) \,\leq\, 2 \mathbb{P}(\mathrm{inrad}(B_{[0,t/2]}) \leq \delta).
\end{equation}
By Brownian scaling, we can set $\delta = 1$ and consider the new time $t' = t/(2\delta^2)$. We proceed by cutting the time interval $[0,t']$ into $N = \lfloor t' \rfloor$ sub-intervals $[k, k+1]$ for $k = 0,\dots, N-1$. The Markov property ensures that on each slice, the Brownian motion has a probability $p > 0$ (independent of $t$ and $\delta$) of forming a loop outside the disc centered around $B_k$ and of radius $1$. (This is because it is an open event for the uniform convergence topology of $\mathcal{C}([0,1], \mathbb{C})$, so its Wiener measure is positive). Finally, on the event $\{\mathrm{inrad}(B_{[0,t']}) \leq 1\}$, the $N$ attempts to create a loop have to fail. Hence,
$$
\mathbb{P} \big( \mathrm{inrad}(B_{[0,t/2]}) \leq \delta \big) \,\leq\, (1-p)^{\lfloor t/(2\delta^2) \rfloor} \,\leq\, C e^{-ct \delta^{-2}}.
$$
We conclude by injecting the previous bound into Equation~\eqref{eq:loop_radon_bm}.
\end{proof}

Recall the definition of $r(x) = (\log (e+x))^{-1}$ from Section~\ref{sec:finiteness_unvisited}. Denote $\mathcal{N}^{\mathrm{tiny}}$ the number of loops $\wp \in \LC$ that are non-separating and that have an inner radius smaller than $r(a(\wp))$, where $a(\wp) = \min \bpi(\wp)$ is the left-longitudinal endpoint of the loop, i.e.
\begin{equation} \label{eq:def_N_tiny}
\mathcal{N}^{\mathrm{tiny}} \,:=\, \# \big\{ \wp \in \LC : \wp \text{ is non-separating and } \mathrm{inrad}(\wp) \leq r(a(\wp)) \big\}.    
\end{equation}
The proof of Theorem~\ref{th:n_delta} relies on the following lemma.

\begin{lemma} \label{lemma:finite_delta_n_number} 
    For all $\theta > 0$, almost surely, $\mathcal{N}^{\mathrm{tiny}}$ is finite.
\end{lemma}

\begin{proof} For $n \geq 1$, since $r(x) = (\log (e+x))^{-1}$ is decreasing, we have $r(a(\wp)) \leq r(n)$ for all loops $\wp$ whose left endpoint satisfies $a(\wp) \in [n,n+1)$. Hence, we denote $B_n$ the event
$$
B_n \,:=\, \big\{ \exists \wp \in \LC :\, a(\wp) \in [n,n+1), \wp \text{ is non-separating } \text{ and } \mathrm{inrad}(\wp) \leq r(n) \big\}.
$$
Denote $\mathcal{M}_n$ the set of non-separating loops of duration $T(\wp) \geq \tau$ whose left endpoint satisfies $a(\wp) \in [n,n+1)$ and whose inner radius is smaller than $r(n)$. Since $\LC$ is a Poisson point process of intensity $\theta \mu_\mathfrak{C}^{\text{loop}}$, we get
\begin{align} \label{eq:bound_P_Bn_10}
\mathbb{P}(B_n) \,=\, 1-e^{-\theta \mu_\mathfrak{C}^{\text{loop}} (\mathcal{M}_n)} \,\leq\, \theta  \mu_\mathfrak{C}^{\text{loop}} (\mathcal{M}_n). 
\end{align}
We now bound the $\mu_\mathfrak{C}^{\text{loop}}$-mass of the set $\mathcal{M}_n$. In order to apply Lemma~\ref{lemma:inrad}, we split the loops of $\mathcal{M}_n$ based on their longitudinal length:
$$
\mathcal{M}_n \,=\, \big( \mathcal{M}_n \cap \{\ell(\wp) > n^2\}\big) \cup \big( \mathcal{M}_n \cap \{\ell(\wp) \leq n^2\}\big).
$$
By the Brownian bubble decomposition of $\mu_\mathfrak{C}^{\mathrm{loop}}$ in \eqref{eq:loop_bubble_decomp} and the asymptotic \eqref{eq:pointwise_mu_tail} established in Lemma~\ref{lem:ppp_naive_projection}, we get
\begin{align} \label{eq:bound_mn_1}
\mu_{\mathfrak{C}}^{\text {loop }}\big(\mathcal{M}_n \cap\{\ell(\wp)>n^2\}\big) \,\leq\, \int_n^{n+1} \dd a\, \mu_{\mathfrak{C}}^{\text {bub}, \tau}(\ell(\wp)>n^2) \,=\, \frac{1}{n^2}+O(e^{-c n^2}).    
\end{align}
Lifted to the complex plane $\mathbb{C}$, any loop on the cylinder $\mathfrak{C}$ started at $z = (x,\varphi)$ is a Brownian bridge on $\mathbb{C}$ started at $x + i \varphi$ and ending at $x + i \varphi + 2i k \pi$ where $k$ is the random winding number of $\wp$. In particular, on the event where $\wp$ is non-separating, we have $k = 0$, so $\wp$ is a real Brownian loop on $\mathbb{C}$. By the restriction property of the loop measure \cite{lawler2003brownianloopsoup}, we get
$$
\mu_{\mathfrak{C}}^{\text{loop}}\big(\mathcal{M}_n \cap \{\ell(\wp) \leq n^2\} \big) \,\leq\, \int_{[n,n^2+n+1) \times [0,2\pi)} \dd z \int_\tau^\infty \frac{\dd t}{t}p_\mathbb{C}(t; z, z) \mathbb{P}_\mathbb{C}^{t, z, z} \big(\mathrm{inrad}(\wp) \leq r(n) \big).
$$
The previous bound comes from the fact that if a loop $\wp$ satisfies $a(\wp) \in [n,n+1)$ and $\ell(\wp) \leq n^2$ then the starting point $z$ of the lifted loop $\tilde{\wp}$ belongs to $[n,n^2+n+1] \times [0, 2\pi)$. We also used the fact that for any non-separating loop $\wp$, we have $\mathrm{inrad}_\mathfrak{C}(\wp) \geq \mathrm{inrad}_\mathbb{C}(\tilde{\wp})$ and $r(x)$ is decreasing. By the small ball estimate of Lemma~\ref{lemma:inrad}, and by translation invariance in $\mathbb{C}$, integrating along the starting position of the loop (whose area is $2 \pi(n^2+1)$) yields,
\begin{align}
\mu_{\mathfrak{C}}^{\text{loop}}\big(\mathcal{M}_n \cap \{\ell(\wp) \leq n^2\} \big) \,&\leq\,  2\pi (n^2+1)\int_\tau^\infty \frac{dt}{t} \frac{1}{2\pi t} \mathbb{P}_\mathbb{C}^{t,0,0} (\mathrm{inrad}(\wp) \leq r(n)) \notag \\
&\leq\,  2\pi(n^2+1) C e^{-c \tau \,r(n)^{-2}} \int_\tau^\infty \frac{dt}{t} \frac{1}{2\pi t} \notag \\
&=\, M_\tau (n^2+1) e^{-c\tau\,r(n)^{-2}}. \label{eq:bound_mn_2}
\end{align}
We have denoted $M_\tau =  C \int_\tau^\infty \frac{dt}{t^2}  < \infty$. Recalling that $r(n) = \log(e+n)^{-1}$, writing $r(n)^{-2} = \log(e+n) \times \log (e+n)$ and injecting the bounds \eqref{eq:bound_mn_1} and \eqref{eq:bound_mn_2} into Equation~\eqref{eq:bound_P_Bn_10}, we get
$$
\mathbb{P}(B_n) \,\leq\, \theta \big(n^{-2} + O(e^{-cn^2})+ M_\tau(n^2+1) (e+n)^{-c\tau \log (e+n)}\big) \,=\, O\big(n^{-2} \big).
$$
In the last line, we have used the fact that for $n$ sufficiently large, $c \tau \log (e+n) \geq 4$. Hence, by the Borel-Cantelli lemma, almost surely, $B_n$ only occurs finitely many times. In addition, since $\mu_\mathfrak{C}^{\mathrm{bub}, \tau}$ is a finite measure, almost surely, there are only finitely many loops $\wp \in \LC$ satisfying $a(\wp) \in [n,n+1)$. Hence, almost surely, $\mathcal{N}^{\mathrm{tiny}}$ is finite.
\end{proof}

\noindent \textbf{Theorem~\ref{th:n_delta}.} 
For all $\theta > 1$, almost surely, $\LC$ contains only finitely many bounded clusters.

\begin{proof}
First, we know from the proof of Theorem~\ref{th:phase_transition_loopsoup} that, almost surely, the separating loops form a unique unbounded cluster (this is because the covering $\Lm$, which is composed only of core intervals of separating loops, contains an unbounded connected component when $\theta > 1$). Hence there exists $X < \infty$ such that $\Lm$ covers $(X, \infty)$. Since there are only finitely many separating loops that satisfy $a(\wp) < X$, there can only be a finite number of bounded connected components formed by the set of separating loops.

Now, let $\{\wp_n\}_{n \geq 0}$ denote the set of non-separating loops. For any loop $\wp$, recall that $a(\wp) = \min \bpi(\wp)$ is defined as its left longitudinal endpoint on the cylinder. We divide these non-separating loops into two categories based on their inner radius. 

By Lemma~\ref{lemma:finite_delta_n_number}, almost surely, only a finite number $\mathcal{N}^{\mathrm{tiny}}$ of loops in $\{\wp_n\}_{n \geq 0}$ satisfy $\mathrm{inrad}(\wp_n) \leq r(a(\wp_n))$. Since there are only finitely many such microscopic loops, they can only create a finite number of isolated clusters.

All infinitely many remaining non-separating loops satisfy $\mathrm{inrad}(\wp) > r(a(\wp))$. Suppose by contradiction that they form infinitely many distinct bounded connected components on the cylinder. Since each compact slice of $\mathfrak{C}$ contains only finitely many loops of $\mathcal{L}_{\mathfrak{C}}^{\theta, \tau}$, these clusters cannot accumulate locally or be infinitely nested. Their longitudinal coordinates must grow to infinity. So we can extract an infinite sequence of disjoint, non-nested, non-separating loops $(\wp_k)_{k \ge 1}$ satisfying $\mathrm{inrad}(\wp_k) > r(a(\wp_k))$ in distinct clusters.

For each $k$, by definition of the inner radius, $\mathrm{Int}(\wp_k)$ contains a disc of radius $r(a(\wp_k))$ centered at some point $z_k = (x_k, e^{i\varphi_k})$. Since $x_k \ge a(\wp_k)$ and $r$ is decreasing, we have $r(x_k) \le r(a(\wp_k))$, which ensures that $D_{\mathfrak{C}}(z_k, r(x_k)) \,\subset\, \mathrm{Int}(\wp_k)$. Because the loops $\wp_k$ are disjoint and non-nested, these discs are also disjoint. In addition, since these clusters are disconnected from the separating loops, none of these discs can intersect any separating loop. 

This gives an infinite family of disjoint discs of the form $D_{\mathfrak{C}}((x, e^{i\varphi}), r(x))$ avoiding all separating loops, which directly contradicts the fact that $\mathcal{N}^{\mathrm{disc}}$ is finite from Proposition~\ref{prop:N_log_finite}. Therefore, all but a finite number of the non-separating loops must intersect, and thus coagulate with, one of the finitely many clusters formed by the separating loops. We conclude that the total number of bounded clusters is almost surely finite.
\end{proof}

\section{Excursions from the boundary of a cluster.} \label{sec:excursions_boundary_proof}

This section is dedicated to the proof of Theorem~\ref{th:intensity}. Let $\mathcal{L}_\mathbb{D}^\theta$ be a Brownian loop soup in the unit disc of intensity $\theta \in (0,\frac{1}{2}]$. Let $\partial O$ be the outer boundary of the outermost cluster surrounding the origin and $O$ be the inner (open) domain delimited by $\partial O$. Conditionally on $\partial O$, define $\mathcal{L}_O$ as the subset of loops $\wp \in \mathcal{L}_\mathbb{D}^\theta$ contained in $O$ and $\mathcal{L}_{\mathbb{D} \setminus O}$ be the loops that intersect $\mathbb{D} \setminus O$. This almost surely defines a partition of the Brownian loop soup. By \cite{Qian_2019,qian2018lawpointprocessbrownian}, conditionally on $O$, $\mathcal{L}_{O}$ is a fresh Brownian loop soup in $O$ and $\mathcal{L}_{\mathbb{D} \setminus O}$ is composed of the loops outside of $O$ together with the loops that touch $\partial O$ (which compose the excursions from the boundary of the cluster). We define the occupation time field of the outer loops $\mathcal{T}$ as 
\begin{equation} \label{eq:def_occupation_O}
\mathcal{T}(f) \:=\; \sum_{\wp \in \mathcal{L}_{\mathbb{D} \setminus O}} \int_0^{T(\wp)} f(\wp_s) \dd s.    
\end{equation}
Note that there is no renormalisation in the above display. Hence, almost surely, for all non-negative test functions $f$, we have $\mathcal{T}(f) = \infty$ if and only if $\mathrm{supp}(f) \cap \overline{O}^c \neq \emptyset$. If $\chi$ is a measure with Lebesgue density $f$, we simply write $\mathcal{T}(\chi)$ for $\mathcal{T}(f)$. 

\begin{lemma} \label{lemma:lim_T_lejan} For $r > 0$, let $\pi_r$ be the uniform measure on $B(0,r)$. Then,
$$
\lim_{r \to 0} \mathbb{E} \big[ \mathcal{T}(\pi_r) \mathbf{1}_{d(0, \partial O) > r} \big] \,=\, -\frac{\theta}{\pi} \mathbb{E} \big[ \log \CR(0, O) \big].
$$
\end{lemma}

\begin{proof}
    We introduce the event of interest $E_r = \{d(0, \partial O ) > r\}$. Since $\partial O$ is the outer boundary of the outermost cluster surrounding the origin, $\partial O$ is almost surely at positive distance from the origin. Hence $\mathbf{1}_{E_r} \to 1$ almost surely as $r \to 0$. For $\mathfrak{L}$ a configuration of loops in the disc, and for $\varepsilon > 0$, we define the (non-normalized, $\varepsilon$-truncated) occupation time field $\widetilde{\mathcal{L}}^\varepsilon(\mathfrak{L})$ of $\mathfrak{L}$ as 
\begin{equation}
\langle \widetilde{\mathcal{L}}^\varepsilon(\mathfrak{L}), f \rangle \,:=\, \sum_{\wp \in \mathfrak{L}} \mathbf{1}_{T(\wp) > \varepsilon} \int_0^{T(\wp)} f(\wp_s) \dd s.
\end{equation}
Since $\mathcal{L}_{O} \cup \mathcal{L}_{\mathbb{D} \setminus O}$ forms a partition of $\mathcal{L}_\mathbb{D}^\theta$, we get that, conditionally on $O$, almost surely, for all test functions $f$,
\begin{align} \label{eq:def_T_eps}
    \mathcal{T}(f) \,:=\, \lim_{\varepsilon \to 0} \mathcal{T}^\varepsilon(f) \,&=\, \lim_{\varepsilon \to 0} \left( \langle \widetilde{\mathcal{L}}^\varepsilon(\mathcal{L}_\mathbb{D}^\theta), f  \rangle - \langle \widetilde{\mathcal{L}}^\varepsilon(\mathcal{L}_O), f \rangle \right). 
\end{align}
We have defined $\mathcal{T}^\varepsilon$ as the term inside the right-hand side $\varepsilon$-limit. Since $\mathcal{T}^\varepsilon(f)$ is also equal to the sum of Equation~\eqref{eq:def_occupation_O} restricted to the loops of duration at least $\varepsilon$, it is also non-increasing in $\varepsilon$. By the monotone convergence theorem, we get, for all $r > 0$,
\begin{equation} \label{eq:montone_convergence_T}
\mathbb{E} [\mathcal{T}(\pi_r) \mathbf{1}_{E_r}] \,=\, \lim_{\varepsilon \to 0} \mathbb{E} [\mathcal{T}^\varepsilon(\pi_r) \mathbf{1}_{E_r}].    
\end{equation}
Setting $f = \pi_r$, we compute the expected value of the two terms inside $\mathcal{T}^\varepsilon(\pi_r)$ in \eqref{eq:def_T_eps}. The event $E_r$ is measurable with respect to $\sigma(O)$, furthermore, by \cite{Qian_2019}, conditionally on $O$, the set of loops inside $O$ is a Brownian loop soup in $O$ of intensity $\theta$. Since $\pi_r$ is deterministic, Campbell's formula applies conditionally on $O$ and the term involving $\mathcal{L}_O$ becomes
\begin{equation} \label{eq:decomp_O_local}
    \mathbb{E} \left[ \langle \widetilde{\mathcal{L}}^\varepsilon(\mathcal{L}_O), \pi_r \rangle \mathbf{1}_{E_r} \,\middle|\, O \right]  \,=\, \theta  \mathbf{1}_{E_r} \int_O \pi_r(\dd x) G_O^{\varepsilon}(x,x),
\end{equation}
almost surely. For a simply connected domain $D$, we have set $G_D^\varepsilon(x,x) = \int_\varepsilon^\infty p_D(t,x,x) \dd t$ the truncated Green's function in $D$. Writing $\mathbf{1}_{E_r} = \mathbf{1} - \mathbf{1}_{E_r^c}$, the term involving $\mathcal{L}_{\mathbb{D}}^\theta$ becomes
\begin{align} \label{eq:decomp_local_times_D}
    \mathbb{E} \left[  \langle \widetilde{\mathcal{L}}^\varepsilon(\mathcal{L}_\mathbb{D}^\theta), \pi_r  \rangle \mathbf{1}_{E_r} \right] \,&=\, \mathbb{E} \left[  \langle \widetilde{\mathcal{L}}^\varepsilon(\mathcal{L}_\mathbb{D}^\theta), \pi_r  \rangle \right] - \mathbb{E} \left[  \langle \widetilde{\mathcal{L}}^\varepsilon(\mathcal{L}_\mathbb{D}^\theta), \pi_r  \rangle \mathbf{1}_{E_r^c} \right] 
\end{align}
Since $\pi_r$ is a deterministic function and $\mathcal{L}_\mathbb{D}^\theta$ is a Brownian loop soup of intensity $\theta$ in $\mathbb{D}$, Campbell's formula yields
\begin{align}
\mathbb{E} \left[  \langle \widetilde{\mathcal{L}}^\varepsilon(\mathcal{L}_\mathbb{D}^\theta), \pi_r  \rangle \right] \,&=\, \theta \int_\mathbb{D} \pi_r(\dd x) G_\mathbb{D}^\varepsilon(x,x)     \\
&=\, \theta \mathbb{E} \Big[\mathbf{1}_{E_r} \int_O \pi_r(\dd x) G_\mathbb{D}^\varepsilon(x,x) \Big] + \theta \mathbb{E} \Big[\mathbf{1}_{E_r^c} \int_\mathbb{D} \pi_r(\dd x) G_\mathbb{D}^\varepsilon(x,x)  \Big]    \notag
\end{align}
In the last line, we have used the trivial identity $1 = \mathbb{E}[\mathbf{1}_{E_r} + \mathbf{1}_{E_r^c}]$, and the fact that on the event $E_r$, $\mathrm{supp} (\pi_r) = B(0,r) \subset O$, so the integral over $\mathbb{D}$ coincides with the integral on $O$. Injecting the previous equality into Equation~\eqref{eq:decomp_local_times_D} yields
\begin{align} 
    \mathbb{E} \left[  \langle \widetilde{\mathcal{L}}^\varepsilon(\mathcal{L}_\mathbb{D}^\theta), \pi_r  \rangle \mathbf{1}_{E_r} \right] \,&=\, \theta \mathbb{E} \Big[\mathbf{1}_{E_r} \int_O \pi_r(\dd x) G_\mathbb{D}^\varepsilon(x,x) \Big] \label{eq:decomp_renormalized} \\
    &\quad\quad\quad -\mathbb{E}\left[ \mathbf{1}_{E_r^c} \Big(  \langle \widetilde{\mathcal{L}}^\varepsilon(\mathcal{L}_\mathbb{D}^\theta), \pi_r  \rangle - \theta \int_\mathbb{D} \pi_r(\dd x) G_\mathbb{D}^\varepsilon(x,x)   \Big) \right]. \notag
\end{align}
Following Le Jan's framework \cite{jan2010markovpathsloopsfields, jan2008markovloopsrenormalization}, we define the renormalized truncated occupation time field $:\mathcal{L}^\varepsilon:$ of the full loop soup $\mathcal{L}_\mathbb{D}^\theta$ as,
$$
:\mathcal{L}^\varepsilon:(f) \,:=\, \sum_{\wp \in \mathcal{L}_\mathbb{D}^\theta} \mathbf{1}_{T(\wp) \geq \varepsilon} \int_0^{T(\wp)} f(\wp_s) \dd s - \theta \mu_\mathbb{D}^{\mathrm{loop}} \Big(  \mathbf{1}_{T(\wp) \geq \varepsilon} \int_0^{T(\wp)} f(\wp_s) \dd s \Big).
$$
Equation~\eqref{eq:decomp_renormalized} can be re-written more compactly as 
\begin{equation} \label{eq:decomp_D_local}
     \mathbb{E} \left[  \langle \widetilde{\mathcal{L}}^\varepsilon(\mathcal{L}_\mathbb{D}^\theta), \pi_r  \rangle \mathbf{1}_{E_r} \right] \,=\, \theta \mathbb{E} \Big[\mathbf{1}_{E_r} \int_O \pi_r(\dd x) G_\mathbb{D}^\varepsilon(x,x)  \Big] - \mathbb{E}\left[  :\mathcal{L}^\varepsilon:(\pi_r) \mathbf{1}_{E_r^c} \right].
\end{equation}
Subtracting \eqref{eq:decomp_O_local} from \eqref{eq:decomp_D_local} we get the exact identity, valid for all $\varepsilon > 0$ and $r > 0$,
\begin{equation} \label{eq:decomp_T_eps_Er}
\mathbb{E} \left[ \mathcal{T}^\varepsilon(\pi_r) \mathbf{1}_{E_r} \right] \,=\, \theta \mathbb{E} \Big[  \mathbf{1}_{E_r} \int_O  \big(G_\mathbb{D}^{\varepsilon}(x,x)- G_O^{\varepsilon}(x,x)\big) \pi_r(\dd x) \Big]  - \mathbb{E}\left[  :\mathcal{L}^\varepsilon:(\pi_r) \mathbf{1}_{E_r^c} \right].    
\end{equation}
We start with the first term. Recall the following fact \cite{jan2010markovpathsloopsfields,aru_lupu_sepulveda_ED_CR_CLE}: for all simply connected domains $D \subset \mathbb{D}$, and for all $x \in D$, the following limit is monotone
$$
\lim_{\varepsilon \to 0} \big( G_\mathbb{D}^\varepsilon(x,x) - G_D^\varepsilon(x,x) \big) \,=\, \frac{1}{\pi}\left(\log \CR(x, \mathbb{D}) - \log \CR(x, D)\right).
$$
By the monotone convergence theorem we get,
\begin{equation*}
\theta \lim_{\varepsilon \to 0} \mathbb{E} \Big[\mathbf{1}_{E_r} \int_O  \big(G_\mathbb{D}^\varepsilon(x,x) - G_O^\varepsilon(x,x) \big) \pi_r(\dd x) \Big] \,=\, \frac{\theta}{\pi} \mathbb{E} \bigg[ \mathbf{1}_{E_r} \int_O \pi_r(\dd x) \log { \frac{\CR(x, \mathbb{D})}{\CR(x, O)}} \bigg].    
\end{equation*}
We have seen that $\mathbf{1}_{E_r} \to 1$ almost surely, and since $x \mapsto \log \CR(x,O)$ is almost surely continuous at $0$ with $\CR(0,\mathbb{D}) = 1$, the integral $\int_O \pi_r(\dd x)\big(\log \CR(x,\mathbb{D}) - \log \CR(x,O)\big)$ converges almost surely to $-\log \CR(0,O)$. Moreover, by monotonicity of the conformal radius and Koebe's quarter theorem, on $E_r$ this integral is bounded uniformly in $r \in (0,1)$ by $-\log \CR(0,O) + C$ for some constant $C > 0$. Since $-\log \CR(0,O)$ is integrable by \cite{Schramm_2009}, Lebesgue's dominated convergence theorem yields
\begin{equation} \label{eq:log_CR_main_term_local}
    \theta \lim_{r \to 0}\lim_{\varepsilon \to 0} \mathbb{E} \Big[\mathbf{1}_{E_r} \int_O \big(G_\mathbb{D}^\varepsilon(x,x) - G_O^\varepsilon(x,x) \big) \pi_r(\dd x) \Big] \,=\, -\frac{\theta}{\pi} \mathbb{E} \left[ \log \CR(0,O) \right]
\end{equation}

Now, by \cite{jan2010markovpathsloopsfields}, for a deterministic test function $\pi_r$, the truncated renormalized occupation time field $:\mathcal{L}^\varepsilon:(\pi_r)$ integrated against $\pi_r$ converges in $L^2(\mathbb{P})$ to the renormalized occupation time field $:\mathcal{L}:(\pi_r)$ as $\varepsilon \to 0$. By continuity of the $L^2$ scalar product, we get
$$
\lim_{\varepsilon \to 0} \left|\mathbb{E}\left[  :\mathcal{L}^\varepsilon:(\pi_r) \mathbf{1}_{E_r^c} \right] \right| \,=\, \left| \mathbb{E}\left[  :\mathcal{L}:(\pi_r) \mathbf{1}_{E_r^c} \right] \right| \,\leq\, \mathbb{E}\big[  :\mathcal{L}:(\pi_r)^2 \big]^{1/2} \mathbb{P}(E_r^c)^{1/2}.
$$
In the last bound, we have used Cauchy-Schwarz's inequality. Le Jan proved in \cite{jan2010markovpathsloopsfields} that the covariance function of $:\mathcal{L}:$ is given by $(x,y) \mapsto \theta G_\mathbb{D}(x,y)^2$. Since $G_\mathbb{D}(x,y) \asymp |\log |x-y||$ as $x \to y$ and $\pi_r(\dd x)$ is the uniform measure on $r \mathbb{D}$, we get
$$
\mathbb{E}\big[ \big(:\mathcal{L}:(\pi_r)) ^2 \big] \,=\, \theta \iint_{\mathbb{D}^2}  G_\mathbb{D}(x,y)^2 \pi_r(\dd x) \pi_r(\dd y)  \,=\, O(|\log r|^2).
$$
By Koebe's quarter theorem, we have $E_r^c = \{ d(0, \partial O) \leq r\} \subset \{\CR(0,O) \leq 4r\}$ and by \cite{Schramm_2009}, $\log \CR(0, O)$ has exponential tails up to some sufficiently small order. Hence, there exists $\alpha > 0$ such that $\mathbb{P}(\CR(0, O) \leq 4r) = O(r^\alpha)$, which yields
\begin{equation} \label{eq:error_towards_local}
    \lim_{r \to 0} \lim_{\varepsilon \to 0} \mathbb{E}\left[  :\mathcal{L}^\varepsilon:(\pi_r) \mathbf{1}_{E_r^c} \right] \,=\, \lim_{r \to 0} O\big( |\log r| r^{\alpha/2} \big) \,=\, 0
\end{equation}
Injecting \eqref{eq:error_towards_local} and \eqref{eq:log_CR_main_term_local} into \eqref{eq:decomp_T_eps_Er} and \eqref{eq:montone_convergence_T}, yields the desired result.
\end{proof}

For a simply connected domain $D \subsetneq \mathbb{C}$ and $\beta >0$, we define $\mathcal{E}^{D,\beta}$ as a point process of first order intensity measure $\beta \mu_D^{\mathrm{exc}}$, where $\mu_D^{\mathrm{exc}}$ is the Brownian excursion measure defined in Section~\ref{sec:bubbles}. We do not suppose that $\mathcal{E}^{D,\beta}$ is Poissonian, and define the occupation time field $\Theta^{D,\beta}$ of $\mathcal{E}^{D,\beta}$ as 
\begin{equation} \label{eq:occupation_excur}
\Theta^{D,\beta}(f) \,=\, \sum_{e \in \mathcal{E}^{D,\beta}} \int_0^{T(e)} f(e(t)) \dd t,    
\end{equation}
where $T(e)$ denotes the duration of the excursion $e$. In order to prove Theorem~\ref{th:intensity}, we need the following lemma, which states that the expected occupation time density of $\mathcal{E}^{D,\beta}$ only depends on the first order intensity measure, and is independent of the domain $D$.

\begin{lemma} \label{lemma:density_occupation}For all simply connected domains $D \subsetneq \mathbb{C}$, we have
    $$
    \mathbb{E} \big[ \Theta^{D,\beta}(\dd x) \big] \,=\, 2 \beta \mathbf{1}_{x \in D} \dd x.
    $$
\end{lemma}

\begin{proof} We recall that in our normalisation, the generator of Brownian motion is $\frac{1}{2} \Delta$, the excursion measure is defined on ordered pairs and that the Poisson kernels are given in Section~\ref{sec:bubbles}. We first work in the case $D = \mathbb{D}$. As noted in \cite{jego2026heightgapplanarbrownian}, a small calculation gives that for any distinct boundary points $a,b \in \partial \mathbb{D}$, the expectation of the occupation measure of a single Brownian excursion from $a$ to $b$ in $\mathbb{D}$ has density w.r.t. Lebesgue on $\mathbb{D}$ given by
    $$
    x\in \mathbb{D} \,\longmapsto\, 2 \frac{H_\mathbb{D}(x,a)H_\mathbb{D}(x,b)}{H_\mathbb{D}(a,b)}.
    $$
(We believe the factor $\frac{1}{\pi}$ in \cite[Lemma 4.5]{jego2026heightgapplanarbrownian} is a typo, that does not affect their results, as this density is not used to compute the height gap of planar Brownian motion). Since, by definition, $\mu_\mathbb{D}^{\mathrm{exc}} = \iint_{(\partial \mathbb{D})^2} H_\mathbb{D}(a,b) \mathbb{P}^{a,b}_\mathbb{D} \dd a \dd b$, we get,
$$
\mathbb{E} \big[ \Theta^{\mathbb{D}, \beta}(f) \big] \,=\, 2\beta \int_\mathbb{D} \dd x f(x) \Big( \int_{\partial \mathbb{D}} H_\mathbb{D}(x,a) \dd a\Big) \Big( \int_{\partial \mathbb{D}} H_\mathbb{D}(x,b) \dd b \Big) \,=\, 2\beta \int_\mathbb{D} f(x) \dd x.
$$
This yields the desired result in $\mathbb{D}$. Let $D \subsetneq \mathbb{C}$ be a simply connected domain and $\varphi : \mathbb{D} \to D$ be a conformal map. By definition, $\mu_D^{\mathrm{exc}}$ is the image of $\mu_\mathbb{D}^{\mathrm{exc}}$ by $\varphi$ with the Brownian time change $\dd s = |\varphi'(e_t)|^2 \dd t $ (it does not depend on the choice of $\varphi$ by conformal invariance). We obtain
\begin{align*}
\mathbb{E} \big[ \Theta^{D, \beta}(f) \big] \,&=\, \beta \int  \mu_\mathbb{D}^{\mathrm{exc}}(\dd e) \Big( \int_0^{T(e)}  f\circ \varphi(e_t) |\varphi'(e_t)|^2 \dd t \Big) \\
&=\, 2\beta \int_\mathbb{D} f \circ \varphi(x) |\varphi'(x)|^2 \dd x \,=\, 2\beta \int_D f(y) \dd y.
\end{align*}
The proof is complete.
\end{proof}

\noindent \textbf{Proof of Theorem~\ref{th:intensity}. } Instead of working with the intensity parameter $\theta$, we choose here to work with $\kappa \in (8/3,4]$ and denote $\theta := \theta(\kappa)$ the associated loop soup intensity. We keep denoting $O$ the (interior of) the outermost cluster surrounding the origin in $\mathcal{L}_\mathbb{D}^{\theta}$. Recall (see Section~\ref{sec:excursions_boundary}) that $\mathcal{E}_{\partial O}$ denotes the set of excursions from the boundary of $O$. Mapping to the unit disc with any $\sigma(O)$-measurable conformal map $\psi$ from $O$ onto $\mathbb{D}$ yields $\mathcal{E}_{\partial \mathbb{D}}^\kappa = \psi(\mathcal{E}_{\partial O})$.

As pointed out in \cite[page 3]{qian2018lawpointprocessbrownian}, the first order intensity measure on the ordered pairs $(x,y) \in (\partial \mathbb{D} \times \partial \mathbb{D})$ of endpoints of the excursions, has to be invariant under Möbius transformations of the disc, forcing it to be a constant multiple of $|x-y|^{-2} \dd x \dd y$ (since the authors also proved that the excursions are locally finite). Recalling that the boundary Poisson kernel in the unit disc is given by $H_\mathbb{D}(x,y) = \frac{1}{\pi} |x-y|^{-2} \dd x \dd y$, this yields the existence of some $\beta(\kappa) > 0$ such that the first order intensity measure of the endpoints equals
\begin{equation} \label{eq:lambda_kappa}
\beta(\kappa) \, H_{\mathbb{D}}(x,y) \dd x\, \dd y.    
\end{equation}
Recalling (see Section~\ref{sec:bubbles}) that the Brownian excursion measure in $\mathbb{D}$ is defined as $\mu_\mathbb{D}^{\mathrm{exc}} = \iint_{\partial \mathbb{D}^2} \dd x \dd y \, H_\mathbb{D}(x,y) \, \mathbb{P}^{x,y}_{\mathbb{D}}$ and that  conditionally on the set of endpoints, the excursions are distributed as independent Brownian excursions linking $x(e)$ to $y(e)$  \cite[Lemma 9]{qian2018lawpointprocessbrownian}, we get 
$$
\mu_{\partial \mathbb{D}}^{\mathcal{E}, \kappa} \,=\, \beta(\kappa) \, \mu_\mathbb{D}^{\mathrm{exc}}.
$$
By conformal invariance in law of $\mathcal{E}_{\partial \mathbb{D}}^\kappa$, this also implies that, conditionally on $\partial O$, $\mathcal{E}_{\partial O}$ is almost surely a point process of first-order intensity measure $\beta(\kappa)\mu^{\mathrm{exc}}_{O} = \beta(\kappa) \psi^{-1} \circ \mu_\mathbb{D}^{\mathrm{exc}}$ for any conformal map $\psi^{-1} : \mathbb{D} \to O$. It remains to identify the intensity $\beta(\kappa)$. 

Denote $E_r = \{d(0, \partial O) > r\}$ and $\pi_r$ the uniform measure on $r\mathbb{D}$. Recall the definition \eqref{eq:def_occupation_O} of the occupation time field $\mathcal{T}$ of the loops that are not contained in $O$. On the event $E_r$, we have $\mathrm{supp}(\pi_r) = \overline{r\mathbb{D}} \subset O$. By connectedness, this implies that any loop $\wp \in \mathcal{L}_{\mathbb{D}  \setminus O}$ intersecting $\mathrm{supp}(\pi_r)$ must also intersect $\partial O$ and its contribution to $\mathcal{T}(\pi_r)$ in \eqref{eq:def_occupation_O} only comes from the times $s$ at which $\wp_s \in O$. In other words, only the excursions from the boundary of $O$ matter. Hence, conditionally on $O$, almost surely on $E_r$,
\begin{equation*} 
\mathcal{T}(\pi_r) \,=\, \sum_{e \in \mathcal{E}_{\partial O}} \int_0^{T(e)} \pi_r(e_t) \dd t\,=\, \Theta^{O,\beta}(\pi_r).
\end{equation*}
Where $\Theta^{O,\beta}$ is the occupation time field \eqref{eq:occupation_excur} of the excursions $\mathcal{E}_{\partial O}$ (whose first order intensity measure conditionally on $O$ is given by $\beta(\kappa) \mu_O^{\mathrm{exc}}$). Consequently, Lemma~\ref{lemma:density_occupation} applies conditionally almost surely on $E_r$ and yields
$$
\mathbf{1}_{E_r} \mathbb{E} \left[ \Theta^{O,\beta}(\pi_r) \,\mid\, O \right] \,=\, 2\beta(\kappa) \mathbf{1}_{E_r} \int_O \pi_r(\dd x) \,=\, 2\beta(\kappa)\mathbf{1}_{E_r},
$$
because $\mathrm{supp} (\pi_r) \subset O$ on the event $E_r$. We get
$$
\lim_{r \to 0}\mathbb{E} \left[ \mathcal{T}(\pi_r) \mathbf{1}_{E_r} \right] \,=\, 2\beta(\kappa) \lim_{r \to 0} \mathbb{P}(E_r) \,=\, 2\beta(\kappa).
$$
But, by Lemma~\ref{lemma:lim_T_lejan}, this limit also equals $\frac{\theta(\kappa)}{\pi} \mu(\kappa)$ where $\mu(\kappa)$ is the expected $|\log|$-conformal radius of a CLE$_\kappa$ loop surrounding the origin. Hence,
$$
2\pi \beta(\kappa) \,=\, \theta(\kappa) \mu(\kappa).
$$
It only remains to prove that $\beta(\kappa)$ has the explicit expression from Equation~\eqref{eq:g_kappa_explicit} and study the function $\beta(\kappa)$ on $(8/3,4)$. By \cite[Equation (2)]{Schramm_2009}, we have $\mu(\kappa) = (\pi \tan(4\pi/\kappa))/(1-\kappa/4)$ and by \cite{sheffield2011conformalloopensemblesmarkovian}, we also have $\theta(\kappa)=(3 \kappa-8)(6-\kappa) /(4 \kappa)$. Combining these equalities yields:
\begin{equation} \label{eq:def_beta_kappa}
\beta(\kappa)  \,=\, \frac{(3 \kappa-8)(6-\kappa)}{2 \kappa(4-\kappa)} \tan \Big(\frac{4 \pi}{\kappa}\Big).        
\end{equation}
As $\kappa \to 8/3^+$, using $\tan(4\pi/\kappa) = \cot\left(\pi(3\kappa-8)/(2\kappa)\right) = \frac{2\kappa}{\pi(3\kappa-8)} + O(3\kappa-8)$ gives
\[
\lim_{\kappa \to 8/3^+} \beta(\kappa) = \lim_{\kappa \to 8/3^+} \frac{6-\kappa}{\pi(4-\kappa)} = \frac{5}{2\pi}.
\]
As $\kappa \to 4^-$, applying the asymptotic relation $\tan(4\pi/\kappa) = \frac{\pi(4-\kappa)}{\kappa} + O((4-\kappa)^3)$ yields
\[
\lim_{\kappa \to 4^-} \beta(\kappa) = \lim_{\kappa \to 4^-} \frac{\pi(3\kappa-8)(6-\kappa)}{2\kappa^2} = \frac{\pi}{4}.
\]
This continuously extends $\beta$ to $[8/3, 4]$. 
Logarithmic differentiation of $\beta(\kappa)$ on $(8/3, 4)$ yields
\begin{equation} \label{eq:kappa_zero}
\frac{\beta'(\kappa)}{\beta(\kappa)} = \frac{3}{3\kappa-8} - \frac{1}{6-\kappa} - \frac{1}{\kappa} + \frac{1}{4-\kappa} - \frac{8\pi}{\kappa^2 \sin\left(\frac{8\pi}{\kappa}\right)},    
\end{equation}
A simple study of the function in \eqref{eq:def_beta_kappa} proves the non-monotonicity of $\beta(\kappa)$ and the fact that $\kappa_0$ is the unique root of Equation~\eqref{eq:kappa_zero}.
$\hfill \square$

\begin{remarque} \label{rem:brownian_area}
Let us conclude with a heuristic observation in the limit $\kappa \to 8/3$, which corresponds to $\theta \to 0$ and therefore to a single Brownian loop. Let $(\wp_t)_{t \in [0,1]}$ be a Brownian loop of duration $1$, with outer boundary $\partial \wp$ and interior domain $O = \mathrm{Int}(\partial \wp)$.

By~\cite{Taylor1964TheEH}, the occupation measure of $\wp$ is a constant multiple of the Hausdorff measure in the gauge $r^2 \log \frac{1}{r} \log \log \log \frac{1}{r}$. Since $\partial \wp$ has Hausdorff dimension $\frac{4}{3} < 2$, the loop $\wp$ spends almost surely zero time on $\partial \wp$. Decomposing $\wp$ into its excursions $\mathcal{E}_{\partial O}$ inside $O$ from $\partial O$, we obtain 
$$1 = T(\wp) = \sum_{e \in \mathcal{E}_{\partial O}} T(e) = \Theta^{O, \beta(8/3)}(\mathbf{1}),$$
where $\Theta^{O, \beta(8/3)}$ is the occupation time field of $\mathcal{E}_{\partial O}$. Conditionally on $O$, our heuristic suggests that the first order intensity measure of $\mathcal{E}_{\partial O}$ is $\beta(8/3) \mu_O^{\mathrm{exc}}$. Applying Lemma~\ref{lemma:density_occupation} conditionally as in the proof of Theorem~\ref{th:intensity} yields $
1 = 2\beta(8/3)\, \mathbb{E}[\mathrm{Area}(\mathrm{Int}(\wp))].$ Since $\beta(8/3) = 5/(2\pi)$, we recover that the expected area enclosed by a Brownian loop of duration $1$ is $\pi/5$, which is the result of Garban and Ferreras Trujillo~\cite{Garban_2006}. Moreover, as noted by Jego, Lupu and Qian~\cite{jego2026heightgapplanarbrownian}, the height gap of planar Brownian motion equals the occupation time density inside $O$, which gives $2\beta(8/3) = 5/\pi$.
\end{remarque}

\appendix
\section{Tool and computational resource disclosure} \label{sec:ai}

In line with the recommendations of the Leiden Declaration on AI and Mathematics, the author discloses the use of automated tools in this work. The model \textbf{Gemini 3.8 Flash} (free web version, no paid subscription) was used as a research assistant for the following tasks:
\begin{itemize}
    \item \textbf{Technical proofs and calculations:} The model helped compute the exact formulas for $K, J, a,$ and $\rho$ in Lemma~\ref{lemma:conformal_circles}, the integrals at the end of Lemma~\ref{lem:ppp_naive_projection} and Appendix~\ref{sec:complex_analysis}, and compute the asymptotic expansion of $|\log \rho(r,y)|$ in Lemma~\ref{lemma:mu_bub_hits_log}. In Lemma~\ref{lemma:hits_D_y_r}, the model also simplified heavy calculations by identifying the Poisson kernel in the expression of $|g'(e^{i\theta})|$. In Theorem~\ref{th:crossing_exponent}, the idea of using the widened random covering $\Lambda^{(K)}$ instead of a direct computation of the conditional distribution of $N_K(y,x)$ given $N(y,x)$ was suggested by the model. It also helped to study the limits of $\beta(\kappa)$ defined in \eqref{eq:def_beta_kappa} when $\kappa \to 8/3$ and $\kappa \to 4$. It also identified the flaw $\frac{1}{\pi}$ in \cite[Lemma 4.5]{jego2026heightgapplanarbrownian} coming from the $\frac{1}{2} \Delta$ normalisation of Brownian motion. 
    \item \textbf{Literature:} The model was used to find references on one-dimensional random coverings, on the Fourier decomposition of the heat kernel on the circle and on the estimate for SLE$_\kappa$ hitting small balls by Beffara. All references were subsequently checked by the author.
    \item \textbf{Language and spelling:} The model helped translate and improve sentences from French to mathematical English. In the early stages of preparing the paper, the model was used to correct spelling mistakes.
\end{itemize}

The main conceptual ideas were conceived entirely by the author. These include the 1D projection reduction, the core intervals, the gap-filling strategy with small loops, the asymptotic density of the separating cluster, the argument showing that non-separating loops cannot hide in unvisited holes, and the occupation time argument for Theorem~\ref{th:intensity}. All calculations, arguments, and references suggested by the tool were checked by the author. The author takes full responsibility for the correctness of the paper.

\section{Proof of Corollaries~\ref{cor:exact_covariance}, \ref{cor:Minkowski} and \ref{cor:polar_sets}} \label{sec:proof_cors}

This section is dedicated to the proofs of the further results, that are all consequences of the results of \cite{Jego_2025} and \cite{jego2023conformallyinvariantfieldsbrownian}. Following their framework, for $\varepsilon > 0$, we set
$$Z_\varepsilon \,:=\, \mathbb{P}\Big(e^{-1} \partial \mathbb{D} \stackrel{\mathcal{L}_{\mathbb{D}}^\theta}{\longleftrightarrow} \varepsilon \partial \mathbb{D}\Big) \,\asymp\, |\log \varepsilon|^{-1+\theta}.$$
Our proofs closely follow the arguments of \cite{Jego_2025} and \cite{jego2023conformallyinvariantfieldsbrownian}, so we omit most details. We only detail the steps that change here, and we strongly encourage the reader to consult those papers for the full, rigorous, proofs. Proposition 4.1 in \cite{Jego_2025} is a crucial component of most of the Corollaries, that can be strengthened using our sharper crossing exponent.

\begin{proposition}[Sharpened Proposition 4.1 in \cite{Jego_2025}] \label{prop:JLQ_prop41} Let $r_0>0$ be a fixed macroscopic radius. There exists $C>0$ such that for all $x, y \in \mathbb{D}$ with $|x-y| \leq 1/2$ and, $\varepsilon \in (0,e^{-1})$,
\begin{equation} \label{eq:prop41_1}
\frac{1}{Z_{\varepsilon}^2} \mathbb{P}\left(D(x, \varepsilon) \stackrel{\mathcal{L}_{\mathbb{D}}^\theta}{\longleftrightarrow} D(y, \varepsilon)\right) \leq C|\log | x-y ||^{2(1-\theta)},
\end{equation}
\begin{equation} \label{eq:prop41_2}
\frac{1}{Z_{\varepsilon}^2} \mathbb{P}\left(D(x, \varepsilon) \stackrel{\mathcal{L}_{\mathbb{D}}^\theta}{\longleftrightarrow} D(y, \varepsilon), D(x, \varepsilon) \stackrel{\mathcal{L}_{\mathbb{D}}^\theta}{\longleftrightarrow} \partial D\left(x, r_0\right)\right) \leq C|\log | x-y ||^{1-\theta}.
\end{equation}
\end{proposition}

\begin{proof}
    As explained by the authors in \cite{Jego_2025} (Proposition 4.1), \eqref{eq:prop41_1} is a direct consequence of \eqref{eq:prop41_2}, we therefore only focus on the second part of the statement. In their proof, the authors used threshold $|x-y| \leq \exp(-|\log \varepsilon|^{1-\delta})$ for some $\delta > 0$. In order to remove the error $\eta$ in the exponents that appeared in the authors' proof, assume first that $|x-y| \leq 2\varepsilon^{1/2}$. In that case, we have $|\log \varepsilon| \leq C_1|\log |x-y||$, which gives directly that 
$$
\frac{1}{Z_{\varepsilon}^2} \mathbb{P}\left(D(x, \varepsilon) \stackrel{\mathcal{L}_{\mathbb{D}}^\theta}{\longleftrightarrow} \partial D\left(x, r_0\right)\right) \leq C_2 \frac{ |\log \varepsilon|^{-1+\theta}}{Z_\varepsilon^2} \,\leq\,  C_3|\log |x-y||^{1-\theta}.
$$
We now assume that $|x-y| > 2\varepsilon^{1/2}$. Using the notations of \cite[Proof of Proposition 4.1]{Jego_2025}, their proof yields that
\begin{equation} \label{eq:bound_jlq_r_z_R}
    \eqref{eq:prop41_2} \leq \frac{1}{Z_\varepsilon^2} \mathbb{E} \bigg[ \prod_{z = x,y} \mathbb{P}\Big({\textstyle\frac{2 \varepsilon}{|x-y|}} \mathbb{D} \stackrel{\mathcal{L}_{\mathbb{D}}^\theta}{\longleftrightarrow} {\textstyle\frac{2 r_z}{|x-y|}} \partial \mathbb{D} \,\mid\, r_z\Big)  \mathbb{P}\Big({\textstyle\frac{R}{2}} \partial \mathbb{D} \stackrel{\mathcal{L}_{\mathbb{D}}^\theta}{\longleftrightarrow} {\textstyle\frac{r_0}{2}} \partial \mathbb{D} \,\mid\, R\Big) \bigg]
\end{equation}
By Theorem~\ref{th:crossing_exponent}, the right-most probability conditionally on $R$ is bounded by a constant multiple (depending only on $r_0$) of $|\log R|^{-1+\theta}$. In addition, in the main regime $\varepsilon < r_z \leq e^{-1} |x-y|/2$ the authors have proved that ($1/Z_\varepsilon$ times) the probability conditionally on $r_z$ is bounded by a constant multiple of 
\begin{equation} \label{eq:fkg_r_z}
\mathbb{P}\left(\varepsilon \mathbb{D} \stackrel{\mathcal{L}_{\mathbb{D}}^\theta}{\longleftrightarrow} {\textstyle\frac{2 \varepsilon}{|x-y|}} \partial \mathbb{D}\right)^{-1} \mathbb{P}\left({\textstyle\frac{2 r_z}{|x-y|}} \mathbb{D} \stackrel{\mathcal{L}_{\mathbb{D}}^\theta}{\longleftrightarrow} e^{-1} \partial \mathbb{D}\right)^{-1}.    
\end{equation}
Setting $u = 2\varepsilon / |x-y|$, since $|x-y| > 2\varepsilon^{1/2}$, we have $u^2 < \varepsilon$. The annulus $u\mathbb{D} \setminus \varepsilon \mathbb{D}$ is contained in $u\mathbb{D} \setminus u^2 \mathbb{D}$. Hence, by \cite[Theorem 1.3]{Jego_2025}
$$
\mathbb{P}\left(\varepsilon \mathbb{D} \stackrel{\mathcal{L}_{\mathbb{D}}^\theta}{\longleftrightarrow} {\textstyle\frac{2 \varepsilon}{|x-y|}} \partial \mathbb{D}\right) \,\geq\, \mathbb{P}\left(u^2 \mathbb{D} \stackrel{\mathcal{L}_{\mathbb{D}}^\theta}{\longleftrightarrow} u \partial \mathbb{D}\right)  \xrightarrow[u \rightarrow 0]{ } f_{\infty}(2)>0.
$$
By our Theorem~\ref{th:crossing_exponent}, the second term in \eqref{eq:fkg_r_z} is bounded by a constant multiple of $\big(\log \frac{|x-y|}{r_z})^{1-\theta}$. The remaining boundary regimes are treated exactly as the authors, with our sharper crossing exponent. Injecting these facts into Equation~\eqref{eq:bound_jlq_r_z_R}, we get,
\begin{equation} \label{eq:bound_esperance_decouples}
\eqref{eq:prop41_2} \,\leq\, C\,\mathbb{E} \left[ \prod_{z=x,y} \Big( \log \frac{|x-y|}{r_z} \Big)^{1-\theta} |\log R|^{-1+\theta} \right].    
\end{equation}
This is precisely \cite[Equation (4.11)]{Jego_2025}, without any $\eta$-error in the exponents. From there, the rest of the argument in \cite[pp.~30--32]{Jego_2025} (discrete integrations by parts and Claim~4.4) still holds with $\eta = 0$. We omit the details and refer to the original proof by the authors.
\end{proof}

\noindent \textbf{Proof of Corollary~\ref{cor:exact_covariance}.} We work here in the unit disc. By Equation (4.6) from \cite[Proposition 4.3]{jego2023conformallyinvariantfieldsbrownian} applied to $F = 1$, together with Equation (4.3),  we get
$$
C_\theta(x,y) \,=\, \lim _{r \rightarrow 0} \frac{1}{Z_r^2} \mathbb{P}\Big(D(x, r) \stackrel{\mathcal{L}_\mathbb{D}^\theta}{\longleftrightarrow} D(y, r)\Big) \,\leq\, C|\log |x-y||^{2(1-\theta)}.
$$
We have applied Proposition~\ref{prop:JLQ_prop41}. Let $\delta > 0$ and $x,y \in \mathbb{D}$ at distance at least $\delta$ from $\partial \mathbb{D}$. By the FKG inequality \cite[Lemma 2.1]{JANSON1984313} together with standard welding arguments, we get
$$
\mathbb{P}\Big(D(x, r) \stackrel{\mathcal{L}_\mathbb{D}^\theta}{\longleftrightarrow} D(y, r)\Big) \,\geq\, p \prod_{z = x,y} \mathbb{P}\Big(D(z, r) \stackrel{\mathcal{L}_\mathbb{D}^\theta}{\longleftrightarrow} \partial D(z, {\textstyle\frac{|x-y|}{2}})\Big),$$
for some positive constant $p \in (0,1)$. Since $x$ is at distance at least $\delta$ from the boundary, we have $D(x, \delta) \subset \mathbb{D}$. By scaling and translation, we get
\begin{align*}
\mathbb{P}\Big(D(z, r) \stackrel{\mathcal{L}_\mathbb{D}^\theta}{\longleftrightarrow} \partial D(z, {\textstyle\frac{|x-y|}{2}})\Big) \,&\geq\, \mathbb{P}\left(\frac{r}{\delta} \partial \mathbb{D} \stackrel{\mathcal{L}_{\mathbb{D}}^\theta}{\longleftrightarrow} \frac{|x-y|}{2 \delta} \partial \mathbb{D}\right) \\
&\geq\, \mathbb{P}\big([u,v] \text{ is covered by } \underline{\Lambda}^{\theta,1}\big).
\end{align*}
Since $|x-y| \to 0$, we have assumed $|x-y| \le \delta$ and we have set $u := |\log \frac{|x-y|}{2\delta}| \ge \log 2$ and $v := |\log \frac{r}{\delta}|$. In the last line, we have applied Lemma~\ref{lem:lm_hypothesis_h} (with $\tau = 1$). Since $Z_r \asymp |\log r|^{-1+\theta} \sim v^{-1+\theta}$ as $r \to 0$ (for fixed $\delta > 0$), Proposition~\ref{prop:random_coverings} yields
\[
\liminf_{r \to 0} \frac{1}{Z_r} \mathbb{P}\Big(D(z,r) \stackrel{\mathcal{L}_{\mathbb{D}}^\theta}{\longleftrightarrow}\partial D(z, {\textstyle\frac{|x-y|}{2}})\Big)
\ge c(\theta) K(u) \ge c'(\theta) u^{1-\theta} = c'(\theta) \Big|\log \frac{|x-y|}{2\delta}\Big|^{1-\theta},
\]
where the last inequality uses $K(u) \asymp u^{1-\theta}$ uniformly for $u \ge \log 2$ (Proposition~\ref{prop:random_coverings}). Combining the bounds for $z = x, y$ yields $C_\theta(x,y) \ge c''(\theta) |\log |x-y||^{2(1-\theta)}$ as $|x-y| \to 0$.

Since the Green's function characterizes the conformal type of any simply connected domain punctured by two points $\{x,y\}$, there exists a measurable function $F : [0,\infty) \to [0, \infty)$ such that for all simply connected bounded domains $D$ and $x,y \in D$, $C_\theta^D(x,y) = F(G_D(x,y))$. Applying this to the unit disc, we get $F(g) \asymp g^{2(1-\theta)}$ as $g \to \infty$. 
$\hfill \square$
\newline

\noindent \textbf{Proof of Corollary~\ref{cor:Minkowski}.} In \cite[Theorem 1.9]{jego2023conformallyinvariantfieldsbrownian}, the authors proved that the gauge function $\phi(r)$ is given by $r \mapsto r^2/Z_r$. The result follows immediately from $Z_r \asymp |\log r|^{-1+\theta}$.
$\hfill \square$
\newline

\noindent \textbf{Proof of Corollary~\ref{cor:polar_sets}.} The proof is identical to the second case of \cite[Theorem 1.5]{Jego_2025}. Suppose $\operatorname{Cap}_{\text {log }^{1-\theta}}(A)>0$, then there exists a probability measure $\mu$ on $A$ such that $ \int|\log | x-y||^{1-\theta} \mu(\mathrm{d} x) \mu(\mathrm{d} y) < \infty $.
Let $\mathcal{C}$ be the outermost cluster surrounding some given point of $D \backslash A$, and define
$$
J_{\varepsilon}:=\frac{1}{Z_{\varepsilon}} \int_A \mathbf{1}_{\{\mathcal{C} \cap D(x, \varepsilon) \neq \varnothing\}} \mu(\mathrm{d} x)
$$
Applying Paley-Zygmund inequality and doing the same steps as the authors it suffices to prove that $\mathbb{E}[J_\varepsilon^2]$ is uniformly bounded in $\varepsilon$. By the second inequality of our sharper Proposition~\ref{prop:JLQ_prop41} (recall $\mathcal{C}$ is an outermost cluster),
\begin{align*}
\mathbb{E} \big[ J_\varepsilon^2 \big] \,&=\, \iint_{A \times A} \mu(\dd x) \mu(\dd y) \frac{1}{Z_\varepsilon^2} \mathbb{P}(\mathcal{C} \cap D(x,\varepsilon) \neq \emptyset, \mathcal{C} \cap D(y,\varepsilon) \neq \emptyset) \\    
&\leq\, C\iint_{A \times A} \mu(\dd x ) \mu( \dd y) |\log |x-y||^{1-\theta}.
\end{align*}
Hence, $A$ is not polar for the clusters of $\mathcal{L}_D^\theta$.
$\hfill \square$ \newline

Finally, let us explain why the converse implication cannot be deduced from~\cite{Jego_2025}. Repeating their proof (just above~(4.19) in~\cite{Jego_2025}) with our sharper crossing exponent $\mu(D(y, e^{-k})) \leq C k^{-1+\theta}$ yields, for $\alpha = 1-\theta$,
\[
\int |\log |x-y||^{1-\theta} \mu(\mathrm{d}x) \leq 1 + C \sum_{k=1}^{\infty} \frac{1}{k}.
\]
This series diverges, even if one knew the exact asymptotic for the crossing probability. The issue is that bounding $\mu(D(y, e^{-k}))$ simply by the probability that $\mathcal{C}_i$ crosses $D(y, \eta) \setminus D(y, e^{-k})$ (which is of order $k^{-1+\theta}$) is too crude, as it ignores the fact that $X$ can lie elsewhere in $\overline{\mathcal{C}}_i \cap A$.

\section{Technical results.} \label{sec:appendix}

Here, we gather some technical lemmas that are used in the proofs of the main results.

\subsection{Estimates on Brownian excursions}

\begin{lemma} \label{lemma:excur_max_gaussian}
Let $(e_t)_{t \geq 0}$ be a Brownian excursion of duration $1$, then,  there exists a positive constant $C$, such that for all $x > 0$, $
\mathbb{P}\big(\sup_{t \geq 0} e_t \geq x\big) \,\leq\, C e^{-x^2}.$
\end{lemma}

\begin{proof}
By \cite[Theorem 1]{max_brownian}, for all $x \geq 1$,
$$
\mathbb{P}\big(\sup_{t \geq 0} e_t \geq x\big) 
\,\leq\, 8 x^2 \sum_{n \geq 1} n^2 e^{-2x^2 n^2} 
\,=\, 8 x^2 e^{-2x^2} \sum_{n \geq 1} n^2 e^{-2x^2(n^2-1)}.
$$
Since $x \geq 1$, we have $e^{-2x^2(n^2-1)} \leq e^{-2(n^2-1)}$, so the series is bounded by the finite constant $K := \sum_{n \geq 1} n^2 e^{-2(n^2-1)} < \infty$.  The map $x \mapsto 8K x^2 e^{-2x^2}$ is bounded on $[1, \infty)$ by some constant $C_1 > 0$. Hence, for all $x \geq 1$,
$$
\mathbb{P}\big(\sup_{t \geq 0} e_t \geq x\big) \,\leq\, 8K x^2 e^{-x^2} e^{-x^2}  \leq C_1 e^{- x^2}.
$$
For $0 < x < 1$, since $\mathbb{P}(\sup e_t \ge x) \leq 1$, the bound holds for all $x > 0$ by choosing a larger constant $C$.
\end{proof}
\medskip

\begin{lemma} \label{lemma:ito_max_tau}
Let $x\geq 1, \tau > 0$ be fixed, and denote $\mathbf{n}$ Itô's excursion measure. Then, for some positive constants $C$ and $c$, depending only on $\tau$, we have
$$
\mathbf{n}\big(\sup _{t \geq 0} X_t \geq x \text { and } \zeta<\tau\big) \,\leq\,  Ce^{-c x^2}.
$$
Where $\zeta$ is the lifetime of $(X_t)_{t \geq 0}$.
\end{lemma}

\begin{proof}
By disintegrating Itô's excursion measure with respect to the duration $\zeta$, we have $\mathbf{n}(\zeta \in ds) = (2\pi s^3)^{-1/2} \dd s$, and conditionally on $\zeta = s$, $(X_t)_{t \geq 0} =  (\sqrt{s} e_{t/s})_{t \geq 0}$, in law, where $e$ is a standard Brownian excursion of duration $1$. Hence,
$$
\mathbf{n}\big(\sup_{t \geq 0} X_t \geq x \text{ and } \zeta < \tau\big) 
\,=\, \int_0^\tau \mathbb{P}\Big(\sup_{t \geq 0} e_t \geq \frac{x}{\sqrt{s}}\Big) \frac{ds}{\sqrt{2\pi s^3}}.
$$
By Lemma~\ref{lemma:excur_max_gaussian}, there exists a positive constant $C_1 > 0$ such that $\mathbb{P}(\sup e_t \geq y) \leq C_1 e^{- y^2}$. Applying this with $y = x/\sqrt{s}$ and performing the change of variables $u = x^2/s$, we get
\begin{align*}
\int_0^\tau s^{-3/2} e^{-x^2/s} \dd s 
\,=\, \frac{1}{x} \int_{x^2/\tau}^\infty u^{-1/2} e^{-u} \dd u 
&\,\leq\, \frac{1}{x} \Big(\frac{x^2}{\tau}\Big)^{-1/2} \int_{x^2/\tau}^\infty e^{-u} \dd u \\
&\,=\, \frac{\sqrt{\tau}}{x^2} e^{-x^2/\tau} 
\,\leq\, \sqrt{\tau} e^{-\frac{1}{\tau} x^2},
\end{align*}
where we bounded $u^{-1/2} \leq \frac{\sqrt{\tau}}{x}$ on $[x^2/\tau, \infty)$, and used $x > 1$ for the last inequality.  Plugging this back into the first equality yields the result with $c = \frac{1}{\tau}$ and $C = \frac{C_1 \sqrt{\tau}}{\sqrt{2\pi}}$.
\end{proof}

\subsection{Heat kernel on the circle.}

\begin{lemma} \label{lemma:_heat_kernel_bound}
Let $\tau>0$, there exist positive constants $C_\tau, C_\tau^{\prime}$ and $c$ such that for all $t > 0$, 
$$\left|2 \pi p_{\mathbb{S}^1}(t)-1\right| \mathbf{1}_{t \geq \tau} \leq C_\tau e^{-c t} \quad \text{ and }\quad p_{\mathbb{S}^1}(t) \mathbf{1}_{t \leq \tau} \leq C_\tau^{\prime} t^{-1 / 2}$$
\end{lemma}

\begin{proof} By \cite[Chapter 4.4]{stein2003fourier}, we have the following Fourier decomposition for the heat kernel on the circle:
$$
p_{\mathbb{S}^1}(t)=\frac{1}{2 \pi} \sum_{n \in \mathbb{Z}} e^{-n^2 t/2} .
$$
Note that in \cite{stein2003fourier}, the authors considered the heat kernel on the circle $\mathbb{R} / \mathbb{Z}$ of perimeter $1$, this explains the mismatch between their terms $e^{-4\pi^2n^2 t}$ and our $e^{-n^2 t / 2}$. We get, for all $t \geq \tau$,
$$
\left|2 \pi p_{\mathbb{S}^1}(t)-1\right|=2 \sum_{n \geq 1} e^{-n^2 t/2} \leq 2 e^{-t/2}\Big(1+\sum_{n \geq 2} e^{-\left(n^2-1\right) \tau/2}\Big) .
$$
The sum in the last bound is a finite constant depending only on $\tau$.
For the second bound, for any $t > 0$, we have $p_{\mathbb{S}^1}(t) = \frac{1}{2 \pi} + \frac{1}{\pi} \sum_{n \geq 1} e^{-n^2 t / 2}$. Since $u \mapsto e^{-u^2 t / 2}$ is decreasing on $[0, \infty)$, an integral comparison and the change of variables $v = u \sqrt{t/2}$ yield
$$
\sum_{n \geq 1} e^{-n^2 t / 2} \leq \int_0^{\infty} e^{-u^2 t / 2} \, \mathrm{d}u = \sqrt{\frac{2}{t}} \int_0^{\infty} e^{-v^2} \, \mathrm{d}v = \sqrt{\frac{2}{t}} \frac{\sqrt{\pi}}{2} = \sqrt{\frac{\pi}{2t}}.
$$
Multiplying by $\frac{1}{\pi}$, and recalling that for all $t \leq \tau$, we have $\frac{1}{2\pi} \leq \frac{\sqrt{\tau}}{2\pi} t^{-1/2}$, we get
$$
p_{\mathbb{S}^1}(t) \leq \Big(\frac{1}{\sqrt{2\pi}} + \frac{\sqrt{\tau}}{2\pi}\Big) t^{-1/2} =: C'_\tau t^{-1/2}.
$$
\end{proof}

\subsection{Complex analysis results} \label{sec:complex_analysis}

\noindent \textbf{Proof of Lemma~\ref{lemma:discs_cylinder_disc}.} Let $z = (x,e^{i\theta}) \in \mathfrak{C}$ and $r < \min(x, \pi)$ (which ensures that $D_\mathfrak{C}(z,r)$ does not touch the boundary nor wraps around the cylinder). This implies that $D_\mathfrak{C}(z,r)$ is isometric to a standard Euclidean disc in the complex plane $D_\mathfrak{C}(z,r) \cong \{w \in \mathbb{C} : |w-z| < r\}$. The image of this disc by $f^{-1}(w) = e^{-w}$ is then given by
$$
f^{-1}\big(D_\mathfrak{C}(z,r) \big) \,=\, \{e^{-w} : |w-z|<r\} \,=\, e^{-z} E_r \,=\, f^{-1}(z) E_r.
$$
We have set $E_r = \{e^{-\xi} : |\xi| < r \}$. Now, the boundary of $E_r$ is exactly the curve parametrised by $\gamma(\theta) = e^{re^{i\theta}}$ for $\theta \in [0, 2\pi)$. The distance to $1$ is given by
$$\big|e^{r e^{i \theta}}-1\big|^2=e^{2 r \cos \theta}-2 e^{r \cos \theta} \cos (r \sin \theta)+1$$
This function is decreasing on $[0,\pi]$ and its minimum is at $\theta = \pi$ which gives $(1-e^{-r})^2$. Consequently, $E_r$ contains the Euclidean ball $B(1,1-e^{-r})$. The proof is complete.
$\hfill \square$

\medskip

\begin{lemma} \label{lemma:conformal_circles}
Let $y \in (0,1)$ and $r > 0$ such that $B(y,r) \subset \mathbb{D}$. Define $$K = \frac{1+y^2-r^2}{2y} \quad \text{ and } \quad J = \frac{1+r^2-y^2}{2r}$$
Define $a = K - \sqrt{K^2-1}$ and $\rho = J  - \sqrt{J^2 -1}$. Then, the map $g(z) = \frac{z-a}{1-az}$ is a conformal automorphism of the unit disc that sends $B(y,r)$ onto $B(0,\rho)$.
\end{lemma}

\begin{proof}
Since $B(y,r) \subset \mathbb{D}$, we have $y+r < 1$ and $y-r > -1$. This implies that $K - 1 > 0$ and $J - 1 > 0$. Furthermore, since $a = K - \sqrt{K^2-1}$ and $\rho = J  - \sqrt{J^2 -1}$, $a$ and $\rho$ are the unique roots in $(0,1)$ of the quadratic equations $a^2 - 2Ka + 1 = 0$ and $\rho^2 - 2J\rho + 1 = 0$. Since $a \in (0,1)$, $g$ is a conformal automorphism of $\mathbb{D}$.

Since $g$ is a Möbius transformation with real coefficients, $g$ preserves the real line and maps any circle orthogonal to $\mathbb{R}$ to a circle orthogonal to $\mathbb{R}$. The boundary $\partial B(y,r)$ meets $\mathbb{R}$ orthogonally at $x_1 = y-r$ and $x_2 = y+r$. Therefore, it suffices to show that $g(x_2) = \rho$ and $g(x_1) = -\rho$. The relation $g(x_2) + g(x_1) = 0$ is equivalent to 
\[
(x_1 + x_2) a^2 - 2(1 + x_1 x_2) a + (x_1 + x_2) = 0 \iff a^2 - 2Ka + 1 = 0,
\]
which yields the correct quadratic expression for $a$. Setting $\rho = g(x_2) = -g(x_1)$, we have $x_2 = \frac{\rho+a}{1+a\rho}$ and $x_1 = \frac{-\rho+a}{1-a\rho}$, hence
\[
2r = x_2 - x_1 = \frac{2\rho(1-a^2)}{1-a^2\rho^2} \quad \text{and} \quad 2y = x_2 + x_1 = \frac{2a(1-\rho^2)}{1-a^2\rho^2}.
\]
Substituting these expressions yields,
\[
J = \frac{1+r^2-y^2}{2r} = \frac{1+\rho^2}{2\rho},
\]
which implies that $\rho^2 - 2J\rho + 1 = 0$, and hence that $\rho = J - \sqrt{J^2-1}$. Consequently, $\partial B(y,r)$ is mapped onto $\partial B(0,\rho)$. Since $g$ is a homeomorphism of $\mathbb{D}$, it maps the interior to the interior, proving $g(B(y,r)) = B(0,\rho)$.
\end{proof}

\bibliographystyle{alpha}
\bibliography{bibliography}
\end{document}